\documentclass[12pt]{article}

    \usepackage{amsmath}
    \usepackage{amsfonts}
\usepackage{color}
\usepackage{amsmath}
\usepackage{amsfonts}
\usepackage{graphicx}
\usepackage{epsf}
\usepackage{mathrsfs}
\usepackage{float}

\usepackage{amsmath,amssymb,latexsym,psfrag}
\usepackage{amsfonts}
\usepackage{graphicx}
\usepackage{epsf,epsfig,color}

\newtheorem{theo}{Theorem}

\newtheorem{prop}{Proposition}
\newtheorem{lemm}{Lemma}
\newtheorem{defn}{Definition}

\def\N{\mathbb{N}}
\def\E{\mathbb{E}}

\def\0{{\bf 0}}
\def\Z{\mathbb{Z}}
\def\R{\mathbb{R}}

\def\cov{{\rm Cov}}

\def\B{{ B}}

\def\en2{{\epsilon_n^2}}
\renewcommand{\E}{\mathbb E \,}
\newcommand{\C}{{\cal C}}

\newcommand{\tod}{\stackrel{{\cal D}}{\longrightarrow}}

\newcommand{\eqco}{\setcounter{equation}{0}}
\newcommand{\thco}{\setcounter{theo}{0}}
\newcommand{\prco}{\setcounter{prop}{0}}
\newcommand{\laco}{\setcounter{lemm}{0}}
\newcommand{\coco}{\setcounter{coro}{0}}
\newcommand{\cjco}{\setcounter{conj}{0}}

\newcommand{\deco}{\setcounter{defn}{0}}
\newcommand{\allco}{\eqco  \thco \prco \laco \coco \cjco \deco}
\newcommand{\X}{{\cal X}}

\def\2e{{\epsilon^2_\la}}
\def\la{{\lambda}}

\def\v{{ w }}

\newcommand{\M}{{\cal M}}

\renewcommand{\P}{{{\cal P}}}

\newcommand{\A}{{\cal A}}
\newcommand{\Cov}{{\rm Cov}}

\newcommand{\Var}{{\rm Var}}
\newcommand{\var}{{\rm Var}}

\newcommand{\card}{{\rm card}}
\newcommand{\Vol}{{\rm Vol}}
\newcommand{\diam}{{\rm diam}}

\newcommand{\liml}{\lim_{\lambda \to \infty} }

\newcommand{\tX}{{\tilde{X}}}
\newcommand{\tY}{{\tilde{Y}}}

\def\bdm{\begin{displaymath}}
\newcommand{\edm}{\end{displaymath}}
\def\benu{\begin{enumerate}}
\def\eenu{\end{enumerate}}
\def\beqn{\begin{equation}}
\def\eeqn{\end{equation}}
\def\be{\begin{equation}}
\def\ee{\end{equation}}
\def\bea{\begin{eqnarray}}
\def\eea{\end{eqnarray}}
\newcommand{\bean}{\begin{eqnarray*}}
\newcommand{\eean}{\end{eqnarray*}}
\newcommand{\bear}{\begin{eqnarray}}
\newcommand{\eear}{\end{eqnarray}}

\def\R{\mathbb{R}}
\def\B^2{\mathbb{D}}

\def\B{\mathbb{B}}

\def\tx{\tilde{\xi} }
\def\hx{{\hat{\xi}}}

\def\de{{\delta}}
\def\A{{\cal A}}

\def\qed{\hfill\hbox{${\vcenter{\vbox{
    \hrule height 0.4pt\hbox{\vrule width 0.4pt height 6pt
    \kern5pt\vrule width 0.4pt}\hrule height 0.4pt}}}$}}

\def\la{{\lambda}}
\def\ka{{\kappa}}

\def\la{{\lambda}}
\def\ka{{\kappa}}

\newcommand{\blue}[1]{\textcolor{blue}{#1}}

\begin{document}
\title{\bf Variance Asymptotics and Scaling Limits for Random  Polytopes}

\author{Pierre Calka$^{*}$,  J. E. Yukich$^{**}$ }

\maketitle

\footnotetext{{\em American Mathematical Society 2010 subject
classifications.} Primary 60F05, 52A20; Secondary 60D05, 52A23} 

\footnotetext{$~^{*}$ Research partially supported by French ANR grant PRESAGE.}\footnotetext{$~^{**}$ Research supported in part by NSF grant
DMS-1406410.}

\begin{abstract}
Let $K$ be a convex set in $\R^d$ and let $K_\la$ be the convex hull of a homogeneous Poisson point process $\P_\la$ of intensity $\la$ on $K$.
When $K$ is a simple polytope, we establish scaling limits as $\la \to \infty$ for the boundary of $K_\la$ in a vicinity of a vertex of $K$
and we give variance asymptotics
for the volume and $k$-face functional of $K_\la$, $k \in \{0,1,...,d-1 \}$, resolving an open question posed in \cite{WW}.  The scaling
limit of the boundary of $K_\la$ and the variance asymptotics are described in terms of a germ-grain model consisting of cone-like grains pinned
to the extreme points of a Poisson point process on $\R^{d-1}\times \R$ having  intensity $\sqrt{d} e^{dh} dh dv$.
\end{abstract}

\section{Main results}\label{INTRO}
\allco
Let $K$ be a convex subset of $\R^d$ with non-empty interior.  For all $\la \in [1, \infty)$,
 let $\P_\la$  denote a homogeneous  Poisson point process of intensity $\la$
on $K$.
Let $K_\la$  be the polytope defined by the convex hull of $\P_\la$, with $f_k(K_\la)$
denoting the number of $k$-faces of $K_\la$, $k \in \{0,...,d-1\}$.
The study of the random polytope $K_\la$ has a long and rich history going back at least until $1864$ and we refer to  the surveys of Weil and Wieacker \cite{WW} and Reitzner \cite{ReBook} for details.
Major papers of R\'enyi and Sulanke \cite{RS,RS2} have played a seminal role in the subject.
When $K$ has a smooth boundary $\partial K$, it has been only recently shown by Reitzner \cite{Re} that $f_k(K_\la)$ and
$\Vol(K_\la)$ satisfy a central limit theorem as $\la \to \infty$.  More recently, for $\partial K$ smooth,  the second order properties and scaling limits of
$\partial K_\la$ have been established in \cite{CY, CSY, SY}.

When $K$ is itself a convex polytope, the analysis of $f_k(K_\la)$ and
$\Vol(K_\la)$ appears more challenging.  The lack of regularity in $\partial K$ as well as the lack of rotational symmetry in $K$
present additional technical obstacles.  Still, the central limit theorem   for
$f_k(K_\la)$ and $\Vol(K_\la)$ was shown in two remarkable papers of  B\'ar\'any and Reitzner \cite{BR2, BR},
who also establish rates of normal convergence for these functionals.    They do not consider scaling limits of $\partial K_\la$
and though they
obtain sharp lower bounds for  $\Var
f_k(K_\la)$ and $\Var \Vol (K_\la)$, their results stop short of determining
precise variance asymptotics  for $f_k(K_\la)$ and $\Vol (K_\la)$ as $\la \to \infty$, an open
problem going  back to the 1993 survey of Weil and
Wieacker (p. 1431 of \cite{WW}).
When $K$ is a simple polytope, we resolve this problem
in Theorems \ref{Th3} and \ref{Th4}, expressing variance asymptotics in terms
of scaling limit functionals of the germ-grain model consisting of  cone-like grains pinned to the `extreme' points
of the  Poisson point process $\P$
on $\R^{d-1} \times \R$ with
intensity \be \label{defP} d\P((v,h)):= \sqrt{d} e^{dh} dh dv, \ \  (v,h) \in \R^{d-1} \times \R.\ee
Along the way, we show that the scaling limit of $\partial K_\la$ near any vertex of $K$ coincides with the boundary of this
germ-grain model.

Our results share some striking similarities with their asymptotic counterparts
for convex hull functionals of i.i.d. uniform samples in the unit ball as well as
for i.i.d. Gaussian samples  in $\R^d$,  as given in \cite{CSY} and \cite{CY2}, respectively.
The remarkable qualitative similarities,
made precise in Remark (ii) below,  help unify both the second order analysis of random polytopes as well as
the scaling limit analysis of their boundaries.

\begin{figure}
  \centering
  \includegraphics[trim= 7cm 16cm 8.2cm 0cm, clip, scale=0.5]{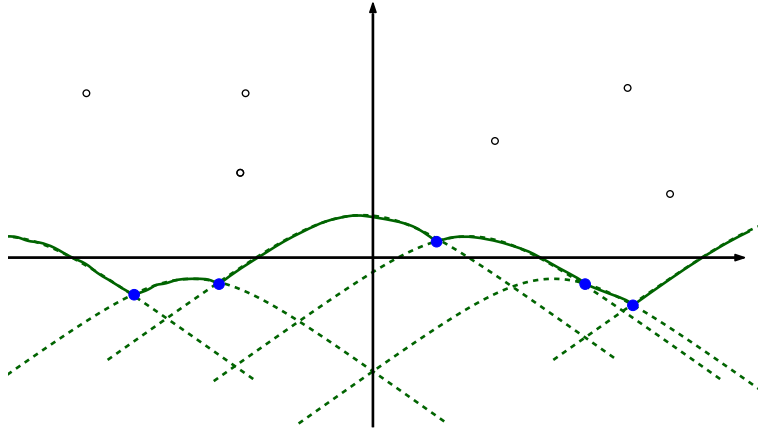}
\caption{The point process $\mbox{Ext}(\X)$ (blue); the boundary  $\partial\Phi(\X)$ of the down-grains containing $\mbox{Ext}(\X)$ (green).}
\label{fig:conespicture1}
\end{figure}

Before stating our  results we require some additional terminology. Henceforth we assume that $K$ is a simple polytope,
namely one whose vertices are adjacent to $d$ facets (faces of dimension $d-1$).
Let $V:=\{(x_1,\cdots,x_d)\in\R^d:\sum_{i=1}^d x_i= 0\}$
and for every $v\in V$ and $1 \leq i \leq d$, we let $l_i(v)$ be the $i$-th coordinate of $v$ in the standard basis with respect to $\R^d$ and $l(v)$ the vector $(l_1(v),\cdots,l_d(v))$ in $\R^d$.
Put  \be \label{defG}G(v):=\log(\frac{1}{d} \sum_{i=1}^de^{l_i(v)}), \ v \in V. \ee
The graph of $G$ has a cone-like structure, as shown in Lemma \ref{lemcompare}, and  gives rise to the cone-like grain
\be \label{Cone}
  \Pi^{\downarrow}:= \{(v,h) \in \R^{d-1} \times \R: \ h \leq -G(v) \}
\ee
opening in the down direction.
For $w:= (v,h) \in \R^{d-1} \times \R$ we
put $\Pi^{\downarrow}(w):= w \oplus
\Pi^{\downarrow}$, where $\oplus$ denotes Minkowski addition.
Given a locally finite  set $\X$ in $\R^d$, the maximal union of grains
 $\Pi^{\downarrow}(w), w \in \R^{d-1} \times \R,$ whose interior  contains no point of $\X$ is
\be \label{hull}
\Phi(\X):= \bigcup_{\left\{ \substack{w\in \R^{d-1} \times\R \\
\X \cap {\rm{int}}(\Pi^{\downarrow}(w)) =\emptyset}\right.} \Pi^{\downarrow}(w).
\ee
Remove points of $\X$  not belonging to $\partial (\Phi(\X))$
and call the resulting thinned point set ${\rm{Ext}}(\X)$. As shown in Figure \ref{fig:conespicture1},
$\partial (\Phi(\X))$ is a union of inverted cone-like surfaces `pinned' to or `suspended' from
${\rm{Ext}}(\X)$.

Extending the logarithmic function to $(0,\infty)^d$ by the formula $\log(z_1,\cdots,z_d)=(\log z_1,\cdots,\log z_d)$,
we consider for all $\la \in [1, \infty)$
\be \label{defT}
T^{(\lambda)}:\left\{\begin{array}{lll}(0,\infty)^d&\longrightarrow &V\times\R
\\(z_1,\cdots,z_d)&\longmapsto &
\left(p_V(\log z),\frac{1}{d}(\log \lambda +\sum_{i=1}^d\log z_i) \right)
\end{array}\right..
\ee
Here $p_V: \R^d \to \R$ denotes the orthogonal projection onto $V$.
Postponing the motivation behind $T^{(\la)}$ until Section 4,  we state our main results. 
Let $K':=[0,\Delta_d]^d$ where $\Delta_d \in [1, \infty)$ is a suitably large constant depending only on $d$, to be further
specified in the sequel (cf. Lemma \ref{hyper}).
Without loss of generality,
re-scaling $K$ if necessary, we make a volume preserving affine transformation such that
the origin is a vertex of $K$, $K' \subset K$, and $K$ is  contained in a multiple of $K'$.
Put \be \label{defr0} \delta_0:= \delta_0(\la) := \exp( - (\log \la )^{1/d})\ee
and let $Q_0 := [0, \delta_0]^d$.



\begin{theo} \label{Th1} Under the transformation $T^{(\la)}$,
the extreme points of $K_\la \cap Q_0$  converge in distribution to the thinned process ${\rm{Ext}}(\P)$ as
$\la \to \infty$.
\end{theo}

Let $B_d(x,r)$ be the closed $d$-dimensional Euclidean ball centered at $x \in \R^{d}$
and with radius $r \in (0, \infty)$.  $\C(B_{d}(x,r))$ is the space of
continuous functions on $B_{d}(x,r)$ equipped with the supremum norm. Let $\0$ denote a point at the origin of
$\R^d$ or $\R^{d-1}$, according to context.

 \begin{theo}\label{Th2} Fix $L \in (0, \infty).$
 As $\la \to \infty$, the re-scaled boundary  $T^{(\la)}( (\partial K_\la) \cap Q_0)$
 converges in probability to  $\partial (\Phi(\P))$ in the space $\C(B_{d-1}(\0,L))$.
 \end{theo}


The transformation $T^{(\la)}$ induces scaling
limit  $k$-face and volume functionals governing the large $\la$
behavior of $f_k(K_\la)$  and $\Vol(K_\la)$ as seen in the next results.


\begin{theo} \label{Th3}
For all $k \in \{0,1,...,d-1\}$, there exists a constant $F_{k,d}
\in (0, \infty)$, defined in terms of averages of covariances of a scaling limit
$k$-face functional on $\P$,  such that
 \be \label{main3a}
\liml (\log \la)^{-(d-1)}  \Var f_k(K_\la) =
F_{k,d} \cdot f_0(K). \ee
\end{theo}

\begin{theo} \label{Th4}
There exists a constant $V_{d}
\in (0, \infty)$, defined in terms of averages of covariances of a scaling limit
volume functional on $\P$,  such that
 \be \label{main4a}
\lim_{\la \to \infty} \la^{2} (\log \la)^{-(d-1)}  \Var \Vol(   K_\la) =
V_{d}\cdot f_0(K). \ee
\end{theo}
\vskip.2cm

\noindent {\em Remarks.} (i)  {\em On the scaling transform}.  Baryshnikov's far reaching work \cite{Ba} uses the scaling transform  $T^{(\la)}$,  though in
a different guise,  to establish Gaussian fluctuations and variance asymptotics  for the number of Pareto extreme points
in a random sample.    Baryshnikov  (cf. Section 2.2.5 of  \cite{Ba})
also mentions that $T^{(\la)}$ could be used to establish the asymptotic normality of $f_0(K_\la)$,
but he does not provide the details.  The present paper, in addition to establishing the scaling
limit of $\partial(K_\la)$ and the asymptotics of $\Var f_k(K_\la)$ and $\Var \Vol(K_\la)$,
makes a three-fold contribution going beyond that in \cite{Ba}.
First, we establish a new, if not crucial,
interpretation of the action of the scaling transform in terms of a dual process defined via cone-like grains called petals.
In the setting of convex hulls of i.i.d. samples in the unit ball, the dual process
has previously featured as a parabolic growth process  \cite{CSY, SY}.
Second, we establish a qualitative link with scaling transforms used previously for different models of random polytopes; see remark (ii) below and Section 4.1.   Lastly,
the transform suitably re-scales the {\em floating bodies} for $K$, showing that their re-scaled images play a central role in the re-scaled
convex hull geometry.    Floating bodies usefully  approximate random polytopes,  as in \cite{BR2,BR},  but here we show that their re-scaled images also play a key role in asymptotic analysis.  The approach surrounding the transform $T^{(\la)}$, together with the
counterpart transforms in  \cite{CY2, CSY, SY}, help unify the
asymptotic analysis of random polytopes.

\vskip.2cm
\noindent (ii) {\em Theorems \ref{Th1} and \ref{Th2} - related work}. The re-scaled point process $T^{(\la)}(\P_{\la})$ converges to the point process $\P$ as seen in Lemma \ref{imagePPP}. The part of the re-scaled boundary of $K_\la$ which is close to a vertex of $K$ converges to a festoon of inverted
cone-like hypersurfaces pinned to ${\rm{Ext}}(\P)$. The situation with $K$ a unit ball involves quantitative differences and similarities.
 When $K$ is the unit $d$-dimensional ball, in the large $\la$ limit, the relevant scaling  transform carries
 $\P_{\la}$  into a homogeneous Poisson point process on the upper half-space $\R^{d-1} \times \R^+$ and it carries the boundary of $K_\la$ into
 a festoon of {\em parabolic} hypersurfaces \cite{CSY, SY}.   On the other hand, if the input is a Poisson point process $\tilde{\P}_\la$ having Gaussian  intensity $\la \phi(x)dx$, with
 $\phi$ being the standard normal density on $\R^d$, then,  as $\la \to \infty$,  the relevant scaling transform carries
 $\tilde{\P}_\la$ into a non-homogeneous Poisson point process $\tilde{\P}$ on $\R^{d-1} \times \R$  with
intensity density $e^{h}dhdv$  and it carries the boundary of  the convex hull of  $\tilde{\P}_\la$ into  a festoon of parabolic hypersurfaces pinned against the extreme points of  $\tilde{\P}$ \cite{CY2}.
\vskip.3cm
\noindent (iii) {\em Theorems \ref{Th3} and \ref{Th4}}.
Variance asymptotics \eqref{main3a} and \eqref{main4a} do not depend on the volume of $K$, but only on the number of its vertices.
 Breakthrough papers of B\'ar\'any and Reitzner \cite{BR2, BR} establish precise growth rates for $\Var f_k(K_\la)$ and $\Var \Vol(K_\la)$.  While these works do not give a closed form expression for the asymptotic constants $F_{k,d}$ and $V_d$, they do insure their strict positivity.   We anticipate that methods given here yield expectation and variance asymptotics for
non-homogenous Poisson point processes having intensity density  $\la \ka$, with $\ka: K \to \R^+$ bounded away from zero and infinity and
continuous on $\partial K$.
\vskip.3cm
\noindent (iv) {\em The locally defined transform  $T^{(\la)}$}.  The  map  $T^{(\la)}$
is local in that it is defined  with respect to $\0$, assumed to be a vertex of $K$.
We are unable to find a suitable global transform for all of $K$.  On the other hand, when $K$ has rotational symmetry, e.g. when $K$ is the $d$-dimensional ball,
then we may globally map $K$ into $\R^{d-1} \times \R^+$ as in \cite{CSY, SY}.  The existence of a global scaling transform brings multiple benefits, leading  to a more regularized re-scaled structure  in $\R^{d-1} \times \R$, including stationarity as $\la \to \infty$ and  local independence
(stabilization) with respect to spatial coordinates.   When $K$ lacks rotational symmetry, as is the case in this paper, then our methods
do not yield any such global scaling transform. Roughly speaking, the obstruction to finding a global scaling transform goes as follows.
The  transform given here, like those  in
 \cite{CY2, CSY},  relies on the construction of a one-parameter family of $(d-1)$-dimensional surfaces interior to $K$
 (boundaries of associated floating bodies for $K$),  in which the  height coordinate is a function of the corresponding parameter and the space coordinate is given by a mapping from a subset of   $\R^{d-1}$ to the surface belonging to the one-parameter family.
It is in general difficult to
construct a global mapping from $\R^{d-1}$ to
a $(d-1)$-dimensional manifold, and thus difficult to find a global scaling transform.

\vskip.2cm
\noindent (v) {\em Approximate additivity of the variance}.  The lack of a global scaling transform necessitates showing that
$\Var f_k(K_\la),  \ k \in \{0,1,...,d -1 \},$ and $\Var \Vol(K_\la)$ are well approximated by the sum of variances of contributions arising from small
neighborhoods around each vertex of $K$.  We show this  decoupling of the variance over the vertex set of $K$ by refining the dependency graph arguments
in \cite{BR2} and applying these arguments to a dyadic collection of Macbeath regions.  These non-trivial technical obstacles are not
present when $K$ is the unit ball \cite{CSY}.
\vskip.2cm
\noindent   (vi) {\em Extension to general polytopes}.  The fundamental work of B\'ar\'any and Buchta \cite{BB} shows that the extreme points of a general polytope $K$ concentrate in regions defined by each flag of $K$.  These `flag regions',  themselves simple polytopes, could individually be treated by the methods of this paper.  If one could show (a) negligibility of covariances of contributions to $f_k(K_\la)$  and $\Vol(K_\la)$ arising from input on
distinct flag regions
as well as (b) negligibility of  contributions to $f_k(K_\la)$  and $\Vol(K_\la)$ arising from
input on the complements of flag regions, then variance would be additive with respect to flags.  This would extend our  results to general polytopes and it would
align our second order results with the first order results of  \cite{BB}, which shows that expectation asympotics
are additive with respect to flags.   However showing additivity of variance with respect to flags seems to be a separate project which would either require
a scaling transform more general than $T^{(\la)}$ or a non-trivial extension of the methods of Section 3.


\vskip.2cm
This paper is organized as follows. Section \ref{Sec2} introduces scaling limit functionals of germ-grain models having cone-like  grains.
These scaling limit functionals feature in Theorem \ref{Th5}, which establishes expectation and
variance asymptotics for the empirical measures for the volume and $k$-face functionals,
thus extending Theorems \ref{Th3} and \ref{Th4}.
In Section \ref{LEMMAS} we establish two propositions which prepare for an effective use of
the crucial transformation $T^{(\la)}$.   We  show that the extreme points near a vertex of $K$ have a preferred normal direction and that
the variance of the $k$-face and volume functional decouples over the vertices of $K$.  
Section \ref{transforms} studies $K_\la$  near a fixed vertex of $K$ and establishes that the image of $\P_\la$ under $T^{(\la)}$ converges in distribution to $\P$ and that $T^{(\la)}$
defines re-scaled $k$-face and volume functionals. We show  that  the scaling transform $T^{(\la)}$
maps the Euclidean convex hull geometry into `cone-like' convex geometry in $\R^{d-1} \times \R$ and  that the extreme points near a vertex of $K$ are with high probability characterized in terms of the geometry of {\em so-called} petals.
Section \ref{props} establishes that the re-scaled $k$-face and volume functionals
localize in space, which is crucial to showing convergence of their means and covariances to the respective means and covariances of
their scaling limits. Section \ref{proofs} provides proofs of the main results whereas  Section \ref{sec:app}, an Appendix, contains proofs of several  lemmas.


\section{General results}\label{Sec2}
\allco

Here we consider  functionals of the germ-grain model $\Phi(\P)$ which
are central to the description of the scaling limits of the $k$-face and volume
functionals $f_k(K_\la)$  and  $ \Vol(K_\la)$. We use their second order correlations to
establish variance asymptotics for the empirical measures induced by the
$k$-face and volume functionals.  This extends Theorems \ref{Th3} and \ref{Th4} to the setting of measures
and yields formulas for the constants $F_{k,d}$ and $V_d$ in \eqref{main3a} and \eqref{main4a}, respectively.
Denote points in $\R^{d-1} \times \R$ by $\v:= (v,h)$.



\vskip.5cm

\noindent{\bf 2.1. Empirical $k$-face  and volume
measures.} Given a finite point set $\X \subset \R^d$, let
$\rm{co}(\X)$ be its convex hull.

\begin{defn} \label{kface} Given $k \in \{0,1,...,d-1\}$
and $x$ a vertex of $\rm{co}(\X)$, define the $k$-face functional
$\xi_k(x, \X)$ to be the product of $(k +1)^{-1}$ and the number of
$k$-faces of $\rm{co}(\X)$ which contain $x$. Otherwise we put
$\xi_k(x, \X) = 0$. Thus the total number of $k$-faces in $\rm{co}(\X)$ is $\sum_{x \in \X}
\xi_k(x, \X)$. Letting $\de_x$ be the unit point mass at $x$, the empirical k-face measure for $\P_\la$ is \be
\label{zerom} \mu^{\xi_k}_\la:= \sum_{x \in \P_\la} \xi_k(x, \P_\la
) \de_x. \ee
\end{defn}

\vskip.5cm

We now consider the defect volume of $K_\la$ with respect to $K$ in a cubical neighborhood of $\0$. 
Recall from \eqref{defr0} that $\delta_0:= \delta_0(\la) := \exp( - (\log \la )^{1/d})$
and $Q_0 := [0, \delta_0]^d$.  

Let $\X \subset [0, \infty)^d$ be finite.
Given $x$ a vertex of $\rm{co}(\X)$, let ${\cal F}^+(x, \X)$ be the
(possibly empty) collection of facets in $\rm{co}(\X)$, included in $Q_0$,
 containing $x$ and having outer normals in $(-\infty, 0]^d$.  Let ${\rm{cone}}(x, \X):=
\{ry:  \ r > 0, y \in {\cal F}^+(x, \X) \}$
 be the (possibly empty) cone generated by ${\cal F}^+(x, \X)$.

\begin{figure}[H]
\label{thelastpictureshow}
\centering
\includegraphics[trim= 0cm 0cm 2cm 2.5cm, clip, scale=0.5]{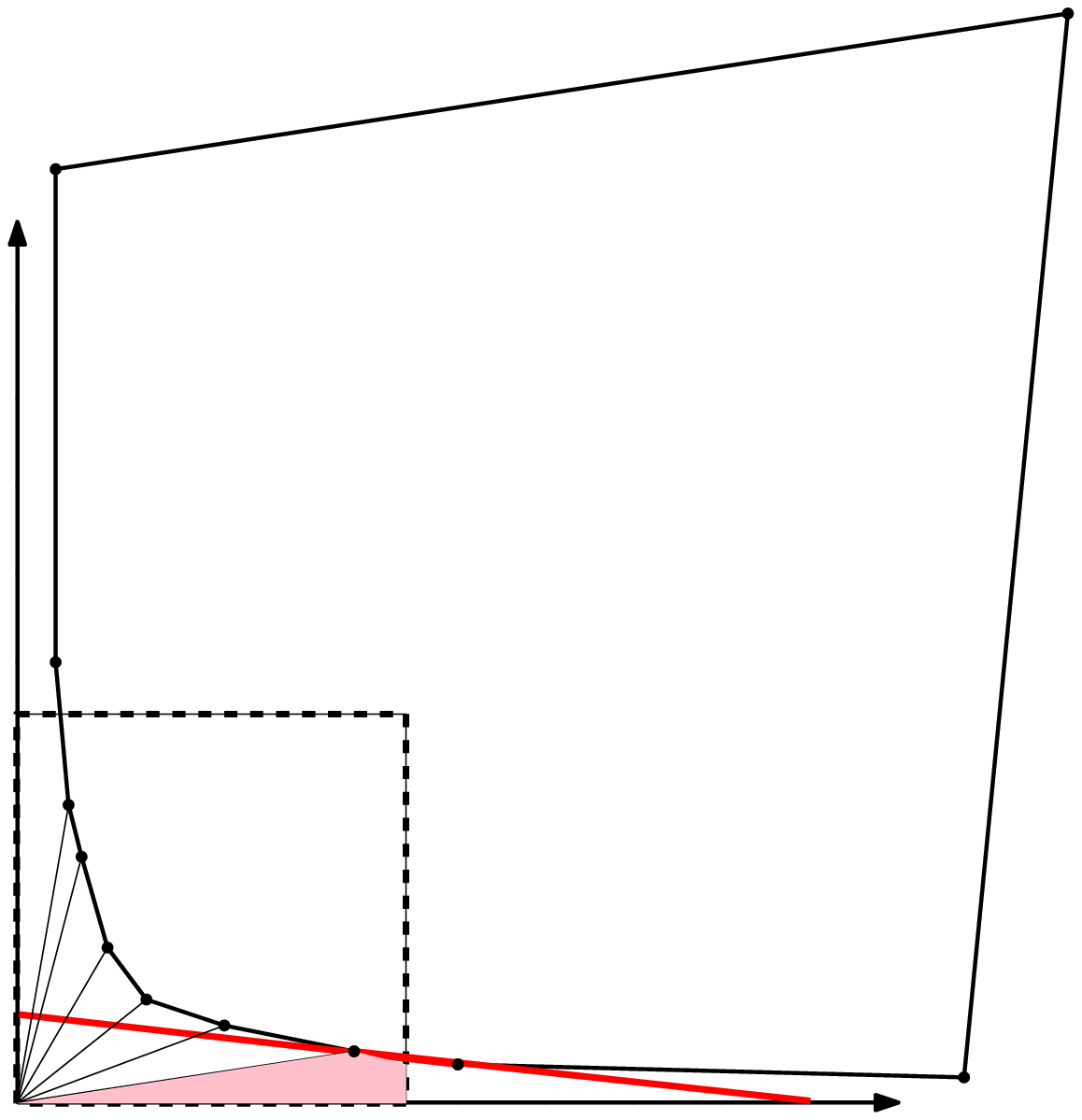}
\caption{The part of the defect volume not counted in the sum of scores (pink); the boundary of the cap containing that part (red).}
\end{figure}

We first define the volume score for points in $\P_\la\cap Q_0$. If $x$ is vertex of $K_\la \cap Q_0$ with ${\rm{cone}}(x, \P_\la) \neq \emptyset$, then
define the defect volume functional
\be\label{volumescore1}
\xi_V(x, \P_\la):= d^{-1} \la \Vol( {\rm{cone}}(x, \P_\la) \cap (K\setminus K_\la) ).
\ee
Otherwise put $\xi_V(x, \P_\la) = 0$. 
In general, the union $\bigcup_x ( {\rm{cone}}(x, \P_\la) \cap (K\setminus K_\la) )$, with $x$ ranging over the
vertices of $K_\la \cap Q_0$, does not equal $(K\setminus K_\la) \cap Q_0$.  See Figure  \ref{thelastpictureshow}.

We now extend the definition of the volume score so that is defined for all points
in $\P_\la$.  This goes as follows.
Let ${\cal V}_K := \{{\mathscr V}_i \}$ denote the vertices of $K$ and recall that we assume $\0 \in {\cal V}_K$. Re-scaling $K$ if necessary,
for each vertex ${\mathscr V}_i \in {\cal V}_K \setminus \{\0\}$, we define an associated volume preserving affine transformation $a_i: \R^d \to \R^d$,  with $a_i({\mathscr V}_i) = \0$,
and such that the facets of $a_i(K)$ containing $\0$ also contain the facets of $K':= [0, \Delta_d]^d$
belonging to the coordinate hyperplanes.  For any $\delta \in (0, \Delta_d)$, define the parallelepiped $p_d({\mathscr V}_i, \delta):= a_i^{-1}( [0, \delta]^d )$.  Put $\P_\la(\delta) := \P_\la \cap \bigcup_{1\le i\le f_0(K)} p_d({\mathscr V}_i,\delta)$. For any convex polytope $K$
 and any face $F$ of $K$, we denote by $C_F(K)$ the cone of outer normal vectors to $F$. For any $x\in p_d({\mathscr V}_i, \delta_0)$ define  ${\mathcal F}^+(x,\P_\la)$ to be the intersection of $p_d({\mathscr V}_i, \delta_0)$ and the (possibly empty) collection of facets in ${\rm co}(\P_\la)$, included in $p_d({\mathscr V}_i, \delta_0)$, containing $x$ and having outer normals in $C_{{\mathscr V}_i}(K)$. We then put  ${\rm {cone}}(x,\P_\la)=\{{\mathscr V}_i \oplus r(y-{\mathscr V}_i): r>0,y\in {\mathcal F}^+(x,\P_\la)\}$.
If $x$ is vertex of $K_\la \cap p_d({\mathscr V}_i, \delta_0)$ with ${\rm{cone}}(x, \P_\la) \neq \emptyset$, define the defect volume functional via the formula \eqref{volumescore1} and put $\xi_V(x, \P_\la) = 0$ otherwise.



Finally, we define the volume score for points in $\P_\la \setminus \P_\la(\delta_0)$. Let ${\cal F}_k(x):= {\cal F}_k(x, K_\la)$ be the (possibly empty) collection of $k$-dimensional faces in
$K_\la$ which contain $x$. For $x$ a vertex of  $K_\la$ and
$x \notin \P_\la(\delta_0)$,
define the defect volume score 
\be\label{volumescore2}
\xi_V(x,\P_\la):=\sum_{k=0}^{d-1}\sum_{F \in {\cal F}_k(x) }\frac{\la}{\card\{F\cap (\P_{\la}\setminus \P_\la(\delta_0))\}} \Vol(  (F \oplus C_F(K_\la)) \cap K  \cap ( D(\P_\la))^c),
\ee
where
$D(\P_\la):=  \bigcup_{x\in \P_\la(\delta_0) } {\rm{cone}}(x,\P_\la).$
Otherwise, put $\xi_V(x,\P_\la)=0$.

\begin{defn} \label{vol} 
Define the empirical defect volume measure of $\P_\la$ 
by \be \label{volmeasure}
\mu^{\xi_V}_\la:= \sum_{x \in \P_\la} \xi_V(x, \P_\la) \de_x, \ee
where $\xi_V(x,\P_\la)$ is given by \eqref{volumescore1} or \eqref{volumescore2} depending on whether $x\in \P_\la(\delta_0)$ or not.
\end{defn}
Notice that in general $\sum_{x \in \P_\la \cap Q_0} \xi_V(x, \P_\la) \leq \la \Vol (Q_0 \cap (K \setminus K_\la) )$
and likewise we have  $ \langle 1, \mu^{\xi_V}_\la \rangle \leq \la \Vol(K \setminus K_\la)$.

\vskip.5cm

\noindent{\bf 2.2. Scaling limit $k$-face and volume functionals}. A set of $(k+1)$ extreme
  points $\{w_1,...,w_{k + 1} \} \subset {\rm{Ext}}(\P)$,  generates a $k$-dimensional {\em face}
 of the  festoon  $\partial (\Phi(\P))$ if there exists a translate $\tilde{\Pi}^{\downarrow}$ of $\Pi^{\downarrow}$
such that $\{w_1,\cdots,w_{k+1} \}=  \tilde{\Pi}^{\downarrow} \cap
{\rm{Ext}}(\P)$.
When $k = d -1$ the face is a {\em hyperface}.

\begin{defn} \label{xiinf}  Let $w \in {\rm{Ext}}(\P)$.
Define the scaling limit $k$-face functional $\xi^{(\infty)}_{k}(w, \P)$,  $k \in \{0, 1,..., d - 1\}$, to be the product of
$(k + 1)^{-1}$ and the number of $k$-dimensional faces of
the festoon $\partial (\Phi(\P))$ which
contain $w$.  The scaling limit defect volume functional is
$$
\xi^{(\infty)}_{V}(w, \P):= \frac{ 1} { d \sqrt{d} }  \int_{ v \in {\rm{Cyl}}(w)}
\exp\{ d \cdot \partial (\Phi(\P))(v) \} dv,
$$
where  ${\rm{Cyl}}(w)$ denotes the projection onto  $\R^{d-1}$ of
the hyperfaces of $\partial( \Phi(\P))$ containing $w$.
When $w \notin {\rm{Ext}}(\P)$ we put
$\xi^{(\infty)}_{k}(w, \P)= 0$ and $\xi^{(\infty)}_{V}(w, \P)= 0$.
\end{defn}


Let $\Xi$ denote the collection of functionals
$\xi_{k}, k \in \{0,1,...,d-1\},$ together with
$\xi_{V}$.
Let $\Xi^{(\infty)}$ denote the collection of scaling limits
$\xi^{(\infty)}_{k}, k \in \{0,1,...,d-1\},$ together with
$\xi^{(\infty)}_{V}$.  A main feature of our approach (cf. Lemma \ref{L2-new}) is that on a high probability set,
the elements of $\Xi^{(\infty)}$ are scaling limits of re-scaled elements of
$\Xi$.  

\vskip.5cm

\noindent{\bf 2.3. Limit theory for  empirical $k$-face
and volume measures.} Define the following second order correlation
functions for $\xi^{(\infty)} \in \Xi^{(\infty)}$.

\begin{defn} For all $\v_1, \v_2 \in \R^d$ and $\xi^{(\infty)} \in \Xi^{(\infty)}$ put
\be \label{SO2} c^{\xi^{(\infty)}}(\v_1,\v_2):=
c^{\xi^{(\infty)}}(\v_1, \v_2, \P):= \ee $$
 \E \xi^{(\infty)}(\v_1, \P
\cup \{\v_2\} ) \xi^{(\infty)}(\v_2, \P \cup \{\v_1\} ) -  \E
\xi^{(\infty)}(\v_1, \P ) \E \xi^{(\infty)}(\v_2, \P )$$
and  \be \label{S03}
 \sigma^2(\xi^{(\infty)}) := \sqrt{d} \int_{-\infty}^{\infty} \E \xi^{(\infty)}((\0,h_0), \P)^2 e^{dh_0} dh_0  \ee
  $$ + \ d\int_{-\infty}^{\infty}  \int_{\R^{d-1}} \int_{-\infty}^{\infty}
   c^{\xi^{(\infty)}}((\0,h_0),(v_1,h_1)) e^{d(h_0 + h_1)} dh_1 dv_1 dh_0.
$$
\end{defn}
\vskip.5cm


Let $\C(K)$ be the class of
bounded functions on $K$ which are continuous on ${\cal V}_K$.
Given $g \in \C(K)$ let $\langle g, \mu_\la^\xi \rangle$ denote its integral with
respect to $\mu_\la^\xi$.  Consider  the regular $d$-dimensional simplex
of edge length $\sqrt{2} d$ given by
\be \label{defSd}
S(d):= \{(x_1,...,x_d)\in (-\infty,1\blue{]}^d: \ \sum_{i=1}^d x_i = 0 \}.
\ee
The following general result is proved in Section \ref{proofs}.
\vskip.5cm

\begin{theo} \label{Th5} For all $\xi\in \Xi$
and $g \in \C(K)$  we have \be \label{main1} \lim_{\la \to \infty} (\log \la)^{-(d-1)}
\E [\langle g, \mu_\la^{\xi} \rangle] =
d^{-d + (3/2)}   \Vol_d(S(d)) \int_{-\infty}^\infty \E \xi^{(\infty)}((\0,h_0), \P) e^{dh_0} dh_0 \sum_{{\mathscr V}_i \in {\cal V}_K }  g({\mathscr V}_i)
 \ee and \be \label{main2} \lim_{\la \to
\infty} ( \log \la)^{-(d-1)}  \Var [\langle g,
\mu_\la^{\xi} \rangle] = d^{-d+1}  \Vol_d(S(d)) \sigma^2(\xi^{(\infty)}) \sum_{{\mathscr V}_i \in {\cal V}_K }   g^2({\mathscr V}_i) . \ee
\end{theo}

\vskip.5cm

\noindent{\bf 2.4. Deducing Theorems \ref{Th3} and \ref{Th4}.} We may deduce the  limits \eqref{main3a} and \eqref{main4a}  from
Theorem \ref{Th5} as follows. Set $\xi$ to be $\xi_k$ and put $g \equiv 1$.
Then $\langle 1, \mu_\la^{\xi_k} \rangle = f_k(K_\la)$ and so  \eqref{main3a} follows from \eqref{main2}, setting $F_{k,d}$ to
be   $d^{-d+1}  \Vol_d(S(d)) \sigma^2(\xi^{(\infty)}_k)$.  Next, set $\xi$ to be $\xi_V$.  By Lemma \ref{L4}(b) we have
$$
\lim_{\la \to \infty} \frac{ \Var \langle 1, \mu_\la^{\xi_V} \rangle } { (\log \la)^{d-1} } = \lim_{\la \to \infty} \frac
{\la^2 \Var( \Vol K_\la) } { (\log \la)^{d-1} }.
$$
Thus \eqref{main4a} follows from \eqref{main2} 
where we set $V_{d}$ to
be  $d^{-d+1}  \Vol_d(S(d)) \sigma^2(\xi^{(\infty)}_V)$.


\section{Decomposition of the variances}\label{LEMMAS}
\allco

Before discussing the scaling transform $T^{(\la)}$ we need some key simplifications.  We show here that variance asymptotics for
 $f_k(K_\la)$ and $\Vol(K_\la)$ are determined by the behavior of these functionals  on points near any fixed vertex of $K$, assumed without loss of generality to be $\0$; {\em that is
 to say the point set $\P_\la \cap Q_0$ determines variance asymptotics. } It is far from clear that this should be the case, as covariances of scores
on subsets of $\P_\la$ near adjacent vertices of $K$ might be non-negligible.   Secondly, variances of scores
on subsets of $\P_\la$ `between' adjacent vertices of $K$ might also be non-negligible.  The purpose of this section is
to address these two issues via  Proposition \ref{Prop2}, showing the negligibility of the afore-mentioned quantities.
This paves  the way for an effective use of $T^{(\la)}$, which is well defined on $\P_\la \cap Q_0$,
and which we use in Section \ref{transforms} to re-scale the scores $\xi \in \Xi$, the input $\P_\la$, as well as $Q_0 \cap \partial K_\la$.

\begin{defn} \label{cone-ext} If the collection $C_F(K_\la \cap Q_0)$ of outward normals to a face $F$ in $K_\la \cap Q_0$ all belong to the normal cone $C_{\0}(K):= \{u = (u_1,...,u_d) \in (-\infty, 0)^d \}$,  then $F$ is called a `cone-extreme' face.
\end{defn}

Before stating Proposition \ref{Prop2} we require an auxiliary result, whose proof is deferred to the Appendix.  We first
assert  that there is a high probability event $A_\la$, to be defined in Section 3.1, such that on $A_\la$  all faces of
$K_\la \cap Q_0$ are cone-extreme. It is precisely the cone-extreme faces which are amenable to
analysis under the transformation $T^{(\la)}$.

\begin{prop} \label{Prop1} There is an event $A_\la$ with $P[A_\la^c] \leq C (\log \la)^{-4d^2}$, such that on
$A_\la$ we have $C_F(K_\la \cap Q_0) \subset C_{\0}(K)$ for any face $F$ of $K_\la \cap Q_0$.
\end{prop}
Sections 4 and 5 develop  the geometry and scaling properties of
cone-extreme faces. In particular  Lemma \ref{lempetal} identifies their collective image under $T^{(\la)}$ with a festoon of inverted cone-like surfaces.

Next, for each $\xi \in \Xi$, put
\be \label{defZ}
Z:= Z_\la := \sum_{x \in \P_\la} \xi(x, \P_\la).
\ee

The contribution to the total score from points in $\P_\la  \cap p_d({\mathscr V}_i, \delta)$ is
\be \label{Wl}
Z_i:= Z_{i}(\delta):= \sum_{x \in \P_\la \cap  p_d({\mathscr V}_i, \delta) } \xi(x, \P_\la), \ \ 1 \leq i \leq f_0(K).
\ee
We now choose  $\delta:= \delta(\la)$ large enough so that $\sum_{i = 1}^{f_0(K)} Z_{i}(\delta)$ captures the bulk of the total score
$Z$, but small enough so that  $Z_{i}(\delta)$ are independent random variables, or at least conditionally so,
given the event $A_\la$ of Proposition \ref{Prop1}.  The next proposition tells us that it suffices to  set $\delta$ to be $\delta_0$
and it shows that $\Var Z$ is essentially a sum of variances of
scores induced by points in $\P_\la$ near each vertex of $K$.  Given two sequences of scalars $\alpha_\la$ and $\beta_\la$,  $\la > 0$,
we write  $\alpha_\la = o (\beta_\la)$ if $\alpha_\la/\beta_\la \to 0$ as $\la \to \infty$.

\begin{prop} \label{Prop2} 
For all $\xi \in \Xi$ we have
\be \label{decomp-e}
\E [Z {\bf{1}}(A_\la) ] = \sum_{{\mathscr V}_i \in {\cal V}_K } \E [Z_{i}(\delta_0) {\bf{1}}(A_\la)] +  o(\E [Z])
\ee
and
\be \label{decomp-a}
\Var [Z   {\bf{1}}(A_\la)  ] = \sum_{{\mathscr V}_i \in {\cal V}_K } \Var [Z_{i}(\delta_0) {\bf{1}}(A_\la)] +  o( \Var [Z]) =
\sum_{{\mathscr V}_i \in {\cal V}_K } \Var [Z_{i}(\delta_0)] +  o( \Var [Z]).
\ee
\end{prop}

Proposition  \ref{Prop2} shows that to prove Theorems \ref{Th3} and \ref{Th4}, it is enough to establish the variance of the
$k$-face and volume functional for that part of $K_\la$ included in $Q_0$. 
The identity \eqref{decomp-e} is essentially a re-phrasing of Theorems 3 and 4 in \cite{BB},  which show that $\E [Z]$ is a sum
of expectations of scores induced by points in $\P_\la$ near each vertex of $K$ (and more generally, near each flag of $K$ when $K$ is
an arbitrary convex polytope).    The methods of \cite{BB} do not appear to extend to variances.

To prove these two propositions, we shall rely heavily on a construction of dyadic Macbeath regions.
The rest of this section is devoted to the proof of Proposition
\ref{Prop2}. The set-up  of the next three subsections closely parallels that of the breakthrough paper  \cite{BR2}.

\vskip.5cm

\noindent{\bf 3.1. A critical annulus and a high probability set}.  As in \cite{BR2}, define  $v: K \mapsto \R$ by 
$$
v(z) := \text{min} \{ V(K \cap H): \ H \ \text{is  a  half-space  and  z} \ \text{in} \ H  \}.
$$
There should be no  confusion between the function $v$, used in this section and
in the Appendix, and  the point $v$, denoting a
generic point in $\R^{d-1}$,  used in subsequent sections.
For $t \in [0, \infty)$, let $K(v = t)$ be the boundary of the floating body $\{z \in K: \ v(z) \geq t\}$, which we abbreviate as
$K(v \geq t)$.
Recall $K':= [0, \Delta_d]^d$, with $\Delta_d \in [1, \infty)$ to be specified.
 Lemma \ref{hyper} in the Appendix shows that for $\Delta_d$ large,   the floating bodies for $K$ and $K'$  coincide
in $[0,1/2]^d$.
Following  \cite{BR2}, put
\be \label{defpara}
s := s_\la:= \frac{1}{ \la (\log \la)^\beta }, \ \ \ \ T
:= \frac{\alpha \log \log \la}{ \la },  \ \ \ T^*:= d6^d T
\ee
with $\beta := 4d^2 + d - 1, \alpha:= (6d)^d \beta$ (in this section,  $T$ denotes the scalar at \eqref{defpara} and there
should be no confusion with $T^{(\la)}$).  Consider the annulus-like set
$$
\A(s,T^*,K):= K(v \geq s) \setminus K(v \geq T^*).
$$
By Lemma 5.2 of \cite{BR2} there is an event $A_\la := A_\la(K)$ such that
on $A_\la$ we have
\be \label{defA} \partial K_\la \subset \A(s, T^*, K), \ {\rm{where}} \
(\log \la)^{- (3d)^{d + 2} } \leq P[A_\la^c] \leq C (\log \la)^{- 4d^2 }.
\ee

\vskip.5cm

\noindent{\bf 3.2. Macbeath regions}. In this subsection we construct Macbeath regions near the origin.
As we shall see, the construction serves as a prototype for constructing Macbeath regions near  vertices ${\mathscr V}_i \in {\cal V}_K \setminus \{\0\}$.

For all $z \in K$, let $M_K(z):= M_K(z, 1/2)$ be the Macbeath region
(M-region for short) with center $z$ and scale factor $1/2$, i.e.,
$$
M_K(z):= M_K(z, 1/2) := z + {1 \over 2}[(K - z) \cap (z - K)].
$$
For $z: =(z_1,...,z_d) \in [0, 1/2]^d$  we have
\be M_{K'}(z) = \prod_{i=1}^d [{z_i \over 2} , {3z_i\over  2}]. \label{subbox} \ee
The inclusion
$K' \subseteq K$  gives for all  $z: =(z_1,...,z_d) \in [0, 1/2]^d$
$$ M_{K}(z) = \prod_{i=1}^d [{z_i \over 2} , {3z_i\over  2}].$$ 
More generally, given $\delta \in (0, 1/2)$ 
and
integers $k_i \in \Z$ with  $3^{k_i} \in (0, 1/(3 \delta)], \ 1 \leq i \leq d$,  the dyadic rectangular solids
\be \label{Mreg}  \prod_{i = 1}^d [ {3^{k_i} \delta \over 2}, {3^{k_i + 1} \delta \over 2} ]
\ee
coincide with  the $M$-regions
\be \label{Mreg1} M_{K'}((3^{k_1} \delta,...,3^{k_d} \delta)) =  M_{K}((3^{k_1} \delta,...,3^{k_d} \delta)), \ \ 3^{k_i} \in (0, 1/(3 \delta)] .\ee

Points $z:= (3^{k_1} \delta,..., 3^{k_d} \delta)$ are centers of dyadic solids. When $\log_3 T/\delta^d \in \Z$,
then $M_K(z)$  has center $z$ belonging to $K(v = T)$ as soon as $\sum_{i=1}^d k_i=\log_3 T/\delta^d$; we shall use such $M$-regions to define a saturated system as in \cite{BR2}.

\begin{figure}
  \centering
  \includegraphics[trim= 0cm 15cm 0cm 0cm, clip, scale=0.45]{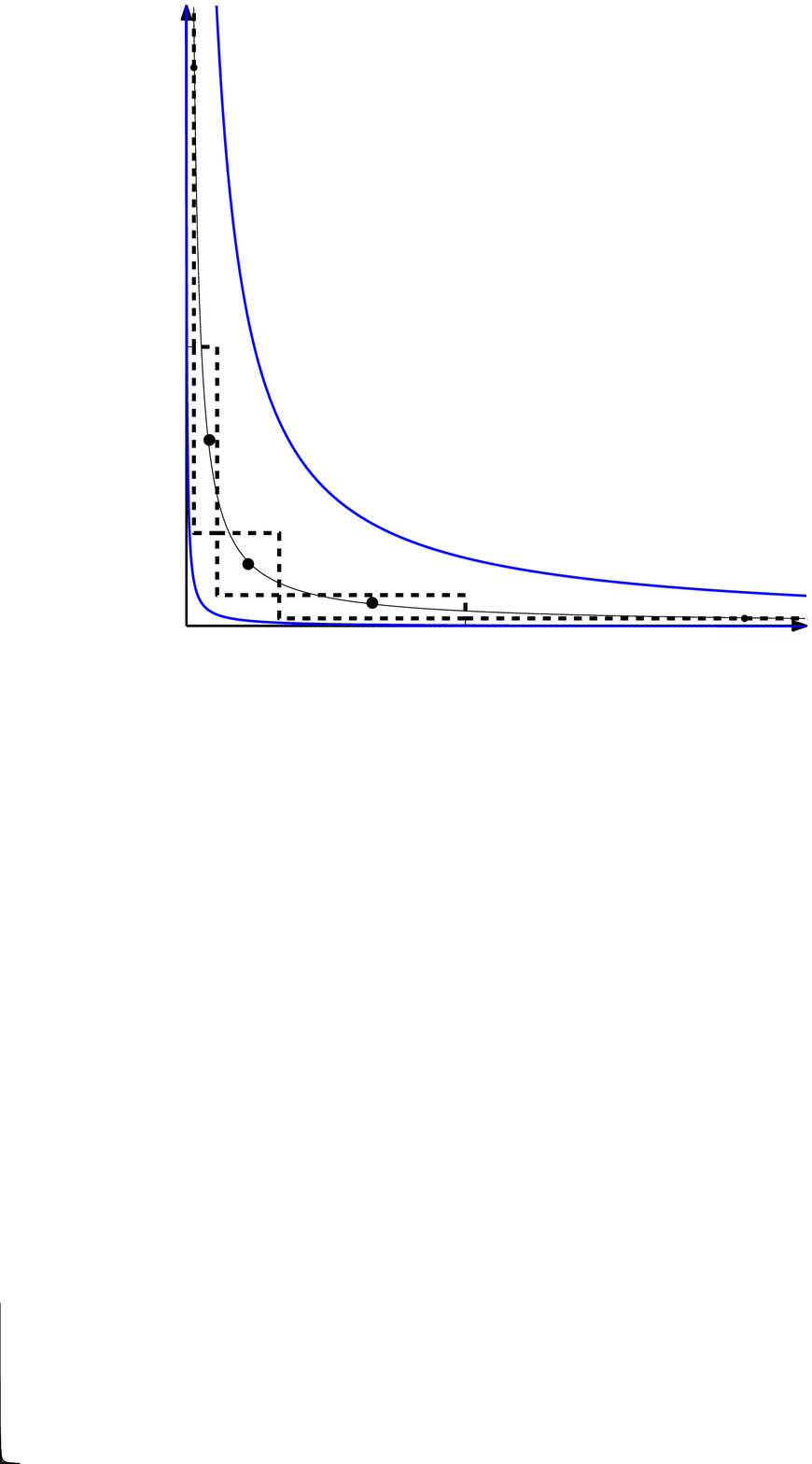}
  \caption{A saturated collection $\M_K(\0,\delta)$ of Macbeath regions}
  \label{macbeathpicture}
\end{figure}

Henceforth, let 
$\delta \in (0, 1/2)$
{\em and with $\log_3 T/\delta^d \in \Z$.}
Consider the collection $\M_K(\0,\delta)$ of dyadic rectangular solids of the type \eqref{Mreg1}
having centers on $K(v=T)\cap [0,1/2]^d$ (see Figure \ref{macbeathpicture}).
The solids in $\M_K(\0,\delta)$
do not cover $K(v = T) \cap [0, 1/2 ]^d$ but they leave  some parts uncovered.
The uncovered part is too small to accommodate another $M$-region with center on  $K(v = T) \cap [0, 1/2 ]^d$.
In other words, the collection $\M_K(\0,\delta)$ of dyadic $M$-regions  defined at \eqref{Mreg1} is maximal in that it can not be enlarged to include another $M$-region with center on
$K(v = T) \cap [0, 1/2 ]^d$. The following is proved in the Appendix.

\begin{lemm} \label{prop1}
The collection  ${\mathcal M}_K(\0, \delta)$ of $M$-regions is maximal.
\end{lemm}

We will use the collection ${\mathcal M}_K(\0, \delta)$ to control the spatial dependence of scores
$\xi \in \Xi$. This is done via supersets of $M$-regions, described below.
\vskip.5cm





\begin{figure}
  \centering
  \includegraphics[trim= 0cm 15cm 0cm 0cm, clip, scale=0.45]{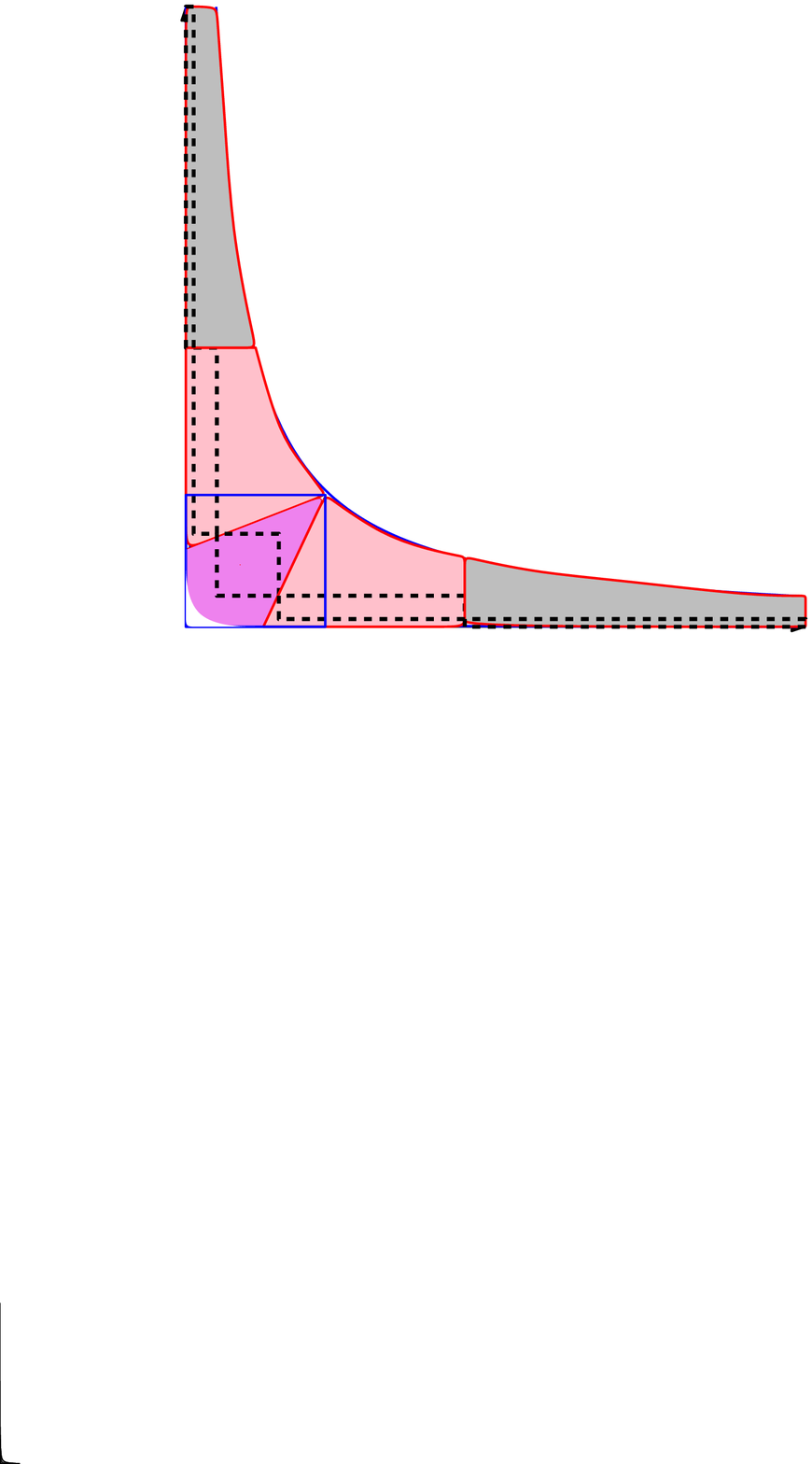}
  \caption{The supersets $S'_j$ (purple, pink and grey) associated with the saturated system of Macbeath regions.}
\label{fig:superset}
\end{figure}

\noindent{\bf 3.3. Supersets of M-regions.}
The collection $\M_K(\0,\delta)$   generates a `dyadic
staircase',  where  the
step width  increases in a geometric progression according to its distance from a coordinate hyperplane $H_l$,
$1 \leq l \leq d.$
Elements of  ${\mathcal M}_{K}(\0,\delta)$  are only pairwise interior-disjoint and not pairwise disjoint. Still, we claim that this collection allows us to reproduce the construction from the economic cap covering theorem (Theorem 2.5 of \cite{BR2}) with the same outcome and to construct a partition of $K(v \leq T^*)$ into supersets $S_j'$ which are also pairwise interior-disjoint (see Figure \ref{fig:superset}).  This goes as follows.

Each $M$-region in $\M_K(\0,\delta)$ produces a superset, called an $S$-region in \cite{BR2},  in the following canonical way.
For  $M$-regions $M_j$ meeting  $[0, (T^*)^{1/d}]^d$  we define the
associated region $S'_j:= S_j'(M_j)$  to be the intersection of
$K(v \leq T^*)$ and the smallest cone with apex at $((T^*)^{1/d},..., (T^*)^{1/d})$ which contains $M_j$.  We call these
the `cone sets $S'_j$'.   The volume of every $M$-region in $\M_K(\0,\delta)$ is $\Pi_{i=1}^d 3^{k_i} \delta = T$, and thus the number of $M$-regions
meeting $[0, (T^*)^{1/d}]^d$ is bounded by a constant depending only on $d$. The `cone sets $S'_j$' are not  contained in
 $[0, (T^*)^{1/d}]^d$ when $M_j$ itself is not  contained in
 $[0, (T^*)^{1/d}]^d$. In this case,  we replace the cone set $S'_j$ with a so-called `cone-cylinder set', defined to
 be the union of  $S'_j \cap  [0, (T^*)^{1/d}]^d$ and the so-called `cylinder set' generated by  $M_j \cap ([0, (T^*)^{1/d}]^d)^c$,  defined as follows.

For  $M$-regions $M_j$ with centers $(3^{k_1} \delta,...,3^{k_d} \delta)$
and such that $M_j$ meets $([0, (T^*)^{1/d}]^d)^c$,  define for $1 \leq l \leq d$, the cylinder
$$C_l(k_1,\cdots,k_d):=  \prod_{i=1}^{l-1} [ {3^{k_i} \delta \over 2}, {3^{k_i + 1} \delta \over 2} ]\times \R\times \prod_{i=l+1}^{d} [ {3^{k_i} \delta \over 2}, {3^{k_i + 1} \delta \over 2} ]  \cap ([0, (T^*)^{1/d}]^d)^c.$$
Note that $C_l(k_1,\cdots,k_d)$ is the smallest cylinder containing $M_j$ and oriented in the direction $n_{H_l}$,
where $n_{H_l}$ is a unit normal vector for  the hyperplane $H_l$.
The $\tilde{S}_j$ region associated with $M_j \cap ([0, (T^*)^{1/d}]^d)^c$  is
$$\tilde{S}_j:=\tilde{S}_j((3^{k_1} \delta,...,3^{k_d} \delta)):=\bigcup_{l:k_l=\min(k_1,\cdots,k_d)}C_l(k_1,\cdots,k_d)
 \cap K.$$
When $k_l$ is the unique minimum, $\tilde{S}_j$ consists solely of a single cylinder $C_l$ and it
simply extends $M_j \cap ([0, (T^*)^{1/d}]^d)^c$ in the direction
 $n_{H_l}$.  Note that $n_{H_l}$ points in the direction of the facet of $K'$ closest to $M_j$.
The union of such $\tilde{S}_j$ does not cover all of $K(v\le T^*) \cap ([0,1/2]^d  \setminus [0, (T^*)^{1/d}]^d)$.   The uncovered parts are rectangular regions produced by precisely one $M$-region having
a cubical face. 
Consequently, we define $S_j$ as the union of $\tilde{S}_j$ and all  rectangular regions  produced by the ties in the minimum of $k_1,\cdots,k_d$.


As on page 1518 of \cite{BR2}, we
define the superset \be \label{Sprimej}S'_j:= S_j \cap K(v \leq T^*).\ee
We  call these the `cylinder sets'.
If $S'_j$ has a facet $F$ which meets $\partial([0, (T^*)^{1/d}]^d)$ then
we adjoin the cylinder set $S'_j$ to the unique cone set containing $F \cap [0, (T^*)^{1/d}]^d$  and we call the resulting set a cone-cylinder set.
By construction, the sets $S'_j$ are disjoint.
Given this way of generating $S'_j, 1 \leq j \leq \text{card}(\M_K(\0,\delta)),$
we may control its diameter in all directions $n_{H_i}, 1 \leq i \leq d$.
The diameter of $M_{K}((3^{k_1} \delta,..., 3^{k_d} \delta))$
in the direction $n_{H_i}$ is  $$
{3^{k_i + 1} \delta \over 2}  - {3^{k_i} \delta \over 2} = 3^{k_i} \delta.$$
The diameter of the corresponding $S_j:= S_j ((3^{k_1} \delta,..., 3^{k_d} \delta))$ in the direction $n_{H_i}, 1 \leq i \leq d,$ is
thus at most $c3^{k_i} \delta+ 3^{k_i} \delta$, with the first term accounting for  the possible uncovered adjoined regions included in $S_j$,  or, if the direction $n_{H_i}$ were
 the direction $i$ corresponding to the smallest $k_i$ from the M-region, then  the diameter would be the distance from the coordinate hyperplane to the pseudo-hyperboloid $K(v=T^*)$ inside the cylinder $C_i$ and thus would be bounded by $c3^{k_i} \delta$.

Let $S'_1,...,S'_J$ be the supersets generated by $M$-regions meeting $[0,1/2]^d \setminus [0, (T^*)^{1/d}]^d$.
For $0 < a < b < \infty$ and $1\le i\le d$, we denote by $H_i[a,b]$ the `parallel slab' between hyperplanes $H_i\oplus an_{H_i} $ and $H_i \oplus bn_{H_i}$. We also  define for any bounded subset $A$ of $\R^d$, the diameter $\rm{diam}_i(A)$ of $A$ in the direction $n_{H_i}$ as the width of the maximal  parallel slab containing $A$.  
If $\bigcup_{j = 1}^{J'} S'_j$ is connected and if it meets $H_{i}[0, \delta]$,
then by the diameter bound  ${\text{diam}}_i (S'_j ((3^{k_1} \delta,..., 3^{k_d} \delta)))  \leq c3^{k_i} \delta$, valid for cone-sets, cone-cylinder sets and cylinder sets, we have
\be \label{diam}  \rm{diam}_i(\bigcup_{j =1}^{J'} S'_j )\le c\delta(1 + 3 + ... + 3^{J -1} ) = c'\delta 3^{J'}.\ee
\vskip.3cm


\noindent{\bf 3.4. Dependency graphs}. The above subsection describes  a collection of   supersets $S_j'$  generated by
the constituent  $M$-regions in ${\cal M}_K(\0, \delta)$ . These sets
are either cone sets, cone-cylinder sets, or cylinder sets, depending on whether the $M$-region lies entirely in
$[0, (T^*)^{1/d}]^d$,  meets the boundary of $[0, (T^*)^{1/d}]^d$, or  lies outside $[0, (T^*)^{1/d}]^d$.
 Given any vertex ${\mathscr V}_i \in
{\cal V}_K \setminus \{ \0\}$ we may likewise construct a collection ${\cal M}_K({\mathscr V}_i, \delta)$ of dyadic $M$-regions in
$p_d({\mathscr V}_i, 1/2)$ and use them to generate  a corresponding collection of $S_j'$ regions.
Here $p_d({\mathscr V}_i, \delta)$ is the parallelepiped defined before \eqref{Wl} and without loss of generality, we may assume that
the parallelepipeds $p_d({\mathscr V}_i, 1/2)$ are disjoint.
We embed
the union of ${\cal M}_K({\mathscr V}_i, \delta), i \leq f_0(K),$ into a  (possibly not unique) larger  collection
${\cal M}_K(m(T, \delta))$ of  $M$-regions having cardinality $m(T,\delta)$ and  which is maximal for the entire surface
$K(v = T)$.  This is possible since among all possible collections containing the union of
$\M_K({\mathscr V}_i,\delta), 1 \leq i \leq f_0(K),$ there is at least one which is maximal.

The integer $m(T, \delta)$ may not be unique, but
in any case it is bounded above and below by integers depending only on $T$, as shown in \cite{BR2}.
Next, let $ {\cal S}'(\delta) := \{S'_j \}_{j = 1}^{m(T, \delta)}$
be the $S'_j$  regions generated by the $M$-regions in ${\cal M}_K(m(T, \delta))$.  The additional $S'_j$ regions which are not associated with
a dyadic $M$-region  are defined exactly as in \cite{BR2}.  The collection ${\cal S}'(\delta)$ partitions the annulus $\A(s, T^*,K)$.
Notice that $m(T, \delta)$ plays the role of $m_\eta:= m(T_\eta)$ in \cite{BR2}.  The choice of a saturated system on $K(v = s)$
is not relevant for our discussion and may be chosen as in \cite{BR2}.

Now we are ready to consider a dependency graph $\cal G := ({\cal V}_{\cal G} , {\cal E}_{\cal G} )$, where ${\cal V}_{\cal G} := {\cal S}'(\delta).$
Following section 6 of \cite{BR2},
define  $L_j,  1 \leq j \leq m(T,\delta),$ to be the union of all $S'_k \in {\cal S}'(\delta)$ such that there are points
$$
a \in S'_j \cap K(v \geq s), \ \ b \in S'_k \cap K(v \geq s)
$$
with the segment $[a,b]$ disjoint from $K(v \geq T^*).$  $L_j$ is not empty since it contains $S'_j$.  Also, $S'_k \subset L_j$
if and only if $S'_j \subset L_k$.  Join vertices $i, j \in {\cal V}_{\cal G}$ with an edge iff $L_i$ and $L_j$ contain at least one
$S'_k$ in common.  Let ${\cal E}_{\cal G}$ be the edges thus defined.

The  first assertion of the next result is proved in  the Appendix.  The second assertion
is Lemma 6.1 of \cite{BR}.  Let $L(\la):= T(K)^{3} (\log \log \la)^{3(d-1)},$  where $T(K)$ is the number of flags of $K$.
Recall that we implicitly assume
 $\delta \in
(0,1/2)$
satisfies $\log_3 T/\delta^d \in \Z$.

\newpage

\begin{lemm} \label{Lem3.2} For any fixed   $\delta \in
(0, 1/2)$ we have
\vskip.1cm
\noindent a.  The geometric properties of sets in ${\cal S}'(\delta)$  fulfill the requirements of [\cite{BR2}, p. 1518, 5 lines before (5.4)], and
\vskip.1cm
\noindent b.   There is a constant $c^* \in (0, \infty)$ such that
for all $1 \leq j \leq m(T, \delta)$ and all $\la \in [1, \infty)$, we have $ {\rm{card}} \{ k: \ S'_k \subset L_j  \} \leq c^*L(\la)$.
\end{lemm}
In other words $L_j, 1 \leq j \leq m(T, \delta),$ contains at most $c^* L(\la)$ sets in ${\cal S}'(\delta)$.  As shown in
\cite{BR2}, neither $L(\la)$ nor the maximal degree of $\cal G$ is a function of $\delta$.
 By \eqref{diam}, if
 $L_j$ has non-empty intersection with $H_l[0, \delta]$, then there exists a constant $c_{\text{diam}} \in (0, \infty)$ such that
its diameter in the direction $n_{H_i}$ satisfies
\be\label{diamLj}
\rm{diam}_i(L_j)\le c_{\text{diam} } \delta 3^{c^*L(\la)}.
\ee

\begin{lemm} \label{indep1}
If $S'_j \subset [\0, \delta]^d$ and if $S'_i \cap H_l[ c_{\text{diam} } \delta 3^{c^*L(\la) + 1 },  {\emph{diam}} (K)] \neq \emptyset $
for all $1 \leq l \leq d$, then there is no edge in ${\cal E}_{\cal G}$ between $j$ and $i$.
\end{lemm}

\noindent{\em Proof.} Case (i). If there were such an $S'_i$, and if such an  $S'_i$ were generated by an $M$-region in ${\cal M}_K(\0, \delta)$,
then by the diameter bound \eqref{diamLj},  it would follow that $L_i \subset H_l[2 c_{\text{diam} }\delta 3^{ c^* L(\la)}, \diam(K)]$
whereas $L_j \subset H_l[ 0,  c_{\text{diam} } \delta 3^{c^* L(\la)} )]$.
 Thus  $L_i  \cap  L_j  = \emptyset$, showing that there is no edge
between $j$ and $i$ in the case when $S'_i$ is produced by an $M$-region in ${\cal M}_K(\0,\delta)$. Case (ii).  When
$S'_i$ is generated by an $M$-region in ${\cal M}_K(m(T,\delta)) \setminus {\cal M}_K(\0, \delta)$ then we proceed by contradiction.
If there were an edge between $j$ and $i$, then there would exist an
$S'_l$ (on the path between $S'_i$ and $S'_j$) such that $S'_l \cap H_l[ 2 c_{\text{diam} }\delta 3^{c^* L(\la)},  {\text{diam}} (K)] \neq \emptyset$, $S'_l$ is
generated by an M-region in ${\cal M}_K(\0, \delta),$ and there is an edge between $l$ and $i$. This contradicts the first case of this proof.
\qed

\vskip.5cm

Next we recall Lemma 7.1 of \cite{BR2} and the discussion on pages 1522-23 of \cite{BR}. Though this lemma is proved
in \cite{BR2} for the volume score $\xi_V$, its proof is general and applies to the scores $\xi_k$ as well.
This result provides conditions for independence of scores on disjoint sets with respect to the graph distance between the
sets.


\begin{lemm} \label{3.4} Let $\xi \in \Xi$ and let ${\cal W}_1$ and ${\cal W}_2$ be disjoint subsets of $\cal V_{\cal G}$ having no edge
between them.  Conditional on $A_\la$, the random variables
$$
\sum_{x \in \P_\la \cap (\cup_{i \in {\cal W}_1 } S'_i) } \xi(x, \P_\la)  \  \text{and} \
 \  \sum_{x \in \P_\la \cap (\cup_{i \in {\cal W}_2 } S'_i) } \xi(x, \P_\la)
$$
are independent.
\end{lemm}

\vskip.3cm

The next  result provides conditions for independence of scores on disjoint sets with respect to the Euclidean distance between the sets.

\begin{lemm} \label{4.4} There exists a constant $c' \in (0, \infty)$ such that if $S'_0 \in {\cal S}'(\delta)$ is a subset of  $[0, \delta]^d$ and if $S' \in {\cal S}'(\delta)$   is at Euclidean distance
at least $c'\delta 3^{ c^*L(\la)}$ from $[0, \delta]^d$, then conditional on $A_\la$ the sum of the scores on $S'_0$ and $S'$ are independent.
\end{lemm}

\noindent{\em Proof.}  This is a consequence of Lemma \ref{indep1}. We can choose $c'>0$ such
 that if the Euclidean distance between $S'_0$ and $S'$ exceeds $c'\delta 3^{c^* L(\la)}$ then the distance in any direction $n_{H_i}$  is greater than $c_{\text{diam} }\delta 3^{c^*L(\la)+1}$. Consequently, by
\eqref{diamLj}, the graph distance exceeds $2$,
 because any edge from the dependency graph would intersect more than
 $c^*L(\la)$ cylinder sets $S'_j$.  \qed

\vskip.3cm

Recall the definition of $Z_i(\delta), 1 \leq i \leq f_0(K),$ at \eqref{Wl}.

\begin{lemm} \label{L5a}  Conditional on $A_\la$, the random variables  $Z_i(\delta),  1 \leq i  \leq f_0(K),$ are independent for $\lambda$ large enough whenever
$\delta:= \delta(\la)$ satisfies $\delta 3^{c^*L(\la)} = o(1).$
\end{lemm}

\noindent{\em Proof.}  Let $S({\mathscr V}_i) \in {\cal S}'(\delta)$ be a subset of  $p_d({\mathscr V}_i, \delta) \cap K, 1 \leq i \leq f_0(K).$
The Euclidean distance
between  $S({\mathscr V}_i)$ and $S({\mathscr V}_j), i \neq j,$ is bounded below  by $||{\mathscr V}_i - {\mathscr V}_j||- 2 \delta$  which exceeds
$c' \delta 3^{ c^* L(\la)}$.  Now apply Lemma \ref{4.4}. \qed

\vskip.5cm

\noindent{\bf 3.5. Variance is additive over vertices of $K$}.
Put $\A(s, T^*, K, \delta) := \A(s, T^*, K) \setminus \bigcup_{i=1}^{ f_0(K)} p_d({\mathscr V}_i, \delta)$ and set $$\P_\la(s, T^*,K, \delta):= \P_\la  \cap {\cal A}(s, T^*,K, \delta).$$  Recall the definition of $Z$ and $Z_i$ at \eqref{defZ} and \eqref{Wl}. Conditional on $A_\la$,  we have for all $\xi \in \Xi$
\be \label{Zsum}
Z = Z_0 + \sum_{i=1}^{ f_0(K)} Z_i,
\ee
where
\be \label{Z0} Z_0:= Z_0(\delta):= \sum_{x \in \P_\la(s, T^*,K, \delta) } \xi(x, \P_\la)
\ee
is the contribution to $Z$ from points in $\P_\la$ which are far from ${\cal V}_K$.

Recall from \eqref{defr0} that $\delta_0 := \exp(- (\log \la)^{1/d})$.  
 We now  put $\delta$ to be $\delta_1:= r(\la, d) \delta_0$,  where
$ r(\la, d) \in [1, 3^{1/d} )$ is chosen so that $\log_3(T/\delta_1^d)\in \Z$.

Roughly speaking, conditional on $A_\la$, we may bound the  number of points in
$\P_\la(s, T^*,K, \delta_1)$ as well as  magnitudes of scores arising from such points.  In this way, 
 the next lemma shows that the contribution to the total score arising from $Z_0(\delta_1)$ is negligible. It also shows that the difference between the sum of the volume scores and the defect volume of $K_\la$ is negligible. The proof is in the Appendix.
\begin{lemm} \label{L4} For $\xi \in \Xi$  we have
\vskip.1cm
\noindent a.
$\var[Z_0(\delta_1) | A_\la] = o(\Var[Z])$.
\vskip.1cm
\noindent b. $\var\left[\frac{1}{\lambda}\sum_{x\in\P_\la}\xi_V(x,\P_\la)\right]=\var[\Vol(K_\la)]+o(\var[\Vol(K_\la)])$.
\end{lemm}

\vskip.3cm

The next lemma, also proved in the Appendix,  shows that the event $A_\la^c$ contributes a negligible amount to the first and second order statistics of $Z$ and $Z_i, 1 \leq i  \leq f_0(K).$


\begin{lemm} \label{Lem9}   Let $Z_i:= Z_i(\delta_0)$ be as at \eqref{Wl}.
We have uniformly for $1 \leq i  \leq f_0(K)$:
$$ \max\{ |\E[Z]-\E[Z | A_\la ]| , \ |\E[Z_i]-\E[Z_i | A_\la ]|,  \ |\E[Z_i]-\E[Z_i {\bf{1}}( A_\la) ]| \} =o(\E[Z]),
$$
and
\begin{align*}
& \max\{  |\Var[Z]-\Var[Z {\bf{1}}( A_\la) ]| , \ |\Var[Z]-\Var[Z | A_\la ]| ,\\
&\hspace*{2cm} \ |\Var[Z_i]-\Var[Z_i | A_\la ]|,  \ |\Var[Z_i]-\Var[Z_i {\bf{1}}( A_\la) ]| \} =o(\Var[Z]).
  \end{align*}
\end{lemm}

\vskip.3cm

Finally we may prove the second  main result of this section.

\vskip.3cm

\noindent{\em Proof of Proposition \ref{Prop2}.}
By Lemma \ref{Lem9} it suffices to show
\be \label{decomp}
\Var [Z | A_\la ] =  \sum_{i=1}^{f_0(K)}  \Var [Z_{i}(\delta_0)| A_\la] +  o( \Var[Z] ).
\ee
To do so, we proceed in two steps: (i) we first show
\be \label{propfordelta1}
\Var [Z|A_\la] =  \sum_{i=1}^{f_0(K)}  \Var [Z_{i}(\delta_1)|A_\la] +  o( \Var[Z] ),
\ee
and (ii) then  for every $1\le i\le f_0(K)$ we show
\be \label{equiv01}
\Var [Z_{i}(\delta_1)|A_\la]=\Var [Z_{i}(\delta_0)|A_\la]+o(\Var[Z]).
\ee
Let us show \eqref{propfordelta1}. Let $\cov((X,Y) | A_\la)$ be short for $\E[(X - \E [X | A_\la] ) (Y - \E [Y | A_\la] )  | A_\la]$.   Recalling \eqref{Zsum},
we have
\begin{align*}\label{sum1}
\Var [Z | A_\la ] &  =  \Var [ Z_0(\delta_1)  +  \sum_i Z_i(\delta_1) | A_\la ] \\
& = \Var [Z_0(\delta_1)| A_\la ] + \Var \left[\sum_i Z_i(\delta_1)| A_\la \right]
+ 2 \cov \left(  (\sum_i Z_i(\delta_1) , Z_0(\delta_1)) | A_\la  \right)
\\
& =  \sum_i \Var[Z_i(\delta_1)| A_\la] + 2 \cov \left(  (\sum_i Z_i(\delta_1) , Z_0(\delta_1)) | A_\la  \right) + o(\Var[Z]),
\end{align*}
where the last equality follows from $\delta_1 3^{c^* L(\la) } = o(1),$  the conditional independence of $Z_i(\delta_1), 1 \leq i  \leq f_0(K)$,
as given by Lemma \ref{L5a}, as well as  Lemma \ref{L4}.

If  random variables $X$ and $Y$ satisfy $\max\{ \Var[X + Y | E], \Var[Y|E]\} = O(\Var[X + Y])$, then
writing $X = (X + Y) - Y$, it follows that $\Var[X|E] = O(\Var[X + Y])$.  We have
$\max\{ \Var[ \sum Z_i(\delta_1) + Z_0(\delta_1) | A_\la ], \ \Var[Z_0(\delta_1) |A_\la] \} = O(\Var [Z])$ by \cite{BR2} and
by Lemma \ref{L4}. It follows that $\Var[\sum Z_i(\delta_1)|A_\la] = O(\Var[Z])$. This estimate, Lemma \ref{L4}
again, and the Cauchy-Schwarz inequality give
\begin{align*}
\cov \left(  (\sum_i Z_i(\delta_1) , Z_0(\delta_1))   | A_\la \right) & \leq  \sqrt{ \Var [ \sum_i Z_i(\delta_1) | A_\la]}
\cdot \sqrt{ \Var[ Z_0(\delta_1) | A_\la ] } \\
& = O(\sqrt{ \Var[ Z] } ) o(\sqrt{ \Var[ Z] }) \\
& = o(\Var[Z]).
\end{align*}
This yields the  decomposition \eqref{propfordelta1}.

To prove \eqref{equiv01} we introduce $\delta_1':=r'(\la,d)\delta_0$  where
 $r'(\la, d) \in (3^{-1/d},1]$
is chosen so that $\log_3(T/\delta_1'^d)\in \Z$.
Methods similar to the proof of Lemma \ref{L4} show that
 $\Var[Z_0(\delta_1')|A_\la]=o(\Var[Z])$ and  $\Var[\sum_{x\in {\mathcal B}}\xi(x,\P_\la)|A_\la]=o(\Var[Z]),$ with ${\mathcal B}$ a subset of $\P_{\la}(s,T^*,K,\delta_1')$. 
  Note that ${\mathcal B}_i:= (p_d({\mathscr V}_i,\delta_1)\setminus p_d({\mathscr V}_i,\delta_0))  \cap \P_\la,
  1\le i\le f_0(K),$ are subsets of $\P_{\la}(s,T^*,K,\delta_1')$.
  Consequently,
 $$\Var[Z_i(\delta_1)-Z_i(\delta_0)|A_\la]=o(\Var[Z]).$$
 Moreover, $$\Cov((Z_i(\delta_1)-Z_i(\delta_0),Z_i(\delta_0))|A_\la)\le \sqrt{\Var[Z_i(\delta_1)-Z_i(\delta_0)|A_\la]}\sqrt{\Var[Z_i(\delta_0)|A_{\la}]}=o(\Var[Z]).$$
 We deduce \eqref{equiv01} from the two previous equalities. This completes the proof of Proposition \ref{Prop2}.  \qed

\section{Re-scaled  convex hull boundaries,  $k$-face, and volume functionals} \label{transforms}

\allco
Section \ref{LEMMAS} showed that variance asymptotics for $f_k(K_\la)$ and $\Vol(K_\la)$ are determined by the
respective behavior
of $\xi_k$ and $\xi_V$  on $\P_\la \cap Q_0$.
We  discuss scaling transforms of $\P_\la \cap Q_0$,  $\partial K_\la \cap Q_0,$ as well as
transforms for $\xi_k$ and $\xi_V$ restricted to input $\P_\la \cap Q_0$.  \\

\noindent{\bf 4.1. Parallel between the scaling transform $T^{(\la)}$ and those in previous work.}
Scaling transforms lie at the heart of our asymptotic analysis.  Before discussing the technical details, we explain their relevant geometric
aspects, comparing $T^{(\la)}$ with counterparts in previous works
on Gaussian polytopes \cite{CY2}, as well as
random polytopes in the unit ball \cite{CSY, SY} and in smooth convex bodies \cite{CY}. \\

\noindent{\it Floating bodies and associated coordinates}. Seminal works of  B\'arany  and Larman \cite{BL89} and  B\'ar\'any \cite{Bar89} show
 the importance of the deterministic approximation of the random polytope inside the mother body $K$ by a floating body $K(v\ge 1/\la)$. Consequently, it makes sense to use the parametric surfaces $\partial K(v\ge t/\la)$, $t>0,$  to associate to any point $z\in K$ a {\it depth coordinate} which is the specific $t$ such that $z\in \partial K(v\ge t/\la)$ and a {\it spatial coordinate} indicating the position of $z$ on the surface $\partial K(v\ge t/\la)$.
When $K$ is the unit ball, the floating bodies are balls $B(\0,r)$, $0<r<1$, and  coordinates coincide with  the usual spherical coordinates. When $K:= (0,\infty)^d$, the floating bodies are pseudo-hyperboloids, as seen in the next subsection.
  We could call the associated coordinates {\it cubical coordinates}. In the case of a general convex
mother body $K$,  there is not necessarily a natural way of globally defining  a spatial coordinate, which explains {\it a posteriori} why we dealt with local spherical coordinates in \cite{CY}.\\

\noindent{\it Extreme points and duality}. This paper, as well as \cite{CY,CY2, CSY}, rely on a dual characterization of extreme points.
 Arguably, it is most natural to define an extreme point as {\it  a point from the input on the boundary of the convex hull}. By duality, we
  may also assert that a point $z_0$ from the input is extreme if it is included in a support hyperplane of the convex hull. In most cases, any hyperplane containing a fixed point $z$ from the input is tangent to exactly one surface $\partial K(v\ge t/\la)$ at one point of tangency. This suggests the idea of considering the {\it petal} of $z$, i.e. the subset $S(z)$ of $K$ whose boundary  $\partial S(z)$ consists of
  all points of tangency of hyperplanes containing $z$. In the case of the unit ball, and when the origin is inside the convex hull, the petal of $z$ is the ball with diameter $[0,z]$. The collection of such balls associated to the points of the input constitutes the so-called {\it Voronoi flower} of the input with respect to the origin, which explains {\it a posteriori} the appellation  {\it petal}. In the case of the orthant $(0,\infty)^d$, when the point $z$ is cone-extreme (recall Definition \ref{cone-ext}), its petal is defined in \eqref{eq:eqpetal} below. This provides the second definition of an extreme point ({\it cone-extreme} in the case of the orthant): {\it A point from the input is extreme iff its petal is not covered by the petals from the other points from the input}.\\

\noindent{\it Scaling transformations}. As in \cite{CY2, CSY}, this paper uses the set of coordinates induced by the floating bodies to build the scaling transformation. With the proper re-scaling of both the spatial and depth coordinates, we get a new picture in a product space $\R^{d-1}\times \R$, the height being the re-scaled depth coordinate. The duality of the two definitions of the extreme points (or cone-extreme points in the case of the orthant) is even more apparent in the re-scaled picture. Indeed, the  re-scaled random polytope
may be described either via the re-scaled boundary of the convex hull given below by $\partial \Phi(\P^{(\la)})$ at \eqref{PhiDown} or
via the re-scaled boundary of the union of petals given by $\partial \Psi(\P^{(\la)})$ at \eqref{PsiUp}.


\vskip.3cm

\noindent{\bf 4.2. A new characterization of cone-extreme points.}
We consider surfaces
\be \label{eq:hyperboloid}
{\mathcal H}_t:=\{(z_1,\cdots,z_d)\in (0,\infty)^d : \ \prod_{i=1}^d z_i=t\}, \ \ t > 0.
\ee

When $d=2$ each  ${\mathcal H}_t$ is a branch of a hyperbola and for this reason we will sometimes refer to ${\mathcal H}_t$ as a pseudo-hyperboloid. The surfaces ${\mathcal H}_t, t > 0,$ coincide with boundaries of floating bodies of the orthant $[0,\infty)^d$, as shown in Lemma \ref{hyper}, and play a key role in the description of cone-extreme points of  input inside $(0,\infty)^d$.

For every $z^{(0)}\in (0,\infty)^d$, we denote by $H(z^{(0)})$ the hyperplane  tangent to the unique surface ${\mathcal H}_t$ containing $z^{(0)}$.
The gradient of the function $f(z) = \Pi_{i=1}^d z_i$ at $z^{(0)}$ is $t {1 \over z^{(0)}}$, where ${1 \over z^{(0)}} := ( {1 \over z_1^{(0)}},...,
{1 \over z_d^{(0)}}).$  It follows that $H(z^{(0)})$ is described by
\begin{equation}
  \label{eq:equationhyperplantangent}
 H(z^{(0)}):= \{(z_1,...,z_d) \in \R^d: \ \sum_{i=1}^d\frac{z_i}{z_i^{(0)}}=d\}.
\end{equation}

To every point $z^{(0)}\in (0,\infty)^d$, we associate the  surface  $${\mathcal S}(z^{(0)}):= \{z \in (0, \infty)^d:  \ z^{(0)} \in H(z) \}.$$
%
The {\em petal} of $z^{(0)}$ is the closed set ${\mathcal S}^-(z^{(0)})$ of points above ${\mathcal S}(z^{(0)})$.
Notice that ${\mathcal S}^-(z^{(0)})$ is also the set of points $z$ such that $z^{(0)}$ lies `below' $H(z)$.
Using \eqref{eq:equationhyperplantangent}, we have 
\begin{equation}
  \label{eq:eqpetal}
   {\mathcal S}^-(z^{(0)}):=  \{(z_1,...,z_d) \in (0,\infty)^d: \ \sum_{i=1}^d\frac{z_i^{(0)}}{z_i}\le d\}.
\end{equation}

\begin{figure}
  \centering
\includegraphics[trim= 0cm 15cm 0cm 407cm, clip, scale=0.5]{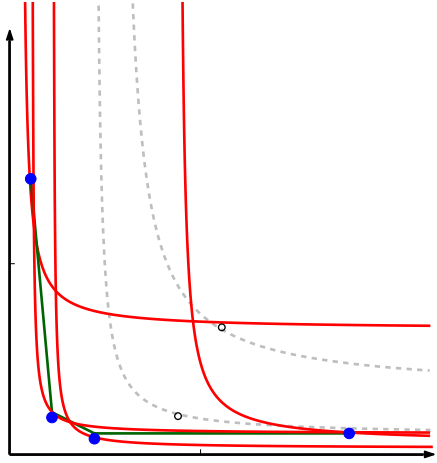}
\caption{The point process $\mbox{Ext}(\P_\la)$ (blue);  the boundary of the associated petals containing the extreme points (red); the boundary of the convex hull (green). Points which are not extreme
are apices of gray petals.}
  \label{fig:petals_polytope}
\end{figure}

The next lemma characterizes cone-extreme points in terms of the geometry of petals (see Figure \ref{fig:petals_polytope}).  We are not
aware of an analogous characterization of extreme points which are not cone-extreme.

\begin{lemm}\label{resultpetal}
Let $\X$ be any point set in $(0,\infty)^d$. Then  $z^{(0)}\in \X$ is cone-extreme
with respect to  $\X$ if and only if ${\mathcal S}^-(z^{(0)})$ is not completely covered by $\bigcup_{z\in \X \setminus\{z^{(0)}\}}{\mathcal S}^-(z)$.
\end{lemm}
\noindent{\em Proof.}
Indeed, ${\mathcal S}^-(z^{(0)})$ is not covered iff there exists $z\in {\mathcal S}(z^{(0)})$ which is below each of the surfaces ${\mathcal S}(z^{(0)})$.  This  is equivalent to saying that the hyperplane $H(\blue{z})$ is a support hyperplane of $\mbox{co}(\X)$ containing $z^{(0)}$ and with outward normal in $C_{\0}(K)$. \qed
\vskip.3cm

\noindent{\bf 4.3. The scaling transform $T^{(\lambda)}$ at \eqref{defT}.}
Here we  describe the image under $T^{(\lambda)}$ of  $\P_\la \cap Q_0$, the image of half-spaces
 with bounding hyperplane  $H(z^{(0)})$, as well as the image of petals.  As a by-product, we find  the image of a face of
 $K_\la \cap Q_0$.

Recall the definitions of $V, l_i(v), l(v),$ and $p_V:\R^d\to V$ introduced before \eqref{defG} and after \eqref{defT}, respectively.
For any $x=(x_1,\cdots,x_d)\in \R^d$ and any function $f: \R \to \R$, we define the vector $f(x):=(f(x_1),\cdots,f(x_d))$.
Note that $e^{l}: V \to {\mathcal H}_1.$

Recall that $w = (v,h)$ denotes a generic point in $V \times \R$.  The inverse of $T^{(\lambda)}$ is
\begin{equation}
  \label{eq:invscaltrans}
[T^{(\lambda)}]^{-1}:\left\{\begin{array}{lll}V\times\R&\longrightarrow &(0,\infty)^d\\(v,h)&\longmapsto &\lambda^{-1/d}e^{h}e^{l(v)}\end{array}\right..
\end{equation}
The expression for $[T^{(\lambda)}]^{-1}$ is justified as follows.
 We have $l(p_V(\log(z)))=\log z -\langle\log z, u^0\rangle\frac{1}{d}u^0$ where $u^0:=(1,\cdots,1)$ and $\langle \cdot,\cdot\rangle$ is the usual scalar product on $\R^d$. Indeed, $\langle\log z,u^0\rangle\frac{1}{d}u^0$ is the projection of $\log z$ onto the line directed by $u^0$. So $\exp(l(p_V(\log z)))= \exp(-\frac{1}{d}\sum_{i=1}^d\log z_i) \cdot z=(z_1.....z_d)^{-1/d} z$ and so
$[T^{(\lambda)}]^{-1} (T^{(\la)} )$ is the identity, as desired.

\begin{defn} For all $\la \in [1, \infty)$ we put
\begin{equation}
  \label{eq:Wlambda}
W_\la := T^{(\la)}(Q_0)=\{(v,h) \in \R^{d-1} \times \R:
\ h \leq - l_i(v) + \log \la^{1/d} \delta_0, \ 1 \leq i \leq d \}
\end{equation}
and
\be \label{TPQ} \P^{(\la)}:= T^{(\la)}(\P_\la \cap Q_0).\ee
When $\la = \infty$ we identify $W_\la$ with $V \times \R$ and $\P^{(\la)}$ with  $\P$
at  \eqref{defP}.
\end{defn}

\begin{lemm}\label{imagePPP}
Let $\P$ be the Poisson point process  at  \eqref{defP}.
Then $\P^{(\la)}$ is equal in distribution to $\P\cap W_\la$ and $
\P^{(\la)} \tod \P$  as $\la \to \infty$.
\end{lemm}

\vskip.3cm
\noindent{\em Remark}.  Since   $d\P^{(\la)}$ is the image under $T^{(\la)}$ of $\la \Vol_d$, with $\Vol_d$ the $d$-dimensional volume measure
on $\R^d$,  this lemma says that  $T^{(\la)} (\la \Vol_d) \tod d\P$.

\vskip.3cm
\noindent{\em Proof.} $\P^{(\la)}$ is a Poisson point process with intensity measure  $T^{(\lambda)}(\lambda dz)$. Endow  $V$ with a direct orthonormal basis ${\mathcal B}$.  Using  \eqref{eq:invscaltrans},  the Jacobian of $[T^{(\lambda)}]^{-1}$ with respect to the direct orthonormal basis of $\R^d$ given by
 $({\mathcal B},\frac{1}{\sqrt{d}}u^{(0)})$ 
equals
$$\sqrt{d}\lambda^{-1}\exp(dh)\exp(\sum_{i=1}^dl_i(v))D,$$
where $D$ is the determinant of the matrix for the change of basis from $({\mathcal B},\frac{1}{\sqrt{d}}u^{(0)})$ to the standard basis of $\R^d$. Since both bases are direct and orthornormal, we have $D=1$. Moreover, we notice that $\sum_{i=1}^dl_i(v)=0$ because $v\in V$.

Consequently, $\P^{(\la)}$ has intensity measure with no $\la$ dependency save that it is carried by the `pyramid-like' set $W_\la:= T^{(\la)}(Q_0)$. In other words,
\begin{equation}
  \label{eq:intenPlambda}
d\P^{(\lambda)}=\sqrt{d}\exp(dh){\bf 1}((v,h) \in W_\la)dv dh.
\end{equation}
 Lemma \ref{imagePPP} follows from  \eqref{eq:intenPlambda} and the convergence  $W_\la \uparrow W=V\times \R$.\qed
~\\~\\


Having considered the behavior of the scaling transform $T^{(\lambda)}$ on $\P_\la$, we now consider the image under $T^{(\lambda)}$ of
surfaces ${\mathcal H}_{c/\lambda}$ with $c>0$, $C_{\0}$-half-spaces, and petals.  A $C_{\0}$-half-space is one having an outward normal in $C_{\0}(K)$,
where $C_{\0}(K)$ is as in Definition \ref{cone-ext}.   Note that
a $C_{\0}$-half-space is bounded by a hyperplane $H(z^{(0)})$ for some $z^{(0)}\in (0,\infty)^d$.  Denote by $H^+(z^{(0)})$ the half-space bounded by $H(z^{(0)})$ and containing $\0$.

Recall the definition of the down cone-like grain $\Pi^{\downarrow}$ given at \eqref{Cone} and its translate $\Pi^{\downarrow}(w):=w\oplus \Pi^{\downarrow}$, $w=(v,h)\in\R^{d-1}\times\R$. Define similarly the up cone-like grain
  \begin{equation}
    \label{eq:defupconegrain}
    \Pi^{\uparrow}:= \{(v,h) \in \R^{d-1} \times \R: \ h \geq G(-v) \}, \ \
  \end{equation}
and the translate $\Pi^{\uparrow}(w)=w\oplus \Pi^{\uparrow}$, $w \in \R^{d-1}\times \R$.

The duality between up and down cone-like grains is expressed through the following equivalence: for all $w,w'\in \R^{d-1}\times\R$,
$$w\in \Pi^{\uparrow}(w')\Longleftrightarrow w'\in \Pi^{\downarrow}(w).$$
The next lemma shows that $T^{(\lambda)}$  sends pseudo-hyperboloids to hyperplanes parallel to $V$,  $C_{\0}$-half-spaces to down cone-like grains,
 and petals to up cone-like grains.
\begin{lemm}  \label{lempetal}
(i) For every $c \in (0, \infty)$, we have
$$T^{(\lambda)}({\mathcal H}_{c/\lambda})=V\times \left\{\frac{1}{d}\log(c)\right\}.$$
(ii) For every $C_{\0}$-half-space  $H^+(z^{(0)})$, $z^{(0)}\in (0,\infty)^d$, we have
$$T^{(\lambda)}(H^+(z^{(0)}))=\Pi^{\downarrow}(T^{(\lambda)}(z^{(0)})).$$
(iii) For every petal $S^-(z^{(0)})$, $z^{(0)}\in (0,\infty)^d$,  we have
$$T^{(\lambda)}(S^-(z^{(0)}))=\Pi^{\uparrow}(T^{(\lambda)}(z^{(0)})).$$
\end{lemm}
\noindent{\em Proof.} For every $(v,h)\in V\times\R$, we have by \eqref{eq:invscaltrans}
\begin{align*}
[T^{(\lambda)}]^{-1}(v,h)\in {\mathcal H}_{c/\lambda}&\Longleftrightarrow \prod_{i=1}^d (\lambda^{-1/d}e^he^{l_i(v)})=c/\la\Longleftrightarrow e^{dh}=c
\end{align*}
which shows (i).

Fix $z^{(0)}\in (0,\infty)^d$ and put   $T^{(\lambda)}(z^{(0)}):=(v^{(0)},h^{(0)})$.
We notice that
\begin{equation}\label{imageinverse}
T^{(\lambda)}(\frac{1}{z^{(0)}})=(-v^{(0)},-h^{(0)}+2\log(\lambda^{1/d})).
\end{equation}
Using the equation of $H^+(z^{(0)})$ implied by \eqref{eq:equationhyperplantangent}, the formula for $[T^{(\lambda)}]^{-1}$ at \eqref{eq:invscaltrans} and \eqref{imageinverse}, we have for any $(v,h)\in V\times \R$
\begin{align*}
[T^{(\lambda)}]^{-1}((v,h))\in H^+(z^{(0)})&\Longleftrightarrow  \langle[T^{(\lambda)}]^{-1}((v,h),\frac{1}{z^{(0)}}\rangle\le d \\
&\Longleftrightarrow \la^{-1/d}e^{h}\langle \exp(l(v)),\frac{1}{z^{(0)}}\rangle\le d \\
&\Longleftrightarrow e^{h-h^{(0)}}\langle \exp(l(v)),\exp(l(-v^{(0)}))\rangle\le d\\
&\Longleftrightarrow h\le h^{(0)}-\log (\frac{1}{d}\langle \exp(l(v-v^{(0)})),u^{(0)}\rangle).
\end{align*}
This last equivalence, coupled with the definition of G at \eqref{defG}, gives (ii).
Similarly, (iii) is a consequence of the equation of the petal ${\mathcal S}^-(z^{(0)})$ at \eqref{eq:eqpetal}, \eqref{eq:invscaltrans} and \eqref{imageinverse}. Indeed, we have for every $(v,h)\in V\times \R$
  \begin{align*}
    [T^{(\lambda)}]^{-1}((v,h))\in {\mathcal S}^-(z^{(0)})&\Longleftrightarrow \langle \lambda^{1/d}e^{-h}e^{-l(v)},\lambda^{-1/d}e^{h^{(0)}}e^{l(v^{(0)})}\rangle \le d\\
&\Longleftrightarrow  e^{h^{(0)}-h}\langle e^{l(v^{(0)}-v)},u^0\rangle \le d\\
&\Longleftrightarrow h\ge h^{(0)}+G(v^{(0)}-v).
  \end{align*}
This completes the proof of Lemma \ref{lempetal}. \qed

\vskip.3cm

\noindent{\bf 4.4. Re-scaled extreme points and scores.}
It is time to define re-scaled scores $\xi^{(\la)}$ on $\P^{(\la)}$.  We use Proposition \ref{Prop1}
to show that on the event $A_\la$ given at \eqref{defA}, the  re-scaled scores $\xi^{(\la)}$ coincide with functionals $\hx^{(\la)}$ defined
in terms of the geometry of the  re-scaled convex hull boundary. This is facilitated with the following definitions.

\begin{figure}
  \centering
  \includegraphics[trim= 7cm 16cm 8.2cm 0cm, clip, scale=0.5]{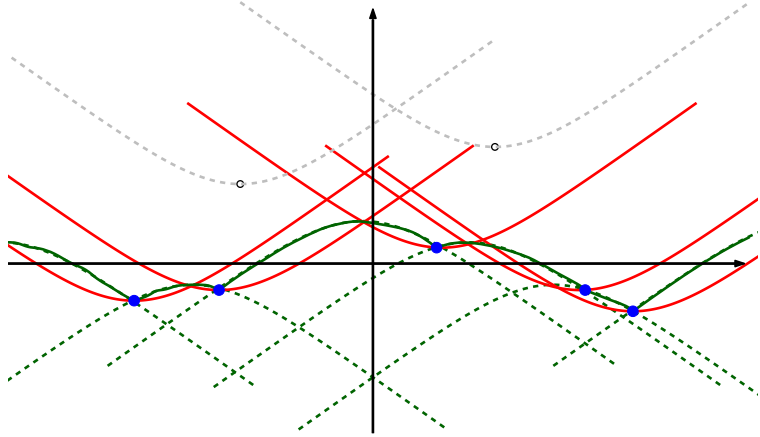}
\caption{The point process $\mbox{Ext}(\P^{(\la)})$ (blue);  the boundary of the up-grains containing the extreme points (red); the boundary  $\partial(\Phi(\P^{(\la)})$ of the down-grains containing $\mbox{Ext}(\P^{(\la)})$ (green). Points which are not extreme
are apices of gray up-grains.}
\label{rescaledpicture2}
\end{figure}

\begin{defn} Write $[\Pi^{\uparrow}( w)]^{(\la)}$ for $\Pi^{\uparrow}( w) \cap W_\la$ and similarly for $[\Pi^{\downarrow}( w)]^{(\la)}$.
Given  $\P^{(\la)}$, $1 \leq \la \leq \infty $, we define the $\P^{(\la)}$-hull
as at \eqref{hull}  with $\X$ set to $\P^{(\la)}$, namely
\be \label{PhiDown}
\Phi(\P^{(\la)}):= \bigcup_{\left\{ \substack{w\in \R^{d-1} \times\R \\
\P^{(\la)} \cap {\rm{int}}(\Pi^{\downarrow}(w)) =\emptyset}\right.} [\Pi^{\downarrow}(w)]^{(\la)}.
\ee
We also put
\be \label{PsiUp}
\Psi(\P^{(\la)}):= \bigcup_{w \in \P^{(\la)}} [\Pi^{\uparrow}(w)]^{(\la)}.
\ee
Abusing notation, we let ${\rm{Ext}}(\P^{(\la)})$ be those points
in $\P^{(\la)}$ which are on the boundary of some down cone-like grain $\Pi^{\downarrow}(w), w \in W_\la$, and
$\text{int}( \Pi^{\downarrow}(w_1)) \cap \P^{(\la)} = \emptyset$.
\end{defn}
Equivalently,
a point $w_0 \in \P^{(\la)}$ is {\em extreme} with respect to $\Psi(\P^{(\la)})$ if
the grain $\Pi^{\uparrow}(w_0)$ is not a subset of the union of the grains $\Pi^{\uparrow}(w), w \in \P^{(\la)} \setminus \{w_0\}$ (see Figure
\ref{rescaledpicture2}).
By Lemma \ref{resultpetal}, Lemma \ref{lempetal}(iii) and Proposition \ref{Prop1}, on the event $A_\la$, the extreme points in $\P_\la \cap Q_0$ are transformed to ${\rm{Ext}}(\P^{(\la)})$.  By  Lemma \ref{lempetal} we also have on the event $A_\la$ that $T^{(\la)}( \partial K_\la \cap Q_0) =
\partial(\Phi(\P^{(\la)})$.



\begin{defn} \label{re-scaled} For $\la \in [1, \infty)$, put $\xi^{(\la)}_k(w, \P^{(\la)}):= \xi_k([T^{(\la)}]^{-1}(w), \P_\la)$. For $w \in {\rm{Ext}}(\P^{(\la)})$
and $\la \in [1, \infty]$ we put
\be \label{hatx}
\hx^{(\la)}_k(w, \P^{(\la)}) :=
(k + 1)^{-1}[ \text{number \ of} \ k\text{-faces  of} \ \Phi(\P^{(\la)})\
\text{containing} \ \ w].
\ee
For  $w \notin {\rm{Ext}}(\P^{(\la)})$ we put $\hx^{(\la)}_k(w, \P^{(\la)}) = 0$.
Similarly,  for $\la \in [1, \infty)$ and  $w \in {\rm{Ext}}(\P^{(\la)})$ we define
$\xi^{(\la)}_V(w, \P^{(\la)}):= \xi_V([T^{(\la)}]^{-1}(w), \P_\la)$
and
\be \label{re-vol}
\hx^{(\la)}_V(w, \P^{(\la)}): =  {1 \over d} \int_{ v\in \rm{Cyl}^{(\lambda)}(w, \P^{(\la)}) }  \int_{-\infty}^{ \partial( \Phi(\P^{(\la)})(v)) }
\sqrt{d} e^{dh} dv dh, \ee
where ${\rm{Cyl}}^{(\lambda)}(w):= {\rm{Cyl}}^{(\lambda)}(w, \P^{(\la)})$ denotes the projection onto $V$ of the hyperfaces of $\Phi(\P^{(\la)})$ containing $w$.
When $w \notin {\rm{Ext}}(\P^{(\la)})$ we define $\hx^{(\la)}_V(w, \P^{(\la)}) = 0$.
\end{defn}
By Proposition \ref{Prop1},  on the event  $A_\la$,
the vertices of $K_\la \cap Q_0$ coincide with cone-extreme points. Putting $w = T^{(\la)}(x)$ gives
\be \label{Key1} \xi_k(x, \P_\la){\bf{1}}(A_\la) =  \hx_k^{(\la)}(w, \P^{(\la)})  {\bf{1}}(A_\la)= \xi_k^{(\la)}(w,\P^{(\la)}){ \bf{1}}(A_\la).
\ee
Lemma \ref{Lem9},  together with \eqref{Key1},
 shows that the  variance asymptotics for $\sum_{x \in \P_\la} \xi_k(x, \P_\la)$  coincide with those for
$\sum_{x \in \P^{(\la)}}  \hx^{(\la)}_k(w, \P^{(\la)}){ \bf{1}}(A_\la)$.
We exploit this fundamental scaling identity in Section 6.
Similarly, by  Proposition \ref{Prop1},
Lemma \ref{resultpetal}, and the remark after Lemma \ref{imagePPP}, we have the analog of \eqref{Key1}, namely
\be \label{Key2} \hx_V^{(\la)}(w, \P^{(\la)}) {\bf{1}}(A_\la)= \xi_V^{(\la)}(w,\P^{(\la)}){ \bf{1}}(A_\la).
\ee
Lemma \ref{Lem9} shows that the  variance asymptotics for  $\sum_{x \in \P^{(\la)}}  \hx^{(\la)}_V(w, \P^{(\la)}){ \bf{1}}(A_\la)$
coincide with those for $\sum_{x \in \P_\la} \xi_V(x, \P_\la).$

Given $\la \in [\blue{1}, \infty),$ let $\Xi^{(\la)}$ denote the collection of re-scaled functionals
$\xi^{(\la)}_{k}, k \in \{0,1,...,d-1\},$ together with
$\xi^{(\la)}_{V}$.  Likewise,  for  $\la \in [1, \infty]$, we let $\hat{\Xi}^{(\la)}$ denote the collection of  functionals
$\hx^{(\la)}_{k}, k \in \{0,1,...,d-1\},$ together with
$\hx^{(\la)}_{V}$.
~\\~\\
\noindent{\bf 4.5. Properties of the function $G$ defined at \eqref{defG}.}
We record two properties of $G$ needed in the sequel.  
Notice that $G$ is an even function only when $d=2$.
\begin{lemm}\label{Gconvex}
$G$ is a positive convex function.
\end{lemm}
\noindent{\em Proof.}
By the convexity of the exponential function, for any $v\in V$, we have
$$G(v)\ge \log(\exp(\frac{1}{d}\sum_{i=1}^d l_i(v)))=\log(1)=0.$$
Let $v,v'\in V$ and $t\in [0,1]$. 
H\"older's inequality gives
\begin{align*}
G(tv+(1-t)v')&=\log(\sum_{i=1}^d[\frac{1}{d}\exp(l_i(v))]^t[\frac{1}{d}\exp(l_i(v')]^{(1-t)})\\
&\le \log([\sum_{i=1}^d \frac{1}{d}\exp(l_i(v))]^t[\sum_{i=1}^d\frac{1}{d}\exp(l_i(v')]^{(1-t)})\\
&=tG(v)+(1-t)G(v').
\end{align*}
Thus $G$ is convex, completing  the proof of Lemma \ref{Gconvex}.\qed
\vskip.3cm
The next lemma shows that the graph of $G$ is sandwiched between circular cones.
\begin{lemm} \label{lemcompare}
There exist $c_1,c_2 \in (0, \infty)$ such that for every $v\in V$,
\be \label{compar}
c_1\|v\|-\log d\le G(v)\le c_2\|v\|.\ee
\end{lemm}
\noindent{\em Proof.}  Since $\max_{1\le i \le d} |l_i(v)|$ is a norm on $V$, it is equivalent to the Euclidean norm $\|\cdot\|$.
It follows that there are constants $c_1,c_2 \in (0, \infty)$  such that $c_1'\|v\|\le \max_{1\le i \le d} |l_i(v)| \le c_2'\|v\|$ for all $v\in V$. We have for every $v\in V$ $$G(v)=\log\left(\frac{1}{d}\sum_{i=1}^d \exp(l_i(v))\right)\le \log\left(\exp(\max_{1\le i \le d} |l_i(v)|) \right)\le c_2'\|v\|.$$
Moreover, one of the $l_i(v)$ is at least equal to $\frac{1}{d-1}\max_{1\le i \le d} |l_i(v)|$. This implies that
$$G(v)\ge \log\left(\frac{1}{d} \exp(\frac{c_1'}{d-1}\|v\|) \right)\ge \frac{c_1'}{d-1}\|v\|-\log d,$$
which establishes \eqref{compar}.\qed

\section{Properties of re-scaled $k$-face and volume functionals}\label{props}
\allco
Section \ref{transforms} introduced re-scaled functionals $\hx^{(\la)}$ of re-scaled input $\P^{(\la)}$. Here we
establish localization properties of the functionals $\hx^{(\la)}  \in \hat{\Xi}^{(\la)}$, bounds on their moments, as well
convergence of their one and two-point correlation functions.
\vskip.3cm

\noindent{\bf 5.1. Stabilization.}  We establish localization of the functionals $\hx^{(\la)}  \in \hat{\Xi}^{(\la)}$ in
both the space and time domains.
Recalling that $B_{d-1}(v,r)$ is the $(d-1)$ dimensional ball
 centered at $v \in \R^{d-1}$ with radius $r$, define the cylinder
 \be \label{cyl}
 C(v,r):= C_{d-1}(v,r):=B_{d-1}(v,r) \times \R. \ee

\vskip.3cm
 We show that  the boundaries of the  germ-grain models  $\Psi(\P^{(\la)})$
and $\Phi(\P^{(\la)})$, $\la \in [1, \infty]$ defined at \eqref{PsiUp} and \eqref{PhiDown}, respectively,  are not far from
$V$.  Recall that $\P^{(\la)}, \la = \infty$, is taken to be $\P$.
If $w \in \rm{Ext}(\P^{(\la)})$ we put $H(w):=H(w,\P^{(\la)})$ to be the maximal height coordinate (with respect to $\R^{d-1}$) of an apex of a down cone-like grain belonging to  $\Phi(\P^{(\la)})$ and containing $w$.  Otherwise, if  $w \notin \rm{Ext}(\P^{(\la)})$ then we put $H(w) = 0$.

\begin{lemm}\label{BLem1}  (a) There is a constant $c$ such that for all $\la \in [1, \infty]$ and  $(v_0, h_0) \in W_\la$
\be \label{fbd-1} P[H((v_0,h_0), \P^{(\la)}) \geq t] \leq c \exp(-\frac {  e^{t/c} } {c}), \ t\ge h_0  \vee 0. \ee
\noindent(b) There is a constant $c$ such that for all $L \in (0, \infty)$  and  $\la \in [1, \infty]$
\be \label{sbd}
 P[ || \partial \Psi(\P^{(\la)}) \cap
C(\0,L) ||_{\infty} >  t] \leq cL^{2(d-1)} \exp(- \frac{t} {c}), \  t \in (0, \infty).
\ee
The bound \eqref{sbd} also holds for $\partial(\Phi(\P^{(\la)})).$
\end{lemm}

\noindent{\em Proof.} We  first prove \eqref{fbd-1}.  Rewrite the event $\{H((v_0,h_0), \P^{(\la)}) \geq  t \}$ as
\begin{align*}
& \ \ \ \ \{H((v_0,h_0), \P^{(\la)}) \geq  t \} \\ & =\{ \exists w_1:=(v_1,h_1) \in \partial[\Pi^{\uparrow}((v_0,h_0))]^{(\la)}:  \ h_1 \in [t,  \infty), [\Pi^{\downarrow}(w_1)]^{(\la)}  \cap \P^{(\la)} = \emptyset \}.
\end{align*}

First consider the case  $\la = \infty$.
Let $w_1:= (v_1,h_1)\in  \partial \Pi^{\uparrow}((v_0,h_0))$ with  $h_1 \in [t, \infty)$.
Recalling \eqref{compar}, the $d\P$ measure of $\Pi^{\downarrow}(w_1)$
is bounded below by the $d\P$ measure of $\{(v,h):h\le h_1-c_2\|v-v_1\|\} \cap  (\R^{d-1} \times [0, \infty))$, which we
generously bound below by
$$c\int_{\frac{h_1}{4}}^{\frac{h_1}{2}}e^{dh}(h_1-h)^{d-1}dh\ge c' e^{h_1/c'}.$$
Here and elsewhere, unless noted otherwise, $c$ and $c'$ denote positive constants which are independent of other parameters
except for dimension and whose value  may change at each occurrence.
Thus the probability that $\Pi^{\downarrow}(w_1)$ does not contain
points in $\P$ is bounded above by $c\exp(-e^{h_1/c})$.

We now discretize $\partial \Pi^{\uparrow}((v_0,h_0))  \cap (\R^{d-1} \times [t, \infty))$.  Notice  that
 if $w_1:=(v_1,h_1) \in \partial\Pi^{\uparrow}((v_0,h_0))$ then \eqref{compar} gives $h_1 =
  h_0 + G(v_0-v_1) \geq h_0 + c_1||v_1 -v_0|| - \log d$, which yields
 $||v_0-v_1|| \leq \frac{1} {c_1} [h_1 - h_0 + \log d]$.
This gives
$$P[  H((v_0,h_0), \P) \geq  t   ]\le
c \int_{t }^{\infty} (h_1 - h_0 + \log d)^{d-2} \exp(- ce^{h_1/c} ) dh_1.
$$
Thus  \eqref{fbd-1} holds.

Next consider the case $\la \in [1, \infty)$. The above argument still holds as soon as we can show for any $(v_1,h_1)\in W_\la$, that the $d\P$ measure of the intersection of $W_\la$ with the down cone-like grain $\{(v,h):h\le h_1-c_2\|v-v_1\|\}$ is
 bounded below by $c\exp(-e^{h_1/c}/c)$. To do so, let $C_{\mbox{\tiny{min}}}$ be the largest circular cone included in the pyramid $W_\la$ and with the same apex as $W_\la$. (Actually, $C_{\mbox{\tiny{min}}}$ does not depend on $\la$ since for $\la'>\la$, $W_{\la'}$ is the image of $W_\la$ by a translation.) Then any {\it vertical} cone with apex in $W_\la$ is such  that its intersection with $W_\la$  is either the cone itself or contains a translate of $C_{\mbox{\tiny{min}}}$ with same apex. Consequently, the intersection $\{(v,h):h\le h_1-c_2\|v-v_1\|\}\cap W_\la$ contains another cone $\{(v,h):h\le h_1-c_3\|v-v_1\|\}$ with $c_3$ depending only on $d$. Its $d\P$ measure is then bounded below by $c\exp(-e^{h_1/c})$.  The proof is concluded  as in the case $\lambda=\infty$.

We now prove \eqref{sbd}. 
We bound the probability of the events
  $$E_3:=\{\partial\Psi(\P^{(\la)})   \cap\{(v,h):||v||\le L,h >  t \}\ne \emptyset\}$$
   and $$E_4:=\{\partial\Psi(\P^{(\la)}) \cap\{(v,h): ||v|| \le L,h< - t \}\ne \emptyset\}.$$

When in $E_3$, there is a point $w_1:= (v_1,h_1)$ with $ h_1 \in [t, \infty), ||v_1|| \leq L,$ and such that $[\Pi^{\downarrow}(w_1)]^{(\la)}\cap \P^{(\la)}=\emptyset$.
As in the proof of  \eqref{fbd-1}, there is a subset of $C(\0,L)$ of volume one and on this subset
 the density of the $d\P^{(\la)}$  measure exceeds $c\exp(h_1/c)$.  Discretize  $\{(v,h): \ ||v|| \leq L, h \in [t, \infty)\}$
into unit volume sub-cubes and bound cross-sectional areas by $cL^{d-1}$ to obtain
 $$P[E_3] \leq cL^{d-1} \int_{c_1 t/3}^{\infty} \exp(dh_1) \exp( -c e^{h_1/c} ) dh_1
\leq cL^{d-1} \exp(-{ e^t\over c}).$$

On the event $E_4$, there exists a point $(v_1,h_1)$ with $||v_1||\le L$ and $h_1 \in (-\infty, - t]$ which is on the boundary of an up cone-like grain with apex in $\P^{(\la)}$. The apex of this up cone-like grain is
contained in the union of all down cone-like grains with apex on
$B_{d-1}(\0,L)\times \{h_1\}$. The $d\P^{(\la)}$  measure of this union is bounded by $cL^{d-1} \exp(h_1/c)$
(here we use that the union is a subset of the union of standard circular cones).
Consequently, the probability that the union contains
points from $\P^{(\la)}$ is less than $1-\exp(-cL^{d-1} e^{h_1/c}) \le cL^{d-1} \exp(h_1/c).$
It remains to discretize and integrate over
$h_1\in (-\infty, t)$.  This goes as follows.

Discretizing $C(\0,L) \times (-\infty,-t]$ into unit volume subcubes
and using  the previous bound,  we find that the probability  there exists
$(v_1,h_1) \in \R^{d-1} \times (-\infty, -t]$ on the boundary of an up cone-like grain
is bounded by
$$
 cL^{2(d-1)}\int_{-\infty}^{- t}e^{h_1/c} e^{dh_1}dh_1.
$$
 This establishes \eqref{sbd}.  Similar arguments apply to $\partial(\Phi(\P^{(\la)}))$.
 \qed

\vskip.5cm

For $(v_0, h_0)  \in {\rm{Ext}}(\P^{(\la)})$ and  $t \in \R$, we define
$${\mathcal U}^{(\la)}(v_0, h_0,t):=\bigcup_{w_1\in [\Pi^{\uparrow}((v_0,h_0))]^{(\la)}\cap (\R^{d-1}\times (-\infty, t])}[\Pi^{\downarrow}(w_1)]^{(\la)}.$$
  The score $\hx^{(\la)}((v_0,h_0), \P^{(\la)})$ depends only on the points of $\P^{(\la)}$ inside ${\mathcal U}^{(\la)}(v_0, h_0, H((v_0,h_0), \P^{(\la)}))$, as this set  contains all faces in $\Phi(\P^{(\la)})$ which contain $(v_0,h_0)$.
Put
\be \label{Rstab}
R:= R^{\hx^{(\la)}}[(v_0,h_0)] :=  \inf\{r > 0: \ \P^{(\la)} \cap  {\mathcal U}^{(\la)}(v_0, h_0, H((v_0,h_0), \P^{(\la)} )   ) \subset  C(v_0,r)\}. \ee
It follows from the definitions that
$$\hx^{(\la)}((v_0,h_0), \P^{(\la)}) = \hx^{(\la)}((v_0,h_0), \P^{(\la)}\cap C (v_0,r) ), \ r \in [R, \infty).
$$
In other words, as in Section 6 of \cite{CSY}, $R^{\hx^{(\la)}}[(v_0,h_0)]$  is a radius of spatial stabilization for  $\hx^{(\la)}$. 
The next lemma shows that $R$
is finite a.s. and in fact has exponentially decaying tails.
Given  $c_1$ as in \eqref{compar}, we put for all $h_0 \in \R$, \be \label{defho}
\tilde{h}_0 := (\frac{6}{c_1}\log d)\vee ((-\frac{6}{c_1}h_0){\bf{1}}(h_0 < 0)).\ee

\begin{lemm} \label{lem4.01}
There is a constant $c >0$ such that for all $\hx \in \hat{\Xi} $, $\la \in [1, \infty]$, $(v_0, h_0)  \in W_\la$, and all $t \in [\tilde{h}_0, \infty)$ we have
  \begin{equation}\label{LocalExpr1}
   P[R^{\hx^{(\la)}}[(v_0,h_0)] >  t] \leq c \exp(- \frac{t}{ c}).
\end{equation}
\end{lemm}
\noindent{\em Proof.}
We show \eqref{LocalExpr1} for  $v_0 = \0$, as the proof is analogous for arbitrary $v_0$.
Put $R:= R^{\hx^{(\la)}}[(v_0,h_0)]$ and write
$$
P[R > t] \leq P[ H((\0,h_0),\P^{(\la)}) \geq {c_1 t\over 6} ]  +
P[ H((\0,h_0),\P^{(\la)}) \in ( - \infty,  {c_1 t\over 6} ], R \geq t].
$$
Lemma  \ref{BLem1}(a) shows that the first term on the right-hand side is bounded by $c \exp(-e^{t/c}/ c)$.
Thus we only need to control the second term.

When $H((\0,h_0),\P^{(\la)}) \in (-\infty, c_1t/6]$, then
 $\hx^{(\la)}((\0,h_0))$ only depends on elements  of $\P^{(\la)}$
in
$${\mathcal U}:=  {\mathcal U}^{(\la)}(\0, h_0,c_1t/6).$$
Let $w=(v,h)\in {\mathcal U}$  and $w_1=(v_1,h_1)$, $h_1\le c_1 t/6,$  be such that $\partial [\Pi^{\downarrow}(w_1)]^{(\la)}$ contains both $(\0,h_0)$ and $w$.
We assert that if $h \in [-c_1t/6, c_1t/6]$, then $||v|| \leq t.$ To see this we first note that $||v - v_1|| \leq t/2$.  This follows
because $h = h_1 - G(v - v_1)$, which in view of \eqref{compar} yields $h \leq h_1 - c_1||v- v_1|| + \log d$, that is to say
\be \label{vv1} ||v - v_1|| \leq {1 \over c_1}[h_1 - h + \log d].\ee
Now  $||v - v_1|| \leq t/2$ since all three quantities $h_1, -h, $ and $\log d$ are bounded by $c_1t/6$.
Using $h_1 = h_0 + G(-v_1) \geq h_0 + c_1||v_1|| - \log d$ we get $||v_1|| \leq t/2$
and thus $||v|| \leq t$ by the triangle inequality.

Consequently, if $\P^{(\la)}\cap {\mathcal U}\cap (\R^{d-1}\times (-\infty, -c_1t/6])=\emptyset$, then
only elements of $\P^{(\la)}\cap {\mathcal U}\cap (\R^{d-1}\times (-c_1t/6, c_1t/6])$ contribute to the score
$\hx^{(\la)}((\0,h_0))$, showing that in this case
$R^{\hx^{(\la)}}[(\0,h_0)]  \in (0, t]$.  Therefore
$$
P[ H((\0,h_0),\P^{(\la)}) \in ( - \infty,  {c_1 t\over 6} ], R \geq t]  \leq
P[\P^{(\la)}\cap {\mathcal U}\cap (\R^{d-1}\times (-\infty,-c_1t/6])\ne\emptyset].
$$
Notice that if $v \in {\mathcal U}$ then $||v|| \leq ||v_1|| + ||v - v_1||  \leq \frac{t}{2} + \frac{t}{3} - \frac{h} {c_1}$, where we use
$||v_1|| \leq t/2$ and \eqref{vv1}.  Similar to the proof of Lemma \ref{BLem1},
discretization methods yield
\begin{align*}
& d\P^{(\la)}(\P^{(\la)}\cap {\mathcal U}\cap (\R^{d-1}\times (-\infty,  \frac{-c_1t} {6}])) \\& \le
c\int^{-c_1t/6}_{-\infty}e^{dh}( \frac{t}{2} + \frac{t}{3} - \frac{h} {c_1} )^{(d-1)}dh  \leq c \exp(- \frac{t}{ c}).
\end{align*}
It follows that
$$P[\P^{(\la)}\cap {\mathcal U}\cap (\R^{d-1}\times (-\infty, {-c_1t \over 6}])\ne\emptyset] \leq c \exp(- \frac{t}{ c}),$$
as desired.
\qed
\vskip.5cm


\begin{lemm}\label{L5.4}
For all $p \in [1, \infty)$ and $\hx\in \hat{\Xi}$, $\hx$ a $k$-face functional,
there is a constant $c>0$
such that for all $(v_0,h_0)  \in W_\la$, $\la \in [1, \infty]$, we have
\begin{equation}\label{LIMITBD2}
      \E [ \hx^{(\la)}( (v_0,h_0)),\P^{(\la)})^p] \leq c ( |h_0|^{c} + 1) \exp( - \frac{e^{(h_0 \vee 0)/c}}{c} ).
\end{equation}
For all $p \in [1, \infty)$ 
there is a constant $c>0$
such that for all $(v_0,h_0)  \in W_\la$, $\la \in [1, \infty]$, we have
\begin{equation}\label{LIMITBD2a}
      \E [ \hx^{(\la)}_V( (v_0,h_0)),\P^{(\la)})^p] \leq c ( |h_0|^{c} + 1) \exp(c(h_0 \vee 0) ) \exp( - \frac{e^{(h_0 \vee 0)/c}}{c} ).
\end{equation}
\end{lemm}

\noindent{\em Proof.}
We first prove \eqref{LIMITBD2} for the  $k$-face functional $\hx^{(\la)}:= \hx_k^{(\la)}$, $k\in \{0,1,...,d-1\}.$
 We start by showing for all $\la \in [1, \infty]$ and $h_0 \in \R$
\begin{equation}\label{LIMITBD1-new}
   \sup_{v_0 \in \R^{d-1}} \E [\hx^{(\la)}((v_0,h_0),\P^{(\la)} )^p] \leq c ( |h_0|^{c} + 1).
\end{equation}

Let $R:=R^{ \hx^{(\la)} }[(v_0,h_0)]$ be as at \eqref{Rstab} and
$N^{(\la)}((v_0,h_0))$   the cardinality of extreme points in $C(v,R)$ which share a common facet with $(v_0,h_0)$.
Clearly
$$
\hx^{(\la)}((v_0,h_0),\P^{(\la)})\leq \frac{1}{k+1}\binom{N^{(\la)}((v_0,h_0))}{k}.
$$
To show \eqref{LIMITBD1-new}, given $p \in [1, \infty)$, it suffices to show there is a constant $c:=c(p,k,d)$ such that for  $\la \in [1, \infty]$
\begin{equation}\label{momentsN}
\E N^{(\la)}((v_0,h_0))^{pk} \leq c ( |h_0|^{c} + 1).
\end{equation}

By \eqref{defP}, for all $r, \ell \in \R$  we have
$$d\P^{(\la)}(C(v_0,r)\cap (\R^{d-1} \times (-\infty,\ell))) \le cr^{d-1} e^{d \ell}.$$
Consequently, with $H:=H((v_0,h_0),\P^{(\la)})$ as defined before  Lemma \ref{BLem1} and with ${\rm{Po}}(\alpha)$ denoting a Poisson random variable
with mean $\alpha \in (0, \infty)$, we have for $\la \in [1, \infty]$
\begin{align*}
& \ \ \ \  \E N^{(\la)}((v_0,h_0))^{pk} \\
& \leq \E[ {\rm{card}}( \P^{(\la)} \cap [C(v,R)\cap ( \R^{d-1} \times (-\infty,H))])^{pk}] \\
& = \sum_{i= 0}^{\infty}\sum_{j=  \lfloor h_0 \rfloor }^{\infty} \E[ {\rm{Po}}(d\P^{(\la)}(C(v,R)\cap ( \R^{d-1} \times (-\infty,H)))^{pk} {\bf{1}}(i
\leq R < i + 1, j \leq H < j + 1) ] \\
& \leq \sum_{i= 0}^{\infty}\sum_{j=  \lfloor h_0 \rfloor }^{\infty}\E[ {\rm{Po}}(c(i + 1)^{d-1}e^{(j+1)d})^{pk} {\bf{1}}(R \geq  i, H \geq j) ].
\end{align*}

We shall repeatedly use the moment bound 
$\E[ {\rm{Po}}(\alpha)^r] \leq c(r) \alpha^r, r \in [1, \infty).$ Using H\"older's inequality, we get
\begin{equation*}
\E N^{(\la)}((v_0,h_0))^{pk} \le c\sum_{i= 0}^{\infty}\sum_{j=  \lfloor h_0 \rfloor }^{\infty}(i + 1)^{pk(d-1)/3}e^{(j+1)dpk/3}P[R\ge i]^{1/3}P[H\ge j]^{1/3}.
\end{equation*}
Splitting the sum on the $i$ indices into $i \in [0, \tilde{h}_0]$ and $i \in [\tilde{h}_0, \infty]$,  with $\tilde{h}_0$ defined at \eqref{defho},  and  splitting the sum
on the $j$ indices into $[\lfloor h_0 \rfloor \wedge 0 ,0]$ and $[0, \infty)$, we  get
$$ E N^{(\la)}((v_0,h_0))^{pk}  \leq  S_1 + S_2 + S_3 + S_4,$$
where
$$ S_1 :=  \sum_{i= 0}^{\tilde{h}_0}\sum_{j= \lfloor h_0 \rfloor \wedge 0 }^{0}  (i + 1)^{pk(d-1)/3}e^{(j+1)dpk/3}P[R\ge i]^{1/3}P[H\ge j]^{1/3} $$
$$S_2 :=  \sum_{i= \tilde{h}_0}^{\infty} \sum_{j= \lfloor h_0 \rfloor \wedge 0 }^{0}  (i + 1)^{pk(d-1)/3}e^{(j+1)dpk/3}P[R\ge i]^{1/3}P[H\ge j]^{1/3} $$
$$S_3:= \sum_{i= 0}^{\tilde{h}_0}\sum_{j= \lfloor h_0\rfloor\vee 0}^{\infty}  (i + 1)^{pk(d-1)/3}e^{(j+1)dpk/3}P[R\ge i]^{1/3}P[H\ge j]^{1/3}$$
$$
S_4:=   \sum_{i= \tilde{h}_0}^{\infty} \sum_{j= \lfloor h_0\rfloor\vee 0}^{\infty} (i + 1)^{pk(d-1)/3}e^{(j+1)dpk/3}P[R\ge i]^{1/3}P[H\ge j]^{1/3}.$$


Now we compute
$$
S_1 \leq c \sum_{i= 0}^{\tilde{h}_0}  (i + 1)^{pk(d-1)/3}  \sum_{j= \lfloor h_0 \rfloor \wedge 0 }^{0} \exp((j+ 1)dpk/3) \leq c ( |h_0|^{c} + 1),
$$
since  the second sum is bounded by a constant and where $c:=c(p,k,d)$.   Next,
$$
S_2 \leq  c \sum_{ i = \tilde{h}_0}^{\infty}  (i + 1)^{pk(d-1)/3} P[R\ge i]^{1/3}  \sum_{j= \lfloor h_0 \rfloor \wedge 0 }^{0} \exp((j+ 1)dpk/3)
\leq c$$
where the first sum converges by the exponentially decaying tail bound for $P[R\ge i]$.  Making use of the super exponentially decaying
tail bound for $P[H\ge j]$  we get  
$$S_3 \leq c\sum_{i= 0}^{\tilde{h}_0}  (i + 1)^{pk(d-1)/3}
 \sum_{j= \lfloor h_0\rfloor\vee 0}^{\infty}e^{(j+1)dpk/3}\exp( - e^{j/c}/3c)\leq c ( |h_0|^{c} + 1).
$$
Finally,
$$
S_4 \leq c \sum_{i= \tilde{h}_0}^{\infty}  (i + 1)^{pk(d-1)/3} P[R\ge i]^{1/3}  \sum_{j= \lfloor h_0\rfloor\vee 0 }^{\infty}e^{(j+1)dpk/3}\exp( - e^{j/c}/3c) \leq c,$$
since both sums are bounded by a constant.   Combining the bounds for $S_1, S_2, S_3$ and $S_4$ gives
the  required bound \eqref{momentsN}.

To deduce \eqref{LIMITBD2}, we argue as follows.
First consider the case $h_0 \in [0, \infty)$.
  By the Cauchy-Schwarz inequality and \eqref{LIMITBD1-new} we have
  \begin{align*}
  & \ \ \ \E  [\hx^{(\la)}((v_0,h_0),\P^{(\la)})^p] \\
  & \leq
  (\E\hx^{(\la)}((v_0,h_0),\P^{(\la)})^{2p})^{1/2} P[
  \hx^{(\la)}((v_0,h_0),\P^{(\la)}) > 0]^{1/2} \\
  & \leq c(|h_0| + 1)^c P[
  \hx^{(\la)}((v_0,h_0),\P^{(\la)}) \neq 0]^{1/2}.
  \end{align*}

  The event $\{\hx^{(\la)}((v_0,h_0), \P^{(\la)} ) \neq
  0\}$ coincides with the event that $(v_0,h_0)$ is extreme in $\P^{(\la)}$ and
  we may now apply \eqref{fbd-1} for $t = h_0$, which is possible since we have
  assumed $h_0$ is positive.  This gives \eqref{LIMITBD2} for $h_0 \in [0, \infty)$.
  When $h_0 \in (-\infty, 0)$ we  bound $P[
  \hx^{(\la)}((v_0,h_0),\P^{(\la)}) > 0]^{1/2}$ by $c \exp( - e^0/c)$, $c$ large,
  which shows \eqref{LIMITBD2} for $h_0 \in (-\infty,0)$.
This concludes the proof of \eqref{LIMITBD2} when $\hx$ is a $k$-face functional.

We now prove \eqref{LIMITBD2} for the volume functional $\hx_V$. We start by proving the analog of \eqref{LIMITBD1-new}.
Without loss of generality we put $(v_0,h_0) = (\0,h_0)$.  Recalling the definition of $H:= H((\0,h_0), \P^{(\la)})$ we have
$$\hx^{(\la)}_V((\0,h_0),\P^{(\la)}) \leq {1 \over d} \int_{v \in {\rm{Cyl}}^{(\la)}( (\0,h_0),\P^{(\la)})} dv \int_{-\infty}^H e^{dh} dh.
$$
Integrating, raising both sides to the $p$th power, taking expectations and applying the Cauchy-Schwarz inequality we get
$$
\E \hx^{(\la)}_V((\0,h_0),\P^{(\la)})^p \leq C ( \E( \Vol {\rm{Cyl}}^{(\la)}( (\0,h_0),\P^{(\la)})^{2p})^{1/2} ( \E e^{2pdH} )^{1/2}.
$$
Now
$$\E( \Vol {\rm{Cyl}}^{(\la)}( (\0,h_0),\P^{(\la)})^{2p})  \leq \E (  R^{\hx_V^{(\la)}}[(\0,h)] )^{2p(d-1)} \leq c (|h_0| + 1)^{c}
$$
by Lemma \ref{lem4.01}.   Lemma  \ref{BLem1}(a) and the formula $\E X = \int_0^{\infty} P[ X \geq t] dt$ imply that
\begin{align*}
\E e^{2pdH}  &  = \int_0^{\infty} P[ e^{2pdH} > t] dt  \\
& = \int_0^{\infty} P[ 2pdH > \log t ] dt \\
& = \int_{ t \leq \exp((h_0 \vee 0)2pd)  } P[ 2pdH > \log t ] dt + \int_{ t \geq \exp((h_0 \vee 0)2pd) } P[ 2pdH > \log t ]dt \\
& \leq \exp((h_0 \vee 0)2pd)  + \int_{t \geq \exp((h_0 \vee 0)2pd) } P[H \geq \log t^{1/2pd}] dt \\
& \leq c \exp((h_0 \vee 0)2pd).
\end{align*}
Thus
\begin{equation}\label{LIMITBD3-new}
\E \hx^{(\la)}_V((\0,h_0),\P^{(\la)} )^p \leq c (|h_0| + 1)^{c} \exp(c(h_0 \vee 0)).
\end{equation}
The bound \eqref{LIMITBD2a} for $\hx_V^{(\la)}$ follows from
\eqref{LIMITBD3-new} in the same way that \eqref{LIMITBD1-new}
implies \eqref{LIMITBD2} for $\hx_k^{(\la)}$.  This completes the proof of
Lemma \ref{L5.4}.
\qed

\vskip.5cm

{\bf \noindent{5.2. Two point correlation function for $\hx^{(\la)}$.}} For
all $h \in \R$, $(v_0,h_0), (v_1,h_1) \in W_\la,$ and $\hx\in\hat{\Xi} $
we extend the definition at \eqref{SO2}  by putting for
all $\la \in [1, \infty]$
\be \label{defcla}
c^{(\la)}((v_0,h_0), (v_1,h_1)) := c^{\hx^{(\la)}}((v_0,h_0), (v_1,h_1), \P^{(\la)} ) :=
\ee
$$
\E[ \hx^{(\la)}((v_0,h_0),\P^{(\la)} \cup (v_1,h_1)) \times
\hx^{(\la)}((v_1,h_1),\P^{(\la)} \cup (v_0, h_0))] - $$
$$\E \hx^{(\la)}((v_0,h_0),\P^{(\la)} ) \E
\hx^{(\la)}((v_1,h_1),\P^{(\la)}).
$$

The first part of the next lemma justifies the assertion that the functionals in $\hat{\Xi} ^{(\infty)}$ are scaling limits of their counterparts in $\hat{\Xi} ^{(\la)}$.

\begin{lemm}  \label{L2-new} (a) For all $(v_0,h_0) \in \R^{d-1} \times\R$ and $\hx\in \hat{\Xi} $
we have
$$\liml \E\hx^{(\la)}((v_0,h_0),\P^{(\la)}) =
 \E\hx^{(\infty)}((v_0,h_0),\P).$$
(b) For all $h_0 \in \R$, $(v_1,h_1) \in \R^{d-1} \times \R $ and $\hx \in \hat{\Xi} $  we have
$$\liml c^{\hx^{(\la)}}((\0,h_0), (v_1,h_1) ) =
c^{\hx^{(\infty)}}((\0,h_0), (v_1,h_1) ).$$
\end{lemm}

\noindent {\em Proof.} We first prove  (a).  Suppose $(v_0,h_0) \notin {\rm{Ext}} (\P).$  Then
 $(v_0,h_0) \notin {\rm{Ext}}( \P^{(\la)}),$
showing that both sides vanish.  Without loss of generality, let $(v_0,h_0) \in {\rm{Ext}}(\P).$
 Put $$B(v_0,h_0) :=C(v_0, R^{\hx^{( \infty)}}[(v_0,h_0), \P]  ) \cap
(\R^{d-1} \times ( -\infty, H((v_0,h_0), \P)])).$$
 We have
$$\hx^{(\infty)}((v_0,h_0),\P) = \hx^{(\infty)}((v_0,h_0),\P \cap B(v_0,h_0)).$$
For $\la$ large we have $B(v_0,h_0) \subset W_\la$.  For such $\la$ it follows that
$$\hx^{(\infty)}((v_0,h_0),\P) = \hx^{(\infty)}((v_0,h_0),\P \cap W_\la) = \hx^{(\la)}((v_0,h_0),\P \cap W_\la),$$
in other words
$$\hx^{(\la)}((v_0,h_0),\P \cap W_\la) \to \hx^{(\infty)}((v_0,h_0),\P) \ \ a.s.
$$
Convergence of expectations follows from the uniform integrability of $\hx^{(\la)}((v_0,h_0),\P \cap W_\la)$,
as shown in Lemma \ref{L5.4}.   This shows  part (a).  Part (b) follows from
identical methods,  since products of scores  $\hx^{(\la)}((\0,h_0),\P \cap W_\la)$  and
$\hx^{(\la)}((v_1,h_1),\P \cap W_\la)$ a.s. converge to their $\hx^{(\infty)}$ counterparts and
the products are uniformly integrable by Lemma \ref{L5.4}. \qed

\vskip.5cm

\begin{lemm} \label{L3} Let $c_1$ be as at \eqref{compar}.
Let $\hx \in \hat{\Xi}$ be a $k$-face functional.
There is a constant $c_3:=c_3(\hx,d) \in (0, \infty)$ such that for all $\la \in [1, \infty]$  and $(v_0, h_0), (v_1,h_1) \in W_{\la}$ satisfying
\be \label{assumv1vo} ||v_1 - v_0|| \geq 2 \max \left( \frac{6}{c_1}\log d, -\frac{6}{c_1}h_0 {\bf{1}}(h_0 < 0),
-\frac{6}{c_1}h_1 {\bf{1}}(h_1 < 0)\right)
\ee  we have 
\be \label{4.7a} |c^{\hx^{(\la)}}((v_0,h_0), (v_1,h_1))| \leq c_3 (|h_0| + 1)^{c_3}
(|h_1| + 1)^{c_3}
\exp\left( {-1 \over c_3} (||v_1 - v_0|| +e^{h_0 \vee 0}+e^{h_1 \vee 0}) \right).
\ee
When
$\hx$ is the volume functional $\hx_V$ we have
\begin{align}
\label{4.7aa} |c^{\hx^{(\la)}}((v_0,h_0), (v_1,h_1))| & \leq c_3 (|h_0| + 1)^{c_3}
(|h_1| + 1)^{c_3}  \exp( c_4((h_0 \vee 0) + (h_1 \vee0 )) )  \nonumber                  \\
& \times \exp\left( {-1 \over c_3} (||v_1 - v_0|| +e^{h_0 \vee 0}+e^{h_1 \vee 0}) \right).
\end{align}
\end{lemm}

\noindent{\em Proof.} We prove this assuming that $\hx$ is the $k$-face functional, as the
proof for the volume functional $\hx_V$ follows from identical methods.
Put
$$X_\la:= \hx^{(\la)}((v_0,h_0),\P^{(\la)} \cup (v_1,h_1)),$$
$$Y_\la:= \hx^{(\la)}((v_1,h_1),\P^{(\la)} \cup (v_0,h_0)),$$
$$\tX_\la:= \hx^{(\la)}((v_0,h_0),\P^{(\la)}) \ \ {\rm{and}} \  \tY_\la:= \hx^{(\la)}((v_1,h_1),\P^{(\la)}).$$
We have
\be
  \label{eq:decompc}
c^{\hx^{(\la)}}((v_0,h_0), (v_1,h_1)) = \E X_\la Y_\la - \E {\tX}_\la \E {\tY}_\la.
\ee
Put  $r := ||v_1 - v_0||/2$ and let $R^{\hx^{(\la)}}[(v_i,h_i)], i \in \{0, 1 \},$ be as at \eqref{Rstab}.
Now
\begin{align*} & |\E X_\la Y_\la - \E X_\la Y_\la {\bf{1}} ( R^{\hx^{(\la)}}[(v_0,h_0)]  \leq r,  R^{\hx^{(\la)}}[(v_1,h_1)] \leq r)| \nonumber\\
& \leq  \E X_\la Y_\la [{\bf{1}} ( R^{\hx^{(\la)}}[(v_0,h_0)]  \geq r) + {\bf{1}}(  R^{\hx^{(\la)}}[(v_1,h_1
)] \geq r)].
\end{align*}

Let $v_1$ and $v_0$ satisfy \eqref{assumv1vo}.
 H\"older's inequality and Lemma \ref{lem4.01} imply that the right hand side of the above is bounded by
\begin{align}
&||X_\la||_3 ||Y_\la||_3 [ P [R^{\hx^{(\la)}}[(v_0,h_0)]  \geq r]^{1/3} +
P [R^{\hx^{(\la)}}[(v_1,h_1)]  \geq r]^{1/3} ]\nonumber\\
& \le c (|h_0| + 1)^{c}
(|h_1| + 1)^{c}
\exp\left( {-1 \over c} (e^{h_0 \vee 0}+e^{h_1 \vee 0}) \right)  \nonumber\\
&
\times [ P [R^{\hx^{(\la)}}[(v_0,h_0)]  \geq r]^{1/3} +
P [R^{\hx^{(\la)}}[(v_1,h_1)]  \geq r]^{1/3} ] \nonumber\\
& \le c (|h_0| + 1)^{c}
(|h_1| + 1)^{c}
\exp\left( {-1 \over c} (||v_1- v_0|| +e^{h_0 \vee 0}+e^{h_1 \vee 0}) \right).  \label{eq:domin} \nonumber\\
\end{align}

Now
\begin{align*}
& \E X_\la Y_\la {\bf{1}} ( R^{\hx^{(\la)}}[(v_0,h_0)]  \leq r,  R^{\hx^{(\la)}}[(v_1,h_1)]  \leq r) \nonumber\\
& = \E \hx^{(\la)}((v_0,h_0),\P^{(\la)}\cap C (v_0,r) ) \hx^{(\la)}((v_1,h_1),\P^{(\la)}\cap C (v_1,r)  ) \nonumber\\
& \times {\bf{1}} ( R^{\hx^{(\la)}}[(v_0,h_0)]  \leq r,  R^{\hx^{(\la)}}[(v_1,h_1)]  \leq r).
\end{align*}
Following the above methods,  the difference of
$$\E X_\la Y_\la {\bf{1}} ( R^{\hx^{(\la)}}[(v_0,h_0)]  \leq r,  R^{\hx^{(\la)}}[(v_1,h_1)]  \leq r)$$
and $$\E \hx^{(\la)}((v_0,h_0),\P^{(\la)}\cap C (v_0,r) ) \hx^{(\la)}((v_1,h_1),\P^{(\la)}\cap C (v_1,r)  )$$
is also bounded by  \eqref{eq:domin}. By independence we have
$$\E \hx^{(\la)}((v_0,h_0),\P^{(\la)}\cap C (v_0,r) ) \hx^{(\la)}((v_1,h_1),\P^{(\la)}\cap C (v_1,r)  )$$
$$ = \E \hx^{(\la)}((v_0,h_0),\P^{(\la)}\cap C (v_0,r) )\E \hx^{(\la)}((v_1,h_1),\P^{(\la)}\cap C (v_1,r) ).$$
Thus we have  shown
\begin{align*}
& |\E X_\la Y_\la  - \E \hx^{(\la)}((v_0,h_0),\P^{(\la)} \cap  C(v_0,r) )\E \hx^{(\la)}((v_1,h_1),\P^{(\la)} \cap  C(v_1,r) )| \nonumber\\
& \le c (|h_0| + 1)^{c} (|h_1| + 1)^{c}
\exp\left( {-1 \over c} (||v_1- v_0|| +e^{h_0 \vee 0}+e^{h_1 \vee 0})\right).\end{align*}
Identical methods give
\begin{align*}
& |\E \tX_\la \E \tY_\la  - \E \hx^{(\la)}((v_0,h_0),\P^{(\la)} \cap  C(v_0,r) )\E \hx^{(\la)}((v_1,h_1),\P^{(\la)} \cap  C(v_1,r) )| \nonumber\\
& \le c (|h_0| + 1)^{c} (|h_1| + 1)^{c}
\exp\left( {-1 \over c} (||v_1- v_0|| +e^{h_0 \vee 0}+e^{h_1 \vee 0})\right).\end{align*}
Combining the last two displays with  \eqref{eq:decompc}, we get \eqref{4.7a}.
 \qed

 \vskip.3cm

 Our last lemma shows that $c^{\hx^{(\la)}}((\0,h_0), (v_1,h_1))e^{dh_0} e^{dh_1}$ is bounded by an integrable function, a fact
  used  in establishing variance asymptotics in the next section.

\begin{lemm} \label{L5} For all $\hx \in \hat{\Xi} $ there is an integrable $g: \  \R \times \R^{d-1} \times \R \to \R^+$  such that
for all $\la \in [1, \infty]$ we have
\be \label{5a} |c^{\hx^{(\la)}}((\0,h_0), (v_1,h_1))|  e^{dh_0} e^{dh_1} \leq g(h_0, v_1, h_1).
\ee
\end{lemm}

\noindent{\em Proof.}  With $c_1$ as at  \eqref{compar}, define  $F: \  \R \times \R^{d-1} \times \R \to \R^+ $
by
\begin{align*}
& F(h_0, v_1, h_1) :=c(|h_0| + 1)^{c} (|h_1| + 1)^{c}  \exp( c_4((h_0 \vee 0) + (h_1 \vee0 )) )  \\
& \times  \left( \exp\left( {-1 \over c} (||v_1|| +e^{h_0 \vee 0}+e^{h_1 \vee 0}) \right)\right.\\&\hspace*{1cm}\left. +
{\bf 1}( ||v_1||
\leq  2 \max \left( \frac{6}{c_1}\log d, -\frac{6}{c_1}h_0 {\bf{1}}(h_0 < 0),
-\frac{6}{c_1}h_1 {\bf{1}}(h_1 < 0)\right)  \right),
\end{align*}
where $c$ is a constant. If $c$ is large enough, then Lemma \ref{L3} gives
$$|c^{\hx^{(\la)}}((\0,h_0), (v_1,h_1))|  e^{dh_0} e^{dh_1} \leq F(h_0, v_1, h_1)  e^{dh_0} e^{dh_1}.$$
Put $g(h_0, v_1, h_1):= F(h_0, v_1, h_1)  e^{dh_0} e^{dh_1}$ and note that $g$ is integrable as claimed.  \qed

\section{Proof of main results} \label{proofs}
\allco

{\bf 6.1. Proof of Theorems \ref{Th1} and \ref{Th2}}.
The next proposition immediately yields  Theorem \ref{Th2}.  It also yields Theorem \ref{Th1},
since it implies that the extreme points of $K_\la \cap Q_0$ converge in law to  ${\rm{Ext}}(\P)$ as $\la \to \infty$.
Recall that $T^{(\la)}(\P_\la \cap Q_0): = \P^{(\la)}$ as at \eqref{TPQ} and
$T^{(\la)}( (\partial K_\la) \cap Q_0): = \partial( \Phi( \P^{(\la)}))$ on the event $A_\la$.


\begin{prop}\label{BLem2} Fix $L \in (0, \infty).$  We have that $\partial \Psi(\P^{(\la)})$
converges in probability as $\la \to \infty$ to  $\partial \Psi(\P)$
in the space $\C(B_{d-1}(\0,L)))$. Likewise,  $\partial \Phi(\P^{(\la)})$ converges in probability  as
$\la \to \infty$ to $\partial(\Phi(\P))$.
\end{prop}

\noindent{\em Proof.}  We prove the first convergence statement as follows.
With $L$
fixed, for all $l \in [0, \infty)$ and $\la \in [1, \infty)$, let
$E(L, l, \la)$ be the event that the heights of
$\partial(\Psi(\P^{(\la)}))$ and $\partial( \Psi(\P))$
belong to $[-l, l]$ over the spatial region $B_{d-1}(\0,L)$.
Lemma \ref{BLem1}(b) shows that $P[E(L, l, \la)^c]$ decays
exponentially fast in $l$, uniformly in $\la$.  It is enough
to show, conditional on $E(L, l, \la)$, that
$\partial(\Psi(\P^{(\la)}))$ and
$\partial(\Psi(\P))$ coincide with high probability in the space $\C(B_{d-1}(\0,L))$, $\la$ large. Indeed, conditional on $E(L,l,\la)$, $\partial \Psi(\P)\cap (B_{d-1}(\0,L)\times [-l,l])$ depends only on points in
$${\mathcal D}:=\bigcup_{w\in B_{d-1}(\0,L)\times [-l,l]}\Pi^{\downarrow}(w).$$
Thus whenever we have equality of $\P \cap {\mathcal D} \cap W_\la$ and $\P \cap {\mathcal D}$, it follows that
$\partial \Psi(\P^{(\la)})$ and $\partial\Psi(\P)$ coincide in $C_{d-1}[\0,L]$.
Since ${\mathcal D} \setminus W_\la$  decreases to $\emptyset$, we have as $\la$ goes to infinity
 $$P[\P\cap {\mathcal D}\cap W_\la\ne \P\cap {\mathcal D}]=P[\P\cap ({\mathcal D} \setminus W_\la)\ne \emptyset]\le d\P [{\mathcal D} \setminus W_\la]\to 0.$$
This completes the proof of the first convergence statement.  The proof of the second convergence statement is nearly identical
and we leave the details to the reader. \qed

\vskip.3cm

\noindent{\bf 6.2. Proof of Theorem \ref{Th5}}.
When $g\equiv 1$,  the decomposition \eqref{decomp-e} shows that  it is enough to find expectation asymptotics for
 $\E [\langle  {\bf 1}(Q_0), \mu_\la^\xi \rangle  {\bf{1}}(A_\la)]$  and multiply the result by $f_0(K)$.
For  arbitrary $g \in {\cal C} (K)$, an identical decomposition holds and so
to show \eqref{main1}, it suffices to find $\liml \E [\langle g{\bf 1}(Q_0), \mu_\la^\xi \rangle {\bf{1}}(A_\la) ]$.
We have
\begin{align*}
& \ \ \ \E [\langle g {\bf 1}(Q_0), \mu_\la^\xi  \rangle {\bf{1}}(A_\la) ]\\
& = \int_{Q_0} g(x) \E [\xi(x, \P_\la) {\bf{1}}(A_\la) ] \la dx \\
& = \sqrt{d}\int_{(v,h) \in W_\la } g( [T^{(\la)}]^{-1} (v,h)) \E [\hx^{(\la)} ((v,h), \P \cap W_\la)  {\bf{1}}(A_\la)] e^{dh} dh dv \\
& = \sqrt{d}\int_{(v,h) \in W_\la } g( [T^{(\la)}]^{-1} (v,h)) \E [\hx^{(\la)} ((\0,h), \P \cap (W_\la - v) )  {\bf{1}}(A_\la)]e^{dh} dh dv
\end{align*}
where the second equality uses \eqref{Key1} and \eqref{Key2}, whereas the last equality uses translation invariance of $\hx^{(\la)}$. Scaling by $(\log \la)^{d-1}$ and making the change of variable
$u = ( {1 \over d} \log \la)^{-1} v$, $dv = d^{-(d-1)} (\log \la)^{d-1} du$,
 we obtain
\be \label{integral}
(\log \la)^{-(d-1)} \E [\langle g {\bf 1}(Q_0), \mu_\la^\xi  \rangle   {\bf{1}}(A_\la)]
= $$ $$ d^{-d+3/2} \int_{(u,h) \in {W_\la}' } g( [T^{(\la)}]^{-1} ( ({1 \over d} \log \la) u  ,h))
\E [\hx^{(\la)} ((\0,h), \P \cap ({W_\la}' - u) \log \la^{1/d} ) {\bf{1}}(A_\la)] e^{dh} dh du,
\ee
where ${W_\la}':= \{ ( {1 \over d} \log \la)^{-1} v, h); \ (v,h) \in W_\la \}$.  Here, for
$B \subset \R^{d-1} \times \R$ and $s \in \R$, we write $sB:= \{(sv,h): \ (v,h) \in B \}.$
We now prove \eqref{main1} via the following three steps.

(i) We first show the almost everywhere convergence
\begin{equation}
  \label{eq:convg}
\liml {\bf 1}( (u,h) \in {W_\la}' )g\left( [T^{(\la)}]^{-1} ( ({1 \over d} \log \la) u  ,h)\right)={\bf 1}( u \in S(d))g(\0),\end{equation}
where $S(d)$ is defined at \eqref{defSd}.  
Indeed, because of \eqref{eq:Wlambda}, the equation of ${W_\la}'$ is  $\ell_i(u) \leq (1 + \frac{\log \delta_0-h}{\log \la^{1/d}})$, $1\le i\le d$. Consequently, in the limit as $\la \to \infty$, we have $\ell_i(v) \leq 1$ for $1\le i\le d$. In other words, the limit of ${W_\la}'$ is a cylinder whose base is  the intersection of $V$ and the pyramid $\{(x_1,\cdots,x_d)\in\R^d: x_i\le 1, 1\le i\le d\}$.  This base is precisely $S(d)$.

Moreover in view of \eqref{eq:invscaltrans}, we have
$$[T^{(\la)}]^{-1} ( ({1 \over d} \log \la) u  ,h)=\la^{\frac{1}{d}(\ell(u)-1)}e^h.$$
If $u\not\in S(d)$, then for $\la$ large enough, the indicator function is equal to $0$. If $u\in \mbox{int}S(d)$ (where $\mbox{int}$ denotes the interior), then  $l_i(u)<1$ for $1\le i\le d$ and thus $\liml [T^{(\la)}]^{-1} ( ({1 \over d} \log \la) u  ,h)=\0$.
 By continuity of $g$ we have $$\liml g([T^{(\la)}]^{-1} ( ({1 \over d} \log \la) u  ,h))=g(\0).$$ This shows  \eqref{eq:convg}.

(ii) We  remove the indicator on the right hand side of \eqref{integral} with small error:
\begin{align}\label{eq:expsimplified}
& (\log \la)^{-(d-1)} \E [\langle g {\bf 1}(Q_0), \mu_\la^\xi \rangle  {\bf{1}}(A_\la)]= d^{-d+3/2} \int g( [T^{(\la)}]^{-1} ( ({1 \over d} \log \la) u  ,h))\nonumber
\\&\hspace*{1cm}
 \E [\hx^{(\la)} ((\0,h), \P \cap ({W_\la}' - u) \log \la^{1/d} )]  {\bf 1} ({(u,h) \in {W_\la}' })
 e^{dh} dh du  + o(1).
\end{align}
Indeed, by the Cauchy-Schwarz inequality and by moment bounds similar to those from Lemma \ref{L5.4},  we have uniformly in  $u$ that
\begin{align*}
& \int \E [\hx^{(\la)} ((\0,h), \P \cap ({W_\la}' - u) \log \la^{1/d})  {\bf{1}}(A_\la^c)] e^{dh} dh \\
&  \leq \int (\E [\hx^{(\la)} ((\0,h), \P \cap ({W_\la}' - u) \log \la^{1/d})]^2)^{1/2} P[A_\la^c]^{1/2}  e^{dh} dh \\
& \leq c (\log \la)^{-2d^2},\end{align*}
where $c$ is a constant not depending on $u$. Equality \eqref{eq:expsimplified} follows  from the estimate above, \eqref{eq:convg}, and the dominated convergence theorem.

(iii) Since $({W_\la}' - u) \log \la^{1/d} \uparrow \R^d$ as $\la \to \infty$, an easy modification of the proof of Lemma \ref{L2-new}
gives
\begin{equation}
  \label{eq:convmeanscore}
 \liml \E \hx^{(\la)} ((\0,h), \P \cap ({W_\la}' - u) \log \la^{1/d} )
= \E \hx^{(\infty)} ((\0,h), \P ).
\end{equation}

Lemma \ref{L5.4} shows that $\E \hx^{(\la)} ((\0,h), \P \cap (W_1 - {1 \over d} u) \log \la ) e^{dh} $ is
dominated by an integrable function on $\R^{d-1} \times \R$.  Combining  \eqref{eq:convg}- \eqref{eq:convmeanscore}
yields \eqref{main1} as desired.


\vskip.5cm
Next we show variance asymptotics \eqref{main2}. By an easy extension of the decomposition \eqref{decomp-a} and Lemma \ref{Lem9}, it suffices
to find  
$$\liml \Var \sum_{x \in \P_\la \cap Q_0 } \xi(x, \P_\la) g(x){\bf{1}}(A_\la).$$

 For $g \in {\cal C}(K),$  the Mecke-Slivnyak formula (Corollary 3.2.3 in \cite{SW}) gives  \be
 \Var  [\langle g{\bf 1}(Q_0), \mu^{\xi}_{\la}  \rangle {\bf{1}}(A_\la)] :=  I_1(\la) + I_2(\la), \label{dis3}
 \ee
 where
 $$
 I_1(\la) :=  \int_{Q_0} g(x)^2 {\Bbb E}\left[\xi(x,\P_{\la})^2  {\bf{1}}(A_\la) \right] \la dx $$
  and
$$
I_2(\la):=
  \int_{Q_0} \int_{Q_0} g(x) g(y)
[ \E \xi(x, \P_\la \cup {y} )\xi(y, \P_\la \cup {x} ){\bf{1}}(A_\la) - \E \xi(x,
\P_\la ) {\bf{1}}(A_\la)\E \xi(y, \P_\la ) {\bf{1}}(A_\la)] \la^2 dy dx.
$$
Replacing  $g$ by $g^2$  in the proof of expectation asymptotics, we obtain
\be \label{VL1}
\liml (\log \la)^{-(d-1)} I_1(\la) = d^{-d+3/2}   \Vol_d(S(d)) \int_{-\infty}^{\infty} \E [\xi^{(\infty)}((\0,h_0), \P)^2] e^{dh_0} dh_0 g^2(\0).
\ee
We next consider $\liml (\log \la)^{-(d-1)} I_2(\la)$. Recalling \eqref{defcla} we have
\begin{align*}
& \ \ \ c^{\hx^{(\la)}{\bf 1}(A_\la) }  ((v_0,h_0), (v_1,h_1))  \\
& =
\E[ \hx^{(\la)}((v_0,h_0),\P^{(\la)} \cup (v_1,h_1)) \times
\hx^{(\la)}((v_1,h_1),\P^{(\la)} \cup (v_0, h_0)) {\bf{1}} (A_\la)] \\
& \ \ - \E [\hx^{(\la)}((v_0,h_0),\P^{(\la)}){\bf{1}} (A_\la)] \E
[\hx^{(\la)}((v_1,h_1),\P^{(\la)}){\bf{1}} (A_\la)].
\end{align*}
By \eqref{Key1}  and \eqref{Key2}, we may rewrite $I_2(\la)$ in terms of $\hx^{(\la)}$:
$$
I_2(\la)  = d \int_{(v_0,h_0) \in W_\la} \int_{(v_1,h_1) \in W_\la}  g([T^{(\la)}]^{-1}(v_0,h_0)) g([T^{(\la)}]^{-1}(v_1,h_1))
$$
$$\cdot
c^{\hx^{(\la)}{\bf 1}(A_\la) }((v_0,h_0), (v_1,h_1))  e^{dh_0} e^{dh_1} dh_0 dh_1 dv_0 dv_1.$$
Translation invariance of $\hx^{(\la)}$ yields
$$
I_2(\la) = d\int_{(v_0,h_0) \in W_\la} \int_{(v_1,h_1) \in W_\la} g( [T^{(\la)}]^{-1}(v_0,h_0) ) g( [T^{(\la)}]^{-1}(v_1,h_1) )
$$$$ \cdot c^{\hx^{(\la)}{\bf 1}(A_\la) } ((\0,h_0), (v_1-v_0,h_1), \P \cap (W_\la - v_0))
e^{dh_0} e^{dh_1} dh_0 dh_1 dv_1 dv_0.
$$
Again, we make the change of variable $u = ( {1 \over d} \log \la)^{-1} v_0$.
  This gives
$$
(\log \la)^{-(d-1)} I_2(\la) = d^{-d+2} \int_{(u,h_0) \in {W_\la}' } \int_{(v_1,h_1) \in W_\la} g( [T^{(\la)}]^{-1}({1 \over d} \log \la \cdot u,h_0) ) g( [T^{(\la)}]^{-1}(v_1,h_1) )
$$$$ \cdot c^{\hx^{(\la)}{\bf 1}(A_\la) } \left((\0,h_0), (v_1-{1 \over d} \log \la \cdot u, h_1), \P \cap (W_\la - {1 \over d} \log \la \cdot u)\right)
e^{dh_0} e^{dh_1} dh_0 dh_1 dv_1 du.
$$
To conclude the proof of  \eqref{main2}, it remains only to compute $\liml (\log \la)^{-(d-1)} I_2(\la)$. We proceed in four steps.
\vskip.3cm
\noindent(i) Similarly to \eqref{eq:convg}, we have
\begin{align}
& \liml {\bf 1}( (u,h_0) \in {W_\la}' ){\bf 1}( (v_1,h_1) \in {W_\la} )g( [T^{(\la)}]^{-1}({1 \over d} \log \la \cdot u,h_0) ) g( [T^{(\la)}]^{-1}(v_1,h_1) )\nonumber\\& = {\bf 1}(u\in S(d))g^2(\0).\label{eq:convvar}
 \end{align}

\noindent (ii) We remove the indicator from $c^{\hx^{(\la)}{\bf 1}(A_\la) }$ at the cost of a small additive error.
Recall the definition  of $c^{\hx^{(\la)}}((v_0,h_0), (v_1,h_1), \P^{(\la)} )$ at \eqref{defcla}. As in the proof of Lemmas \ref{L3} and \ref{L5}
 we show that the difference
\begin{align*}
&  |c^{\hx^{(\la)}{\bf 1}(A_\la) } \left((\0,h_0), (v_1-{1 \over d} \log \la \cdot u, h_1), \P \cap (W_\la - {1 \over d} \log \la \cdot u)\right)\\
&-c^{\hx^{(\la)}} \left((\0,h_0), (v_1-{1 \over d} \log \la \cdot u, h_1), \P \cap (W_\la - {1 \over d} \log \la \cdot u)\right)|e^{dh_0}e^{dh_1}
\end{align*}
is bounded above by $P[A_{\la}^c]^{1/4}G(h_0,v_1,h_1)$ where $G$ is a function which is integrable with respect to $(h_0,v_1,h_1)\in \R\times\R^{d-1}\times\R$
and which does not depend on $u$.
Consequently, the integrated error is bounded,  uniformly in $u$:
\begin{align}\label{eq:interror}
\int_{h_0\in\R} \int_{(v_1,h_1) \in W_\la} P[A_{\la}^{c}]G(h_0,v_1,h_1)dh_0 dh_1 dv_1
& \leq C ( \log \la)^{-d^2} \nonumber\\
& = o((\log \la)^{(d-1)}).
\end{align}
Combining  \eqref{eq:convvar} and \eqref{eq:interror} with the  dominated convergence theorem, we find  that the removal of the indicator does not modify the asymptotics of the variance.
\vskip.3cm
\noindent(iii) Given a fixed $u$, for all $(v_1,h_1) \in W_\la$, we make the change of variable $v' = v_1-{1 \over d} \log \la \cdot u$, $dv' = dv_1$. This
transforms $$c^{\hx^{(\la)}} \left((\0,h_0), (v_1-{1 \over d} \log \la \cdot u, h_1), \P \cap (W_\la - {1 \over d} \log \la \cdot u)\right) e^{dh_0} e^{dh_1}$$ into
$$c^{\hx^{(\la)}} \left((\0,h_0), (v', h_1), \P \cap ({W_\la}' - u) \log \la^{1/d} \right) e^{dh_0} e^{dh_1}.$$ By
Lemma \ref{L5} the last expression is bounded by an integrable function of $h_0, v'$ and $h_1$, uniformly in $\la$,
and by Lemma \ref{L2-new} it converges to $c^{\hx^{(\infty)}}((\0,h_0), (w, h_1), \P)$ as $\la \to \infty$.
\vskip.3cm
\noindent(iv) We make the change of variable $v' = v_1-{1 \over d} \log \la \cdot u$.  The integration domain $W_\la$ transforms to
 $\{(v',h_0) \in  ({W_\la}' - u) \log \la^{1/d}\}$, which increases
up to $\R^{d-1} \times \R$.
\vskip.3cm
Combining  observations (i)-(iv) with  the dominated convergence theorem yields
\begin{align} \label{VL2}
&\liml (\log \la)^{-(d-1)} I_2(\la) =\nonumber\\& \hspace{.3cm}d^{-d+2}   \Vol_d(S(d)) g^2(\0) \int_{-\infty}^{\infty}  \int_{\R^{d-1}} \int_{-\infty}^{\infty}
   c^{\xi^{(\infty)}}((\0,h_0),(v',h_1)) e^{d(h_0 + h_1)} dh_0 dv' dh_1.
\end{align}
Combining  \eqref{VL1} and  \eqref{VL2} and recalling the definition of $\sigma^2(\xi^{(\infty)})$ at
\eqref{S03} gives
$$
\liml (\log \la)^{-(d-1)}  \Var  [\langle g, \mu^{\xi}_{\la} {\bf 1}(Q_0) \rangle]= d^{-d+1}   \Vol_d(S(d)) g^2(\0) \sigma^2(\xi^{(\infty)}).
$$
We repeat this computation for each vertex of $K$.  Proposition \ref{Prop2} yields
\eqref{main2}, as desired.  \qed

\section{Appendix }\label{sec:app}

\allco

 We establish the unproved assertions of Section 3. Our first lemma shows that, near the origin, the boundary of the floating body for $K$ is a pseudo-hyperboloid.
\begin{lemm} \label{hyper}
There exists $\Delta_d\in [1,\infty)$ depending only on $d$ such that when $K$ contains $[\0,\Delta_d]^d$ and is contained in some multiple of that cube, then 
\begin{equation}
  \label{eq:floatinghyperb}
K(v=t)\cap  \left[0, \frac{1}{2}\right]^d   =\left\{(z_1,\cdots,z_d)\in \left[0,\frac{1}{2}\right]^d : \ \prod_{i=1}^d z_i=  \frac{d!}{d^d}t \right\},  \ \ t \in (0, \infty).
\end{equation}
\end{lemm}
\noindent{\em Proof.} Put $\tilde{t}:=d! t/d^d$. Recall the definition of the surface ${\mathcal H}_{\tilde{t}}$  at \eqref{eq:hyperboloid}.
We start by proving that for every $z^{(0)}=(z_1^{(0)},\cdots,z_d^{(0)})\in {\mathcal H}_{\tilde{t}}$, we have
\begin{equation}
  \label{eq:volcaphyp}
\Vol([0,\infty)^d\cap H^+(z^{(0)}))=t
\end{equation}
 where we recall that $H^+(z^{(0)})$ is the half-space containing the origin and bounded by the hyperplane tangent to ${\mathcal H}_{\tilde{t}}$ at $z^{(0)}$.
By  \eqref{eq:equationhyperplantangent},  we have
 \begin{equation}
   \label{eq:eqhypertangent}
 {H}^+(z^{(0)}):= \{ (z_1,...,z_d): \  \sum_{i=1}^d\frac{z_i}{z_i^{(0)}}\le d\}.
 \end{equation}
Then
  \begin{align*}
\Vol([0,\infty)^d\cap H^+(z^{(0)}))&=\int_{[0,\infty)^d}{\bf 1}(\sum_{i=1}^d\frac{z_i}{z_i^{(0)}}\le d)dz_1\cdots dz_d\\&=\tilde{t}\int_{[0,\infty)^d}{\bf 1}(\sum_{i=1}^dy_i\le d)dy_1\cdots dy_d=t
  \end{align*}
 where we use the change of variable $y_i=\frac{z_i}{z_i^{(0)}}$, $1\le i\le d,$ and the identity $\prod_{i=1}^dz_i^{(0)}=\tilde{t}$. The proof of \eqref{eq:volcaphyp} is complete.

 We now prove that the boundary of the floating body for $[0,\infty)^d$ at level $t$ satisfies
 \begin{equation}
   \label{eq:floatingquadrant}
 [0,\infty)^d(v=t)={\mathcal H}_{\tilde{t}}.
 \end{equation}
 Indeed, let $H^+$ be a half-space with a boundary denoted by $H$ such that $H$ contains $z^{(0)}$.
We may assume that $H$ has a normal vector with all strictly positive coordinates (otherwise, we would have $\Vol([0,\infty)^d\cap H^+)=\infty$). Then $H$ is tangent to exactly one ${\mathcal H}_{\tilde{s}}$ where $s=\Vol([0,\infty)^d\cap H^+)$. Since $z^{(0)}\in {\mathcal H}_{\tilde{t}}\cap H$,
then necessarily $t\in (0, s)$. Consequently, we  obtain
 $z^{(0)}\in [0,\infty)^d(v=t)$. In other words, ${\mathcal H}_{\tilde{t}}\subset [0,\infty)^d(v=t)$. That this holds for any $t \in (0,s]$,
combined with the fact that all ${\mathcal H}_{\tilde{t}}$, $t>0$, and all $[0,\infty)^d(v=t)$, $t>0,$ form a partition  $(0,\infty)^d$,
is enough to yield \eqref{eq:floatingquadrant}.

It remains to show the validity of  \eqref{eq:floatingquadrant}  when replacing the orthant $[0,\infty)^d$ with $K$ whenever $K$ contains a
 large enough cube $[0, \Delta_d]^d$.  We remark that as soon as $\Delta_d \ge d/2$, we have for every $t>0$ and every $z^{(0)}\in  [0,\frac{1}{2}]^d \cap {\mathcal H}_{\tilde{t}}$, the equality $$\Vol([0,\infty)^d\cap H^+(z^{(0)}))=\Vol(K\cap H^+(z^{(0)}))=t.$$
Indeed, we observe that $[0,\infty)^d\cap H^+(z^{(0)})\subset [0,\frac{d}{2}]^d\subset K$ because the equation \eqref{eq:eqhypertangent} defining $H^+(z^{(0)})$ implies that every $z\in H^+(z^{(0)})$ satisfies $z_i\le z_i^{(0)}d\le d/2$ for $1\le i \le d$.

We now fix again $z^{(0)}\in  [0,\frac{1}{2}]^d  \cap {\mathcal H}_{\tilde{t}}$ and consider a half-space $H^+$ with boundary $H$ such that $H$ contains $z^{(0)}$. We have to show that
\begin{equation}
  \label{eq:captropgrand}
\Vol(K\cap H^+)\ge t.
\end{equation}
When $H^+\cap (0,\infty)^d$ is a subset of $[0, \Delta_d]^d$ then $K\cap H^+\supset [0, \Delta_d]^d \cap H^+$
and \eqref{eq:captropgrand} follows.  If $H^+\cap (0,\infty)^d$ is not a subset of $[0, \Delta_d]^d$
then there is at least one point $z=(z_1,\cdots,z_d)$ from $H$ with a coordinate greater than $\Delta_d$, say $z_1$.
In particular, $H^+\cap [0, \Delta_d]^d$ contains a simplex which is the convex hull of $(\{z_1^{(0)}\}\times \prod_{i=2}^d[0,z_i^{(0)}])\cup \{z\}$. Consequently, $\Vol(K\cap H^+)$ is bounded from below by $c \Delta_d \prod_{i=2}^dz_i^{(0)}$ where $c$ is a multiplicative constant depending only on $d$. When $\Delta_d \ge {d^d}/{2c d!}$, we have
$$c\Delta_d \prod_{i=2}^dz_i^{(0)}\ge \frac{d^d}{c d!}z_1\prod_{i=2}^dz_i^{(0)}=t,$$
which completes the proof of Lemma \ref{hyper}. Notice that these arguments show that  we can take $\Delta_2=1$.  \qed

\vskip.3cm

\noindent{\em Proof of Proposition \ref{Prop1}.} Recall from Section 3.5 that  $\delta_1:= r(\la, d) \delta_0$,  where
$ r(\la, d) \in [1, 3^{1/d} )$ is chosen so that $\log_3(T/\delta_1^d)\in \Z$.
We show the slightly stronger result that
$C_F(K_\la \cap [0, \delta_1]^d) \subset C_{\0}(K)$ holds on $A_\la$.
 Assume there is a normal $u \in C_F(K_\la \cap [0, \delta_1]^d)$ with $u \notin  C_{\0}(K)$.
Thus there is some  $j \in \{1,...,d-1 \}$ such that the first $j$ coordinates
 of $u$ are positive and the last
$(d- j)$ coordinates are negative.

Let  $H_u$ be the support hyperplane containing $F$ and let $z \in H_u \cap [0, \delta_1]^d \cap {\cal A}(s, T^*,K).$ The existence of $z$ is guaranteed on the event $A_\la$.
The definition of $u$ shows that $\langle z' - z, u \rangle \leq 0$ must hold for
all $z' \in K_\la$.  However this is not case  
and we assert there is a $z':= (z'_1,...,z'_d) \in \P_\la$ such that
\be \label{comps} z'_1 - z_1 > 0,...,z'_j - z_j > 0;  \ \ z'_{j + 1} - z_{j + 1} < 0,...,z'_d - z_d < 0.\ee
We  prove \eqref{comps} as follows.

Since $z \in [0, \delta_1]^d$, the point $z$ must belong to an $M$-region
$\Pi_{i = 1}^d [3^{k_i}\delta_1/2, 3^{k_i+1}\delta_1/2]$.
 We shall show there exist integers $k'_1,...,k'_d$ with
$k'_i \leq \log_3(\delta_1^{-1}), 1 \leq i \leq d$,  such that
\be \label{claim1}
k'_i > k_i \ {\rm{for}}  \ 1 \leq i \leq j \ {\rm{whereas}}  \ \ k'_i <k_i \ {\rm{for}} \ j +1 \leq i \leq d,\ee
and $\sum_{i = 1}^d k'_i = \log_3(T/\delta_1^d).$  The last equality implies that the $M$-region
 $$M:= \prod_{i = 1}^d [  \frac{ 3^{k'_i}\delta_1} {2}, \frac{3^{k'_i+1}\delta_1} {2}]$$
is an element of ${\cal M}_K(\0,\delta_1)$.  On $A_\la$ we
 know that $M$ contains at least one element of $\P_\la \cap {\cal A}(s, T^*,K) \cap [0, 1/2]^d$, say $z'$, with coordinates
$z'_i, 1 \leq i \leq d,$ satisfying \eqref{comps}.  This shows $\langle z' - z, u \rangle \geq 0$ as desired.   To show \eqref{claim1} we consider two cases.  Note that
$\sum_{i = 1}^d k_i \in [\log_3 s/\delta_1^d,  \log_3 T^*/\delta_1^d].$

\vskip.3cm
\noindent Case (i).  $\sum_{i = 1}^d k_i \in [\log_3 s/\delta_1^d,  \log_3 T/\delta_1^d].$   We choose $k'_i< k_i$ for $j + 1 \leq i \leq d$ such that $$j\le \sum_{i=j+1}^d(k_i-k'_i) \le d.$$  Then for
$1 \leq i \leq j$, we choose $k'_i \in [k_i + 1, \log_3 \delta_1^{-1} ]$ such that
$$
\sum_{i = 1}^j (k_i' - k_i) = \left(\log_3 T/\delta_1^d-\sum_{i = 1}^d k_i\right)+  \sum_{i=j+1}^d(k_i-k'_i).$$
Such $k'_i$ exist since $\log_3 T/\delta_1^d- \sum_{i = 1}^d k_i$  is bounded by the maximum allowable value
in the range of $k'_i$, that is to say it is bounded by $\log_3 T/\delta_1^d- \sum_{i = 1}^d k_i = o( \log_3 \delta_1^{-1}) $.
Thus \eqref{claim1} holds in this situation.

\vskip.3cm
\noindent Case (ii).  $\sum_{i = 1}^d k_i \in [\log_3 T/\delta_1^d,  \log_3 T^*/\delta_1^d].$  We choose $k'_i>k_i$
 for $1 \leq i \leq j$ such that $$(d-j) \le \sum_{i=1}^j(k'_i-k_i)\le d.$$ Then  for
$j + 1 \leq i \leq d$, we choose $k'_i \leq k_i - 1$ such that $$\sum_{i = j + 1}^d (k_i - k'_i) =
\left(\sum_{i = 1}^d k_i - \log_3 T/\delta_1^d\right)+ \sum_{i=1}^j(k'_i-k_i).$$
 This shows \eqref{claim1}, completing the proof of Proposition \ref{Prop1}.
 \qed

\vskip.5cm

\noindent{\em Proof of Lemma \ref{prop1}.}   We start with a preliminary  observation about sums of scalars.
Given  $u_1,\cdots,u_d\in \R$ and  $\sum_{i=1}^du_i=k$, with $k$ an integer, we assert there exists $v_1,\cdots,v_d$ such that $u_i\le v_i<u_i+1$ for every $1\le i\le d$ and $\sum_{i=1}^d\lfloor v_i\rfloor=k$.
We prove this assertion  for $k=0$ as the proof is similar for any other integer $k$. We can see that
$$\sum_{i=1}^d \lfloor u_i\rfloor\ge \sum_{i:u_i\mbox{ \tiny{integer}}} u_i+\sum_{i:u_i\mbox { \tiny{not integer}}}(u_i-1)=-\#\{i:u_i\mbox{ not integer}\}.$$
That number is at most equal to $-d$. 
Let us say that this number is equal to $-k$. Then it suffices to modify exactly $k$ of the $u_i$ which are not integers into $v_i$ with $u_i\le v_i<u_i+1$ so that the integer part will grow by $1$ exactly.

Given this assertion, we now prove Lemma \ref{prop1}.  Any point $(z_1,\cdots,z_d)$ of $[0,1)^d$ is coded by a $d$-tuple of integers $(k_1,\cdots,k_d)$ such that
 $\frac{1}{2}\delta 3^{k_i}\le z_i\le \frac{1}{2}\delta 3^{k_i + 1}$, i.e.
\begin{equation}\label{eq:defki}
k_i=\lfloor\log_3(\delta^{-1}2z_i)\rfloor.
\end{equation}

The point $(z_1,\cdots,z_d)$ belongs to an $M$-region in the collection ${\mathcal M}_K(\0,\delta)$ iff $k_1+\cdots+k_d=\log_3(T/\delta^d)$. Indeed, $(k_1+\cdots+k_d=\log_3(T/\delta^d))$ means that $(z_1,\cdots,z_d)\in M((3^{k_1}\delta,\cdots, 3^{k_d}\delta ))$, which is an $M$-region centered at a point on $K(v=T)$, where we recall that $\log_3(T/\delta^d) \in \Z$  by assumption.

Now, let $(z_1^{(0)},\cdots,z_d^{(0)})\in K(v=T)$. Let us prove that $M((z_1^{(0)},\cdots,z_d^{(0)}))$ intersects an $M$-region in ${\mathcal M}_K(\0,\delta)$. To do this, we rewrite the equation of $M((z_1^{(0)},\cdots,z_d^{(0)}))$ in terms of $(k_1,\cdots,k_d)$-coordinates.  We observe
\begin{align}\label{eq:macbeathequiv}
&(z_1,\cdots,z_d)\in M((z_1^{(0)},\cdots,z_d^{(0)}))\nonumber\\
&\Longleftrightarrow \frac{z_i^{(0)}}{2}\le z_i\le \frac{3}{2}z_i^{(0)}
\ \ \forall\,1\le i\le d\nonumber\\
&\Longleftrightarrow \log_3(\delta^{-1}z_i^{(0)})\le     \log_3(\delta^{-1}2z_i) < \log_3(\delta^{-1}z_i^{(0)})+1 \ \ \forall\,1\le i\le d.
  \end{align}
  Since $(z_1^{(0)},\cdots,z_d^{(0)})\in K(v=T)$, this implies $\Pi_{i = 1}^d z_i^{(0)} =T$.
We  appeal to our assertion, setting  $u_i:=\log_3(\delta^{-1}z_i^{(0)})$ and $k:=\log_3(T/\delta^d)$.  We  choose $v_1,\cdots,v_d$ such that $\sum_{i=1}^d\lfloor v_i\rfloor=k=\log_3(T/\delta^d)$ and $u_i\le v_i< u_i+1, 1 \leq i \leq d.$  Then we take $\tilde{z}_i$ such that $v_i=\log_3(\delta^{-1}2\tilde{z}_i)$
and we put $\tilde{k}_i:= \lfloor\log_3(\delta^{-1}2\tilde{z}_i)\rfloor.$
Now $(\tilde{z}_1,\cdots,\tilde{z}_d)\in M((z_1^{(0)},\cdots,z_d^{(0)}))$ because $(\tilde{z}_1,\cdots,\tilde{z}_d)$ satisfies \eqref{eq:macbeathequiv}.
On the other hand, $(\tilde{z}_1,\cdots,\tilde{z}_d)$ belongs to an $M$-region in ${\mathcal M}_K(\0, \delta)$
 because $\sum_{i=1}^d \tilde{k}_i=\sum_{i=1}^d\lfloor v_i\rfloor=\log_3(T/\delta^d)$. Thus $M((z_1^{(0)},\cdots,z_d^{(0)}))$ intersects an $M$-region in ${\mathcal M}_K(\0,\delta)$,
 showing the desired maximality of ${\mathcal M}_K(\0,\delta)$. \qed

\vskip.3cm


\noindent{\em Proof of Lemma \ref{Lem3.2}}. Similar to Section 2 of \cite{BR2}, we use the notation  $C(z^{(0)})=  K \cap H^+(z^{(0)})$
to signify a cap of $K$ at $z^{(0)}$,
where $H^+(z^{(0)})$  is at   \eqref{eq:eqhypertangent}.  
Also, if $(3^{k_1}\delta,\cdots,3^{k_d}\delta)$ denotes the center of the $M$-region $M_j$, then
for $\gamma > 0$ we
define $K_j^{\gamma}:=C^{6\gamma}((3^{k_1}\delta,\cdots,3^{k_d}\delta))$, where  $C^{6\gamma}(z^{(0)})$ is
the enlarged  cap 
\begin{equation}
  \label{eq:equationcap}
C^{6\gamma}(z^{(0)}):= \{ (z_1,...,z_d): \  \sum_{i=1}^d \frac{z_i}{z_i^{(0)}}\le 6d\gamma\}.
  \end{equation}

 Looking closely at Section 5 of \cite{BR2}, it suffices to show the set inclusions\\
(i) $S_j' \subset K_j^{\gamma}$ where  $\gamma \in (\frac{1}{6d}(2^{d-1}d6^d+\frac{3}{2}(d-1)), \infty)$  and\\(ii)   $K'_j  \subset S_j$, where
 $K'_j:= M_{K}((3^{k_1}\delta,\cdots,3^{k_d}\delta))\cap C((3^{k_1}\delta,\cdots,3^{k_d}\delta)).$

 When $S'_j$ is either a cone set or cone-cylinder set, the cardinality of which is bounded independently of $\la$,
 these set inclusions are satisfied for large $\gamma$.  The only challenge is to show these inclusions for
 the cylinder sets $S'_j$.

We start by showing the first inclusion. Let $z^{(0)}\in K(v=T)$. 
The aim is to show that the explicit regions $S_j'$
defined at \eqref{Sprimej} satisfy the requirement from \cite{BR2}, i.e. that there exists an explicit $\gamma$ depending only on dimension $d$ such that $S_j'\subset K_j^{\gamma}$ for every $j$. Actually, our explicit value of $\gamma$ will be larger than the one used in \cite{BR2} (see display before (5.4) therein) but we claim that this does not affect any of the results from \cite{BR2} and in particular it does not modify the construction of the dependency graph.



We now describe the set $S_j'$ constructed from the $M$-region containing $z^{(0)}$.
Recalling \eqref{Sprimej}, assume $z_d^{(0)}=\min_{1\le i\le d}z_i^{(0)}$ and that there is no tie for sake of simplicity (the case of a tie would be treated analogously).  We have
$$S_j' = K(v\le T^*)\cap [\frac{1}{2}z_1^{(0)},\frac{3}{2}z_1^{(0)}]\times\cdots\times [\frac{1}{2}z_{d-1}^{(0)},\frac{3}{2}z_{d-1}^{(0)}]\times \R.$$
Putting $\tilde{T^*}:={d!}/{d^d}T^*$,  the height of $S_j'$  above the point $(z_1,\cdots,z_{d-1},0)$  is
\begin{equation}
  \label{eq:hauteurSj}
z_d=\frac{\tilde{T^*}}{\Pi_{i=1}^{d-1}z_i}.
\end{equation}
In particular, \eqref{eq:equationcap} and \eqref{eq:hauteurSj} imply that $S_j'\subset K_j^{\gamma}$ as soon as for every $(z_1,\cdots,z_{d-1})\in [\frac{1}{2}z_1^{(0)},\frac{3}{2}z_1^{(0)}]\times\cdots\times [\frac{1}{2}z_{d-1}^{(0)},\frac{3}{2}z_{d-1}^{(0)}]$, we have
\begin{equation}\label{eq:diffhauteur}
\frac{\tilde{T^*}}{\Pi_{i=1}^{d-1}z_i}\le z_d^{(0)}(6d\gamma-\sum_{i=1}^{d-1}\frac{z_i}{z_i^{(0)}}).
\end{equation}
Noticing that on the one hand,
$$\frac{\tilde{T^*}}{\Pi_{i=1}^{d-1}z_i}\le \frac{2^{d-1}d6^d d!T}{d^d\Pi_{i=1}^{d-1}z_i^{(0)}}=2^{d-1}d6^dz_d^{(0)}$$
and that on the other hand,
$$z_d^{(0)}(6d\gamma-\sum_{i=1}^{d-1}\frac{z_i}{z_i^{(0)}})\ge z_d^{(0)}(6d\gamma-\frac{3}{2}(d-1)),$$
we conclude that \eqref{eq:diffhauteur} is satisfied as soon as $\gamma \in [\frac{1}{6d}(2^{d-1}d6^d+\frac{3}{2}(d-1)), \infty)$.

We show now the second inclusion (ii). 
In particular, because of its definition, $S_j$ contains $M_{K'}((3^{k_1}\delta,\cdots,3^{k_d}\delta))\cap C((3^{k_1}\delta,\cdots,3^{k_d}\delta))=K'_j$. 
This concludes the proof of Lemma \ref{Lem3.2}.   \qed

\vskip.5cm

\noindent{\em Proof of Lemma \ref{L4}.}
   We first assert that
\be \label{volbd}
\Vol( {\cal A}(s, T^*,K, \delta_1) ) = O( \log \log \la (\log \la)^{d -2 + 1/d}  \la^{-1} ).
\ee
Indeed, we notice that ${\cal A}(s, T^*,K, \delta_1)\subset K(v \leq T^*).$   Next we apply the bound in display (4.1)  of
\cite{BB} with the $\varepsilon$ and $\varphi$ of that bound set to $T^*$ and a constant multiple of $\delta_1$, respectively. This is possible
because without loss of generality the  parallelepiped $p_d(v_i, \delta_1):= a_i^{-1}([0, \delta_1]^d )$ contains the intersection of $K$ with a slab of thickness proportional to $\delta_1$.  Note that
$\varphi^d \geq \text{const} \cdot \varepsilon,$ which gives \eqref{volbd}.

We prove part (a) of Lemma \ref{L4} for the $k$-face functional $\xi_k$ and then treat  the
volume functional $\xi_V$.   We first show $\var[Z_0(\delta_1){{\bf 1}(A_\la)}] = o(\Var Z)$, which goes as follows.
Put $\tx(x, \P_\la) := \xi_k(x, \P_\la){{\bf 1}(A_\la)}$.
Letting $\P^y_\la := \P_\la \cup \{y\}$ we have  
\be \label{July21}
\Var \sum_{x \in \P_\la(s, T^*,K, \delta_1) } \tx(x, \P_\la)
= V_1 + V_2,
\ee
where
$$
V_1 := \E \sum_{x \in \P_\la(s, T^*,K, \delta_1) } \tx(x, \P_\la)^2
$$
and
$$
V_2:= \E \sum_{x,y \in \P_\la(s, T^*,K, \delta_1); \ x \neq y} [\tx(x, \P^y_\la)\tx(y, \P^x_\la)  - \E \tx(x, \P_\la)\E \tx(y, \P_\la)].$$
We  bound $V_1$ as follows.
Each  $x \in {\cal A}(s, T^*, K, \delta_1)$ belongs to some  $S'_i$ region in the collection
$\{S_j'\}_{j=1}^{m(T, \delta_1)}$. 
Let ${\cal S}_x$ denote the union of those $S'_j$ such that there is an edge between $i$ and $j$.
By Theorem 6.2 of \cite{BR2}, we have ${\rm{card}} ({\cal S}_x) \leq D(\la)$, where $D(\la)= O(\log \log \la^{6(d-1)})$ is the maximal degree of the
dependency graph $({\cal V}_{\cal G} , {\cal E}_{\cal G} )$, where ${\cal V}_{\cal G} := {\cal S}'(\delta_1).$
On $A_\la$ we have $\max_{j \leq m(T,\delta_1)} \text{card} (S_j \cap \P_\la) \leq c(d) \log \log \la$,
as explained  two lines after display (5.4) in \cite{BR2}.
It follows that  on $A_\la$ at most $O((\log \log \la)^{6(d-1) + 1})$ points in  $\P_\la$ can
potentially contribute to a $k$-face containing $x \in \P_\la(s, T^*, K, \delta_1)$.
By McMullen's bound \cite{Mc}, the number of $k$-faces on an $n$ point set is
bounded by $Cn^{d/2}$. The score at $x$ thus satisfies
\be \label{JeY1}
\sup_{x, y\in {\cal A}(s, T^*,K, \delta_1)    } | \tx(x, \P^y_\la) | = O((\log \log \la)^{(6(d-1) + 1)d/2}).
\ee

Combining \eqref{volbd}, \eqref{JeY1}, and using the Slivnyak - Mecke formula, we find the desired bound for the first term
in \eqref{July21}:
\begin{align*}
V_1 & = \la \int_{ {\cal A}(s, T^*,K, \delta_1) } \E [\tx(x, \P_\la)^2] dx\\
& = O( \la \Vol( {\cal A}(s, T^*,K, \delta_1) ) (\log \log \la)^{(6(d-1) + 1)d} )\\
& =  o((\log
\la)^{d-1}) \\
& = o(\Var [Z]).
\end{align*}

Now we bound $V_2$.
We treat separately  the sum over $x\in \P_\la(s,T^*,K, \delta_1)$ and $y\not\in S_x$ and the sum over $x\in \P_\la(s,T^*,K, \delta_1)$ and $y\in S_x$. When $x\in {\cal A}(s,T^*,K)$ and $y\not\in S_x$, we have  $$\E[\xi_k(x, \P^y_\la)\xi_k(y, \P^x_\la)|A_\la]  - \E[\xi_k(x, \P_\la)|A_\la]\E[\xi_k(y, \P_\la)|A_\la]=0.$$
Consequently,
\begin{align}\label{covzero}
&\E[\tx(x, \P^y_\la)\tx(y, \P^x_\la)]-\E[\tx(x, \P_\la)]\E[\tx(y, \P_\la)]\nonumber\\&=\E[\xi_k(x, \P^y_\la)\xi(y, \P^x_\la)|A_\la]P[A_\la]-\E[\xi_k(x, \P_\la)|A_\la]\E[\xi_k(y, \P_\la)|A_\la]P[A_\la]^2\nonumber\\
&=\E[\xi_k(x, \P_\la)|A_\la]\E[\xi_k(y, \P_\la)|A_\la]P[A_\la]P[A_\la^c].
\end{align}
Combining \eqref{defA}, \eqref{JeY1} and  \eqref{covzero} and applying the Slivnyak-Mecke formula, we get
$$\E[\sum_{x,y \in \P_\la(s, T^*,K, \delta_1); \ y \notin {\cal S}_x } [\xi_k(x, \P^y_\la)\xi_k(y, \P^x_\la)  - \E \xi_k(x, \P_\la)\E \xi_k(y, \P_\la)] = o(\Var Z).$$
Now we prove that
\be \label{July21a}
\E \sum_{x,y \in \P_\la(s, T^*,K, \delta_1); \ y \in {\cal S}_x } [ \tx(x, \P^y_\la)\tx(y, \P^x_\la)  - \E \tx(x, \P_\la)\E \tx(y, \P_\la)] = o(\Var Z).\ee
By \eqref{JeY1} we also have
\be \label{JY1}
\sup_{x, y} | \tx(x, \P^y_\la)\tx(y, \P^x_\la)  - \E \tx(x, \P_\la)\E \tx(y, \P_\la) | = O((\log \log \la)^{(6(d-1) + 1)d}).
\ee
Moreover, we deduce from (5.4) in \cite{BR2} that
\begin{align}\label{JY2}
\sup_{x\in {\cal A}(s,T^*,K)}\Vol({\mathcal S}_x)&
\le \sup_{x\in {\cal A}(s,T^*,K)}\mbox{card}({\mathcal S}_x)\cdot\sup_{S'_j\in {\cal S}'(\delta_1)}\Vol(S'_j)= O\left(\frac{\log\log\la^{6(d-1)+1}}{\la}\right).
\end{align}
Consequently, using  the Slivnyak - Mecke formula, \eqref{volbd}, \eqref{JY1} and \eqref{JY2}, we get
\begin{align*} \label{July21aa}
& \E \sum_{x,y \in \P_\la(s, T^*,K, \delta_1); \ y \in {\cal S}_x } [ \tx(x, \P^y_\la)\tx(y, \P^x_\la)  - \E \tx(x, \P_\la)\E \tx(y, \P_\la)]\\
& = O( \la^2 \Vol(\A(s, T^*, K, \delta_1)) \sup_{x\in {\mathcal A}(s,T^*,K)}\Vol({\mathcal S}_x)(\log \log \la)^{(6(d-1) + 1)d})\\
&= o(\Var Z).
\end{align*}

Now we show $\var[Z_0(\delta_1) | A_\la] = o(\Var Z)$. Notice that
\begin{align*}
\E[Z_0^2(\delta_1) {\bf 1}(A_\la) ]   & \leq  \la^2 \int_{ {\cal A}(s, T^*,K, \delta_1) } \int_{ {\cal A}(s, T^*,K, \delta_1) } \E [\tx(x, \P^y_\la) \tx(y, \P^x_\la)] dydx\\
& = O( \la^2 (\Vol( {\cal A}(s, T^*,K, \delta_1) ))^2 (\log \log \la)^{(6(d-1) + 1)d} )\\
& =  O( (\log \la)^{2(d-2) + 2/d}  (\log \log \la)^{(6(d-1) + 1)d + 1}),
\end{align*}
where we use \eqref{volbd}.  Thus by \eqref{defA} we get $\E[Z_0^2(\delta_1) {\bf 1}(A_\la) ]P[A_\la^c] = o(1).$
The desired bound $\var[Z_0(\delta_1) | A_\la] = o(\Var[Z])$ follows from this estimate and the identity
$\var[Z_0(\delta_1) | A_\la] = P[A_\la]^{-2} (\var[Z_0(\delta_1){{\bf 1}(A_\la)}] - \E [Z_0^2(\delta_1){{\bf 1}(A_\la)}](1 - P[A_\la]) ).$

Now we show $\var[Z_0(\delta_1){{\bf 1}(A_\la)}] = o(\Var[Z])$ when $\xi$ is the volume score.
Recall from \eqref{volumescore2} that  ${\cal F}_{d-1}(x)$ is  the collection of facets in
$K_\la$ which contain $x$.
Regardless of whether we use  \eqref{volumescore1} and \eqref{volumescore2}, we have on $A_\la$ that
\begin{align*}
\xi_V(x, \P_\la )&\le \lambda\mbox{card}(\mbox{${\cal F}_{d-1}(x)$})\\&\hspace*{.5cm}\times \sup \{ \Vol\mbox{({Cap} of $K$): \ Cap of $K$  tangent to $K(v=s')$ with $s\le s'\le T^*$} \}.
\end{align*}
 Indeed, for any face $F$ and facet $F'$ containing $F$, let $H$ be the hyperplane containing $F'$. Then $C_F(K_\la)\cap K$ is
included in the cap of $K$ bounded by $H$. Since $H$ meets ${\mathcal A}(s,T^*,K)$ but not $K(v\ge T^*)$, it is tangent to some
 some $K(v=s')$ with $s\le s'\le T^*$. Consequently, Lemma 2.4 in \cite{BR2} yields
\be \label{volcap}
\sup \{ \Vol\mbox{({Cap} of $K$):\ Cap of $K$ tangent to $K(v=s')$ with $s\le s'\le T^*$} \}=O(\frac{\log\log\la}{\la}).\ee
 Moreover, \eqref{JeY1} implies that ${\rm{card}}({\cal F}_{d-1}(x)) = O((\log\log\la)^{(6(d-1)+1)d/2})$.  Thus on $A_\la$ we have
 $$\sup_{x \in Q_0^c \cap {\cal A}(s, T^*, K)} |\xi_V(x,\P_\la)| = O\left((\log\log\la)^{(6(d-1)+1)\frac{d}{2}+1}\right).$$
Using  this bound in place of the bound \eqref{JeY1} and following the method for the $k$-face functional verbatim, we obtain $\var[Z_0(\delta_1){{\bf 1}(A_\la)}] = o((\log
\la)^{(d-1)}) = o(\Var [Z])$.
Following the discussion for the $k$-face functional, we also have $\E[Z_0^2(\delta_1) {\bf 1}(A_\la) ] = o(1)$.
This gives $\var[Z_0(\delta_1) | A_\la] = o(\Var[Z])$ as explained  for the $k$-face score.

To show part (b) of Lemma \ref{L4}, we recall Figure \ref{thelastpictureshow}  and notice that
  \begin{align*}
&0\le \Vol(K\setminus K_\la)-\frac{1}{\lambda}\sum_{x\in\P_{\la}}\xi_V(x,\P_{\la})\\&\le \mbox{card}(
\{\mbox{facets intersecting one of the $p_d({\mathscr V}_i, \delta_0)$ and its complement}\})\\
&\hspace*{.5cm}\times \sup \{ \Vol\mbox{({Cap} of $K$): \ Cap of $K$  tangent to $K(v=s')$ with $s\le s'\le T^*$} \}\\
&\le c\sum_{x\in \P_\la(s,T^*,K,\delta_1)}\frac{\log\log\la}{\la},
    \end{align*}
where the last inequality uses the estimate \eqref{volcap}. We now apply the exact same method as for the proof of part (a) in order to bound the variance of the difference $\Vol(K\setminus K_\la)-\frac{1}{\lambda}\sum_{x\in\P_{\la}}\xi_V(x,\P_{\la})$ and show the required statement (b). This concludes the proof of Lemma \ref{L4}.
 \qed

\vskip.5cm

\noindent{\em Proof of Lemma \ref{Lem9}}. Recall the definition of $Z$ at \eqref{Zsum}. For  $1 \leq i \leq |{\cal V}_K|$
 the assertions
$$ \max\{ |\E[Z]-\E[Z | A_\la ]| , \ |\E[Z_i]-\E[Z_i | A_\la ]| \} =o(\E[Z])
$$
and
$$ \max\{ |\Var[Z]-\Var[Z | A_\la ]| , \ |\Var[Z_i]-\Var[Z_i | A_\la ]| \} =o(\Var[Z])
$$
follow from Lemmas 8.2 and 8.3 of \cite{BR2}. We now show
\be \label{assert1}  |\Var[Z]-\Var[Z{\bf 1}({A_\la})]|=o(\Var[Z]); \ \   |\Var[Z_i]-\Var[Z_i{\bf 1}({A_\la})]|=o(\Var[Z]).
\ee
We only prove the first assertion, as the second follows from identical methods.  We prove the first assertion
when $Z$ is the number of $k$-dimensional faces of $K_\la$  and then treat the case when $Z$  is the defect volume of $K_\la$.
We have
\begin{align*}
  \Var[Z]&=\Var[Z{\bf 1}(A_\la))+Z{\bf 1}(A_\la^c)]\\
&=\Var[Z{\bf 1}({A_\la})]+\Var[Z{\bf 1}({A_{\la}^c})] + 2\Cov(Z{\bf 1}({A_\la}),Z{\bf 1}({A_\la^c})).
\end{align*}
Consequently,
\begin{equation}
  \label{eq:diffvar}
|\Var[Z]-\Var[Z{\bf 1}({A_\la})]| \le \E[Z^2{\bf 1}({A_{\la}^c})]+2\sqrt{\Var[Z{\bf 1}({A_\la})]}\sqrt{\E[Z^2{\bf 1}({A_{\la}^c})]}.
\end{equation}
We first estimate $\E[Z^2{\bf 1}({A_{\la}^c})]$ as follows. Using the event $B_{\la}$ provided by \cite{BR2} (and denoted by $B$ there, see p. 1519 of
\cite{BR2}), we write
\begin{equation}
 \label{eq:varneglig}
\E[Z^2{\bf 1}({A_{\la}^c})]  =  \E [Z^2{\bf 1}({A_\la^c\cap B_\la})]+ \E[Z^2{\bf 1}({A_\la^c\cap B_\la^c})].
\end{equation}
We treat separately each term  on the right hand side of \eqref{eq:varneglig}. Let us start with the first one:  On the event $B_{\la}$, we know that only the points inside $K(v\le d6^d \la^{-1} \log\la )$ are needed to construct $K_\la$  and that their cardinality is $O((\log\la)^d)$. Consequently,  by McMullen's bound, we have $Z = O((\log \la)^{d^2/2})$.  It follows from \eqref{defA} that
\begin{equation}
  \label{eq:term1}
  \E[Z^2{\bf 1}({A_\la^c\cap B_\la})]=O((\log \la)^{d^2}P[A_\la^c])=O((\log \la)^{-3d^2}).
\end{equation}
To estimate   the second term of \eqref{eq:varneglig}, we proceed as in the proof of Lemma 8.2 in \cite{BR2}:
\begin{align}\label{eq:decomp}
\E[Z^2{\bf 1}({A_\la^c\cap B_\la^c})]
&=\sum_{m=0}^{\infty}\E[Z^2{\bf 1}({A_\la^c\cap B_\la^c})| \ \mbox{card}(\P_\la)=m]P[\mbox{card}(\P_\la)=m] \nonumber\\
&=\sum_{m=0}^{\lfloor 3\Vol(K)\la \rfloor }\E[Z^2{\bf 1}({A_\la^c\cap B_\la^c})|\ \mbox{card}(\P_\la)=m]P[\mbox{card}(\P_\la)=m]\nonumber
\\&\hspace*{1cm}+\sum_{m= \lfloor 3\Vol(K)\la \rfloor + 1}^{\infty}\E[Z^2{\bf 1}({A_\la^c\cap B_\la^c})| \ \mbox{card}(\P_\la)=m]P[\mbox{card}(\P_\la)=m].
\end{align}
When $\mbox{card}(\P_\la)=m$, we have $Z = O(m^{d/2})$. In particular, when $m \in \{0, 1,...., \lfloor 3\Vol(K)\la) \rfloor \}$, we bound $Z^2$ by $O(\la^d)$. Consequently, the identity  \eqref{eq:decomp} gives
\begin{align*}
\E[Z^2{\bf 1}({A_\la^c\cap B_\la^c})] & = O \left( \la^d P[{A_\la^c\cap B_\la^c}]+ \sum_{m =  \lfloor 3\Vol(K)\la \rfloor + 1}m^{d}P[\mbox{card}(\P_\la)=m] \right)\\
& = O \left( \la^d P[B_\la^c|A_\la^c]P[A_\la^c]+\E[(\mbox{card}(\P_\la))^d{\bf 1}(\mbox{card}(\P_\la)\ge 3\Vol(K)\la) ] \right ).
\end{align*}
In view of the bound (8.2) in \cite{BR2}, the first term is $O(\la^{-2d+1}(\log\la)^{-4d^2})$ whereas  the second one is decreasing exponentially fast in $\la$. Consequently, we have
\begin{equation}
  \label{eq:term2}
  \E[Z^2{\bf 1}({A_\la^c\cap B_\la^c})]=o(\la^{-2d+1}).
\end{equation}
Inserting \eqref{eq:term1} and \eqref{eq:term2} into \eqref{eq:varneglig}, we get
$\E[Z^2{\bf 1}({A_{\la}^c})]=o((\log \la)^{-3d^2})$. This fact combined with \eqref{eq:diffvar} and Theorem 1.3 from \cite{BR} implies \eqref{assert1}. 

We now adapt the above proof when $Z$ is the defect volume. Regarding the first term of \eqref{eq:varneglig},  Theorem 2.7 in \cite{BR2} implies that on $B_\la$, we  have $Z\le \Vol(K(v\le d6^d \log\la/\la))=O((\log\la/\la) (\log \la)^{d-1}))$. Consequently, \eqref{eq:term1} is replaced by
 \begin{equation}
   \label{eq:term1_2}
   \E[Z^2{\bf 1}({A_\la^c\cap B_\la})]= O\left(
\frac{(\log \la)^{2d}}{\la^2}
P[A_\la^c] \right)
=O\left(\frac{(\log \la)^{-4d^2+2d}}{\la^2}\right).
 \end{equation}
For the second term of \eqref{eq:varneglig}, we simply bound $Z^2$ by a constant  and we get from \eqref{eq:decomp} that
 \eqref{eq:term2} also holds for the defect volume. Now, inserting \eqref{eq:term1_2} and \eqref{eq:term2} into \eqref{eq:varneglig}, we find  $
\E[Z^2{\bf 1}({A_{\la}^c})] = o((\log \la)^{-4d^2+2d}\la^{-2})$ and so \eqref{assert1}  also holds when $Z$ is  the defect volume.
This completes the proof of Lemma \ref{Lem9}.   \qed

\vskip.5cm

\noindent{\em Acknowledgements}.   P. Calka is  grateful to the Department of Mathematics at Lehigh University for its support and great help during his several stays.  J. Yukich likewise  thanks the Department of Mathematics at the Universit\'e de Rouen for its kind hospitality and support.

Pierre Calka, Laboratoire de Math\'ematiques Rapha\"el Salem, Universit\'e de Rouen,
 Avenue de l'Universit\'e, BP.12, Technop\^ole du Madrillet, F76801 Saint-Etienne-du-Rouvray France;
\ \ \  {\texttt pierre.calka@univ-rouen.fr}

\vskip.5cm

J. E. Yukich, Department of Mathematics, Lehigh University,
Bethlehem PA 18015;\\
 \ \ {\texttt joseph.yukich@lehigh.edu}

\end{document}